\documentclass[a4paper,12pt,twoside]{article}
\usepackage{graphicx}\usepackage{a4wide}
\usepackage{amssymb,amsmath,amsthm}
\usepackage{siunitx}
\usepackage{subfigure,color,colordvi,graphicx,xcolor}
\usepackage[only,llbracket,rrbracket,interleave]{stmaryrd}
\usepackage[superscript,sort,compress]{cite}

\makeatletter
\newcommand*\wt[2][0.2ex]{%
        \begingroup
        \mathchoice{\wt@helper{#1}{#2}{\displaystyle}{\textfont}}
                   {\wt@helper{#1}{#2}{\textstyle}{\textfont}}
                   {\wt@helper{#1}{#2}{\scriptstyle}{\scriptfont}}
                   {\wt@helper{#1}{#2}{\scriptscriptstyle}{\scriptscriptfont}}%
        \endgroup
        #2%
}
\newcommand*\wt@helper[4]{%
        \def\currentfont{\the#41}%
        \def\currentskewchar{\char\the\skewchar\currentfont}%
        \setbox\tw@\hbox{\currentfont#2\currentskewchar}%
        \dimen@ii\wd\tw@
        \setbox\tw@\hbox{\currentfont#2{}\currentskewchar}%
        \advance\dimen@ii-\wd\tw@
        \rlap{\raisebox{-#1}{$\m@th#3\kern\dimen@ii\widetilde{\phantom{#2}}$}}%
}
\makeatother

\newenvironment{keywords}{\begin{quote}\emph{\textbf{Keywords:}}}{\end{quote}}
\newtheorem{lemma}{Lemma}

\newcommand{\bm}[1]{\text{\boldmath $#1$\unboldmath}}
\newcommand{\abs}[1]{\lvert#1\rvert}

\newcommand{\vect}[1]{\mathbf{#1}}
\newcommand{\mat}[1]{\mathbf{#1}}

\newcommand{\grad}{\bm{\nabla}}

\newcommand{\eltwo}{\ensuremath{\mathcal{L}_2}}

\newcommand{\nen}  {\ensuremath{\texttt{n}_{\texttt{en}}}}
\newcommand{\nfn}  {\ensuremath{\texttt{n}_{\texttt{fn}}}}

\newcommand{\nsd}  {\ensuremath{\texttt{n}_{\texttt{sd}}}}
\newcommand{\numel}{\ensuremath{\texttt{n}_{\texttt{el}}}}
\newcommand{\numfa}{\ensuremath{\texttt{n}_{\texttt{fa}}}}

\newcommand{\nit}{\ensuremath{\texttt{n}_{\texttt{i}}}}

\newcommand{\hu}{\hat{u}}
\newcommand{\bu}{\bm{u}}

\newcommand{\bhu}{\widehat{\bu}}
\newcommand{\bq}{\bm{q}}

\newcommand{\bn}{\bm{n}}

\newcommand{\bL}{\bm{L}}

\newcommand{\qe}{\vect{q}_e}

\newcommand{\uHi}{\hat{\text{u}}_i}
\newcommand{\uHj}{\hat{\text{u}}_j}

\newcommand{\Le}{\mat{L}_e}
\newcommand{\bue}{\vect{u}_e}
\newcommand{\pe}{\text{p}_e}

\newcommand{\buHj}{\hat{\vect{u}}_j}

\newcommand{\Vone}{\mathcal{V}^1}
\newcommand{\Vzero}{\mathcal{V}^0}
\newcommand{\VzeroHat}{\mathcal{\hat{V}}^0}

\newcommand{\VoneDomain}{\Vone(\Omega)}
\newcommand{\VzeroDomain}{\Vzero(\Omega)}
\newcommand{\VzeroHatSkeleton}{\VzeroHat(\Gamma)}

\newcommand{\VoneElem}{\Vone(\Omega_e)}
\newcommand{\VzeroElem}{\Vzero(\Omega_e)}
\newcommand{\VzeroFaces}{\VzeroHat(\Gamma\cup\Gamma_N)}

\newcommand{\Insd}{\mat{I}_{\nsd}}

\newcommand{\jump}[1]{\llbracket #1\rrbracket}

\newcommand{\Pu}{\mathbb{P}_0}

\newcommand{\bw}{\bm{w}}

\newcommand{\Aset}{\mathcal{A}_e}
\newcommand{\Bset}{\mathcal{B}_e}
\newcommand{\Dset}{\mathcal{D}_e}
\newcommand{\Nset}{\mathcal{N}_e}

\newcommand{\projc}{p_{e,i}}
\newcommand{\projv}{\vect{p}_{e,i}}
\newcommand{\projm}{\mat{P}_{e,i}}

\newtheorem{remark}{Remark}

\RequirePackage{algorithm}[0.1] 
\usepackage{algorithmic} 

\graphicspath{{figsPDF/}{figsPNG/}}

\begin{document}
\title{A second-order face-centred finite volume method for elliptic problems}

\author{
\renewcommand{\thefootnote}{\arabic{footnote}}
			  L.M. Vieira\footnotemark[1]\textsuperscript{ \ ,}\footnotemark[2] \ ,
			  M. Giacomini\footnotemark[2], \
			  R. Sevilla\footnotemark[1]\textsuperscript{ \ ,}*, \
             A. Huerta\footnotemark[2]
}

\date{\today}
\maketitle

\renewcommand{\thefootnote}{\arabic{footnote}}

\footnotetext[1]{Zienkiewicz Centre for Computational Engineering, College of Engineering, Swansea University, Wales, UK}
\footnotetext[2]{Laboratori de C\`alcul Num\`eric (LaC\`aN), ETS de Ingenieros de Caminos, Canales y Puertos, Universitat Polit\`ecnica de Catalunya, Barcelona, Spain
\vspace{5pt}\\
* Corresponding author: Ruben Sevilla. \textit{E-mail:} \texttt{r.sevilla@swansea.ac.uk}
}

\begin{abstract}
A second-order face-centred finite volume method (FCFV) is proposed. Contrary to the more popular cell-centred and vertex-centred finite volume (FV) techniques, the proposed method defines the solution on the faces of the mesh (edges in two dimensions). The method is based on a mixed formulation and therefore considers the solution and its gradient as independent unknowns. They are computed solving an element-by-element problem after the solution at the faces is determined. The proposed approach avoids the need of reconstructing the solution gradient, as required by cell-centred and vertex-centred FV methods. This strategy leads to a method that is insensitive to mesh distortion and stretching. The current method is second-order and requires the solution of a global system of equations of identical size and identical number of non-zero elements when compared to the recently proposed first-order FCFV. The formulation is presented for Poisson and Stokes problems. Numerical examples are used to illustrate the approximation properties of the method as well as to demonstrate its potential in three dimensional problems with complex geometries. The integration of a mesh adaptive procedure in the FCFV solution algorithm is also presented. 
\end{abstract}

\begin{keywords}
finite volume method, face-centred, second-order convergence, hybridisable discontinuous Galerkin
\end{keywords}

\section{Introduction}
	
Cell-centred and vertex-centred second-order finite volume (FV) methods are still the predominant techniques used in commercial and industrial computational fluid dynamics (CFD) solvers due to their robustness, easy implementation and relatively low cost~\cite{MR1925043,MR2731357,MR2417929,BO2004,EYMARD2000713,droniou2014review}. Both, cell-centred and vertex-centred, FV techniques require a reconstruction of the gradient of the solution to ensure second-order convergence of the unknown and first-order convergence of the fluxes~\cite{diskin2010comparison,diskin2011comparison,chassaing2013moving,aguirre2015upwind}. The accuracy of the scheme is therefore dependent on the accuracy of the reconstruction technique, which in turns depends on the quality of the mesh. In particular, FV methods are known to suffer an important loss of accuracy, and sometimes even a loss of second-order convergence, when unstructured meshes are used or highly stretched and deformed cells are present in the computational mesh~\cite{diskin2010comparison,diskin2011comparison}.

In~\cite{FCFV2018}, the authors proposed a novel methodology, called face-centred finite volume (FCFV) method. The technique can be seen as a particular case of the hybridisable discontinuous Galerkin (HDG) method by Cockburn and co-workers \cite{cockburn2004characterization,Jay-CG:05,Jay-CG:05-GAMM,Jay-CGL:09}, by considering an element-by-element constant degree of approximation. Contrary to other FV methods, the resulting FCFV method is insensitive to mesh distortion and stretching. In addition, being a mixed method, it provides first-order convergence of the gradient of the solution without the need of a reconstruction. Contrary to other mixed methods, in the context of incompressible flows, the FCFV method passes the so-called Ladyzhenskaya-Babu{\v s}ka-Brezzi
(LBB) condition. However, the FCFV method only provides first-order convergence for the solution.

This paper proposes a second-order FCFV with a computational cost almost identical to the original, first-order, FCFV. The method is also insensitive to mesh distortion and stretching and provides a first-order approximation of the gradient of the solution without the need of a reconstruction. Similar to the first-order FCFV method, the proposed method can be seen as a particular case of the HDG method, in which the space of approximation used for the solution is linear whereas constant degree approximation spaces are used for the gradient of the solution and the solution on the cell faces~\cite{oikawa2015hybridized,oikawa2016analysis,qiu2016superconvergent}. Therefore, the FCFV method inherits the convergence properties of HDG and it passes the LBB condition in the context of incompressible flows. It is worth noting that the degrees of freedom of the solution and its gradient inside each cell can be eliminated via a Schur complement procedure. The only global coupled degrees of freedom correspond to the value of the solution on the cell faces. For both the original and the second-order FCFV methods, the solution on the cell faces is approximated with a constant degree approximation. Hence both methods have a global matrix with the same size and same number of non-zero elements. Due to the extra accuracy of the proposed method, compared to the original FCFV~\cite{FCFV2018,FCFVelas}, this work also proposes a combination of first and second-order methods to produce an error indicator that is used to drive an $h$-adaptivity process. 

The remaining of the paper is organised as follows. In Section~\ref{sc:Poisson} the proposed second-order FCFV method is presented for the solution of the Poisson equation. Its extension to the Stokes problem is presented in Section~\ref{sc:Stokes}. The ability to combine first-order and second-order FCFV methods to perform an automatic mesh adaptive process is discussed in Section~\ref{sc:adaptivity}.
Section~\ref{sc:studies} presents a number of numerical experiments to validate the optimal approximation properties of the method and to compare the accuracy of the first-order and second-order FCFV methods in terms of the computational cost. The insensitivity to mesh distortion and stretching is also demonstrated using numerical experiments. Section~\ref{sc:examples} presents more challenging problems to show the potential of the proposed second-order FCFV method and its application in an automatic mesh adaptive process. Finally, Section~\ref{sc:Conclusion} summarises the conclusions of the work that has been presented.

\section{Second-order FCFV for the Poisson equation}
\label{sc:Poisson}

\subsection{Problem statement}
\label{sc:PoissonStatement}

Let us consider an open bounded domain $\Omega\in\mathbb{R}^{\nsd}$  with boundary $\partial\Omega=\Gamma_D\cup\Gamma_N$, $\Gamma_D \cap \Gamma_N = \emptyset$ and $\nsd$ the number of spatial dimensions. The strong form of the Poisson problem is
\begin{equation} \label{eq:Poisson}
\left\{\begin{aligned}
-\grad\cdot\grad u &= s       &&\text{in $\Omega$,}\\
u &= u_D  &&\text{on $\Gamma_D$,}\\
\bn\cdot\grad u &= t        &&\text{on $\Gamma_N$,}\\
\end{aligned}\right.
\end{equation}
where $s$ denotes a source term, $\bn$ is the outward unit normal vector to the boundary $\partial\Omega$ and $u_D$ and $t$ are the Dirichlet and Neumann data. 
	
The domain $\Omega$ is assumed to be partitioned in $\numel$ disjoint triangular or tetrahedral cells $\Omega_e$ in two and three dimensions respectively
\begin{equation}
\Omega =  \bigcup_{e=1}^{\numel} \Omega_e, \quad 
\Omega_e \cap \Omega_l = \emptyset \text{ for } e\neq l ,
\end{equation}
with boundaries $\partial\Omega_e$, defining an internal interface 
\begin{equation}\label{eq:Gamma}
\Gamma := \Big[ \bigcup_{e=1}^{\numel} \partial\Omega_e \Big]\setminus\partial\Omega
\end{equation}

The boundary of a cell is also expressed as the union of a set of edges or faces in two and three dimensions respectively, namely
\begin{equation}
\partial\Omega_e :=  \bigcup_{j=1}^{\numfa^e} \Gamma_{e,j},
\end{equation}
where $\numfa^e$ denotes the number of edges/faces of the cell $\Omega_e$. For triangular cells $\numfa^e=3$ and for tetrahedral cells $\numfa^e=4$.

The proposed FCFV method uses the mixed form of the Poisson problem in the so-called \emph{broken} computational domain, namely 
\begin{equation} \label{eq:PoissonBrokenFirstOrder}
\left\{\begin{aligned}
\bq+\grad u &= \bm{0} &&\text{in $\Omega_e$, and for $e=1,\dotsc ,\numel$,}\\	
\grad\cdot\bq &= s          &&\text{in $\Omega_e$, and for $e=1,\dotsc ,\numel$,}\\
u &= u_D     &&\text{on $\Gamma_D$,}\\
\bn\cdot\bq &= -t         &&\text{on $\Gamma_N$,}\\
\jump{u\bn} &=\bm{0}  &&\text{on $\Gamma$,}\\
\jump{\bn\cdot \bq} &= 0  &&\text{on $\Gamma$,}\\
\end{aligned} \right.
\end{equation}
where, following~\cite{AdM-MFH:08}, the jump operator is defined as the sum of the values from the left and right of an interface, that is $\jump{\odot} = \odot_e + \odot_l$. 

It is worth noting that the last two equations in~\eqref{eq:PoissonBrokenFirstOrder} impose the continuity of the solution and the normal flux across the interface $\Gamma$, respectively.

\subsection{Strong form of the local and global problems}
\label{sc:PoissonStrong}

As it is usual in HDG methods~\cite{Jay-CGL:09,Nguyen-NPC:09,Nguyen-NPC:10,RS-SH:16} and FCFV methods, the strong form of the problem is split into the so-called local problem, defined element-by-element,
\begin{equation} \label{eq:Dlocal-strong}
\left\{\begin{aligned}
\bq_e + \grad u_e &=\bm{0}  &&\text{in $\Omega_e$, }		\\
\grad\cdot\bq_e &= s          &&\text{in $\Omega_e$,}		\\		
u_e &= u_D     &&\text{on $\partial\Omega_e\cap\Gamma_D$,}	\\
u_e &=\hu  &&\text{on $\partial\Omega_e\setminus\Gamma_D$,}
\end{aligned} \right.
\end{equation}
for $e=1,\dotsc ,\numel$, and the global problem, defined over the interface $\Gamma$ and the Neumann boundary
\begin{equation} \label{eq:Dtransmission}
\left\{\begin{aligned}
\jump{u\bn} &=\bm{0}  &&\text{on $\Gamma$,}\\
\jump{\bn\cdot \bq} &= 0  &&\text{on $\Gamma$,}\\
\bn\cdot\bq &= -t         &&\text{on $\Gamma_N$.}\\
\end{aligned} \right.
\end{equation}

The local problem of Equation~\eqref{eq:Dlocal-strong} is a pure Dirichlet problem defined on each cell and introduces the value of the solution at the cell faces, $\hu$, as an independent variable. This problem is used to provide an explicit expression of the solution and its gradient, $u$ and $\bq$, in terms of the solution at the faces, $\hu$. When this explicit expression is introduced in Equation~\eqref{eq:Dtransmission} a global problem with the only unknown being $\hu$ is obtained. This means that the dominant cost of the method is associated to solving the global problem. In a second step the solution and its gradient are computed in each cell using the explicit expression derived from the local problem but this step can be easily parallelised and its cost is negligible compared to the cost of the global problem.

\begin{remark}
The first equation in~\eqref{eq:Dtransmission} is automatically satisfied due to the imposition of the Dirichlet boundary condition in the local problem and the unique value of $\hu$ on the interior faces
\end{remark}

\subsection{Second-order FCFV weak formulation}
\label{sc:PoissonWeak}

Let us denote by $\VoneDomain$ the space of $\eltwo(\Omega)$ functions that are, at most, linear in each cell, $\VzeroDomain$ the space of $\eltwo(\Omega)$ functions that are constant in each cell and  $\VzeroHatSkeleton$ the space of $\eltwo(\Gamma)$ functions that are constant on each cell face. With these definitions, the discrete weak formulation of the local problem is: find $(u_e^h,\bq_e^h) \in  \VoneElem \times [\VzeroElem]^{\nsd}$ such that
\begin{gather} 
- \int_{\Omega_e} \bq_e^h d\Omega = \int_{\partial\Omega_e \cap \Gamma_D} u_D \bn_e  d \Gamma + \int_{\partial\Omega_e\setminus\Gamma_D} \hu^h \bn_e  d \Gamma , \label{eq:weakPoisson1Q}
\\
- \int_{\Omega_e} \grad v \cdot \bq_e^h d\Omega + \int_{\partial\Omega_e}  v (\bn_e\cdot\widehat{\bq}_e^h ) d \Gamma = \int_{\Omega_e} v s d\Omega \label{eq:weakPoisson1U}
\end{gather}
for all $v \in \VoneElem$ and for $e=1,\dotsc ,\numel$. It is worth noting that in Equation~\eqref{eq:weakPoisson1Q}, a constant test function has been arbitrarily chosen in the space $[\VzeroElem]^{\nsd}$ and it has been used that $\bq_e^h \in [\VzeroElem]^{\nsd}$, that is $\grad \cdot \bq_e^h = 0$.

The so-called numerical flux, $\widehat{\bq}_e^h$, is defined as 
\begin{equation} \label{eq:EBENumFlux}
\bn_e\cdot\widehat{\bq}_e^h := \begin{cases}
\bn_e\cdot\bq_e^h + \tau_e (\Pu u_e^h - u_D    ) & \text{on $\partial\Omega_e\cap\Gamma_D$,} \\
\bn_e\cdot\bq_e^h + \tau_e (\Pu u_e^h - \hu^h) & \text{elsewhere,}  
\end{cases}
\end{equation}
where $\Pu$ denotes the projection operator over the space of constant functions~\cite{oikawa2015hybridized,oikawa2016analysis} and $\tau_e$ is a stabilisation parameter that is selected to ensure stability, accuracy and convergence of the resulting scheme~\cite{Jay-CGL:09,Cockburn-CDG:08,Nguyen-NPC:09,Nguyen-NPC:09b,Nguyen-NPC:10,Nguyen-NPC:11}.

\begin{remark} \label{rk:projection}
The definition of the numerical flux considered in Equation~\eqref{eq:EBENumFlux} follows the rationale of the hybridised DG method with reduced stabilisation~\cite{oikawa2015hybridized}. If the projection operator is not considered in Equation~\eqref{eq:EBENumFlux}, the second-order convergence is lost and the resulting method is only first-order~\cite{oikawa2015hybridized}, providing no advantages with respect to the original FCFV. 
\end{remark}

Introducing the expression of the numerical flux in Equation~\eqref{eq:weakPoisson1U} and integrating by parts the first term, leads to the following discrete weak  formulation: find $(u_e^h,\bq_e^h) \in  \VoneElem \times [\VzeroElem]^{\nsd}$ such that
\begin{gather} \label{eq:weakPoisson}
- \int_{\Omega_e} \bq_e^h d\Omega = \int_{\partial\Omega_e \cap \Gamma_D} u_D \bn_e  d \Gamma + \int_{\partial\Omega_e\setminus\Gamma_D} \hu^h \bn_e  d \Gamma ,
\\
\int_{\partial\Omega_e}  v \tau_e \Pu u_e^h d \Gamma  = \int_{\Omega_e} v s d\Omega + \int_{\partial\Omega_e \cap \Gamma_D}  v \tau_e u_D d \Gamma + \int_{\partial\Omega_e\setminus\Gamma_D}  v \tau_e \hu^h d \Gamma
\end{gather}
for all $v \in \VoneElem$ and for $e=1,\dotsc ,\numel$.

The discrete weak form of the global problem is obtained following an analogous process. It becomes: find $\hu^h \in \VzeroFaces$ such that
\begin{equation} \label{eq:HDG-Poisson-Dglobal1}
\sum_{e=1}^{\numel} 
\int_{\partial\Omega_e\setminus\Gamma_D} \bn_e\cdot \widehat{\bq}_e^h d\Gamma
= -\sum_{e=1}^{\numel}  \int_{\partial\Omega_e\cap\Gamma_N} t  d\Gamma,
\end{equation}
a constant test function has been arbitrarily chosen in the space $\VzeroFaces$.

Introducing the definition of the numerical flux in Equation~\eqref{eq:HDG-Poisson-Dglobal1} leads to the following weak form of the global problem: find $\hu^h \in \VzeroFaces$ such that
\begin{equation}\label{eq:HDG-Poisson-Dglobal}
\sum_{e=1}^{\numel} 
\int_{\partial\Omega_e\setminus\Gamma_D} \left(  \bn_e\cdot\bq_e^h + \tau_e (\Pu u_e^h - \hu^h) \right) d\Gamma 
= -\sum_{e=1}^{\numel}  \int_{\partial\Omega_e\cap\Gamma_N} t  d\Gamma.
\end{equation}

\subsection{Second-order FCFV discretisation}
\label{sc:PoissonDiscrete}

To simplify the notation, the following sets of faces are introduced
\begin{equation}
\begin{aligned}
\Aset &:= \{1, \ldots, \numfa^e \}, \\
\Dset &:= \{j \in \Aset \; | \; \Gamma_{e,j} \cap \Gamma_D \neq \emptyset \}, \\
\Nset &:= \{j \in \Aset \; | \; \Gamma_{e,j} \cap \Gamma_N \neq \emptyset \}, \\
\Bset &:= \Aset \setminus \Dset = \{j \in \Aset \; | \; \Gamma_{e,j} \cap \Gamma_D = \emptyset \},
\end{aligned}
\end{equation}
corresponding to all faces of a cell, the faces on the Dirichlet boundary, the faces on the Neumann boundary and the faces not on the Dirichlet boundary, respectively. It is also convenient to denote the set of nodes of a cell $\Omega_e$ belonging to a face $\Gamma_{e,j}$ as $\mathcal{F}_{e,j}$. Finally, the indicator function of a set $\square$ is defined as
\begin{equation}
\chi_{\square}(l) = 
\bigg\{
\begin{array}{ll}
1 & \text{ if } \ l\in\square \\
0 & \text{ otherwise}.
\end{array}
\end{equation}

With this notation, the discrete local problem becomes 
\begin{subequations}\label{eq:HDG-Poisson-DlocalK0}
	\begin{equation}\label{eq:HDG-Poisson-DlocalK0Q}
	- |\Omega_e| \qe =  \sum_{j\in\Dset} |\Gamma_{e,j}| \bn_j u_{D,j}  + \sum_{j\in\Bset} |\Gamma_{e,j}| \bn_j \uHj  ,
	\end{equation}
	\begin{equation}\label{eq:HDG-Poisson-DlocalK0U}
	\mat{m}_e  \bue =  \mat{f}_e + \sum_{j\in\Dset} \tau_j \vect{d}_j + \sum_{j\in\Bset} \tau_j \vect{r}_j \uHj,
	\end{equation}
\end{subequations}
where $\tau_j$ denotes the value of the stabilisation parameter on the $j$-th face, assumed constant, $\qe$ contains the value of $\bq$ in the cell and $\bue$ contains the nodal values of the solution in the cell. The matrices and vectors in Equation~\eqref{eq:HDG-Poisson-DlocalK0U} are given by
\begin{align} 
(m_e)_{IJ} & := \sum_{k=1}^{\numfa^e} \frac{1}{\nfn^{e,i}} |\Gamma_{e,k}| \tau_k \chi_{\mathcal{F}_{e,k}}(I) (\projc)_J, &  (f_e)_I & := \frac{1}{\nen} s_e |\Omega_e|,      \label{eq:matsPoisson1}  \\
(d_j)_I & := \frac{1}{\nfn^{e,i}} u_{D,j} |\Gamma_{e,j}|, &  (r_j)_I & := \frac{1}{\nfn^{e,i}} |\Gamma_{e,j}| \delta_{Ij},    \label{eq:matsPoisson2}
\end{align}
where $\nen$ is the number of cell nodes, the vector $\projv$ is introduced to compute the projection of the solution, i.e. the average of the nodal values of $\bue$ on face $\Gamma_{e,i}$. Formally it is defined as
\begin{equation} \label{eq:vectPoissonProj}
(\projc)_l = \frac{1}{\nfn^{e,i}} \chi_{\mathcal{F}_{e,i}}(l)
\end{equation}
with $\nfn^{e,i}$ being the number of nodes of the face $\Gamma_{e,i}$.

The discrete local problem allows to obtain an explicit expression of both the solution and its gradient in terms of the solution at the cell faces/edges, namely
\begin{subequations}\label{eq:HDG-Poisson-DlocalK0Exp}
	\begin{equation}\label{eq:HDG-Poisson-DlocalK0QExp}
	 \qe =  - |\Omega_e|^{-1}\vect{z}_e  - |\Omega_e|^{-1} \sum_{j\in\Bset} |\Gamma_{e,j}| \bn_j \uHj  ,
	\end{equation}
	\begin{equation}\label{eq:HDG-Poisson-DlocalK0UExp}
	\bue =  \mat{m}_e^{-1} \mat{b}_e + \mat{m}_e^{-1} \sum_{j\in\Bset} \tau_j \vect{r}_j \uHj,
	\end{equation}
\end{subequations}
where 
\begin{equation} \label{eq:poissonPrecomp}
\mat{b}_e := \mat{f}_e  + \sum_{j\in\Dset} \tau_j \vect{d}_j, \qquad
\vect{z}_e := \sum_{j\in\Dset} |\Gamma_{e,j}| \bn_j u_{D,j} .
\end{equation}

It is worth noting that the Equation~\eqref{eq:HDG-Poisson-DlocalK0UExp} involves the solution of a $3 \times 3$ system of equations for triangular cells and a $4 \times 4$ system for tetrahedral cells. Given the size of the system and the  definition of $\mat{m}_e$, its inverse can be analytically computed to substantially reduce the computational cost of this operation.

The discretisation of the global problem of Equation~\eqref{eq:HDG-Poisson-Dglobal} leads, for $i \in \Bset$, to 
\begin{equation}\label{eq:HDG-Poisson-DglobalK0}
\sum_{e=1}^{\numel}\Bigl\{
|\Gamma_{e,i}| \bn_i \cdot \qe + |\Gamma_{e,i}| \tau_i \projv \cdot \bue - |\Gamma_{e,i}| \tau_i \uHi \Bigr\}
= -\sum_{e=1}^{\numel} \Bigl\{|\Gamma_{e,i}| t_i \, \chi_{\Nset}(i) \Bigr\} \quad.
\end{equation}

By inserting the explicit expressions of Equation~\eqref{eq:HDG-Poisson-DlocalK0Exp}, in the global problem of Equation~\eqref{eq:HDG-Poisson-DglobalK0}, a linear system of equations involving only the solution on the faces as an unknown is obtained, namely
\begin{equation} \label{eq:globalSystemPoisson}
\mat{\widehat{K}}\vect{\hu}=\vect{\hat{f}}.
\end{equation}
The global matrix $\mat{\widehat{K}}$ and right hand side $\vect{\hat{f}}$ are obtained by assembling the in the cell contributions given by
\begin{subequations}\label{HDG-Poisson-globalSystem}
	\begin{align}
	&
	{\widehat{K}}^e_{i,j} :=  |\Gamma_{e,i}| \Big(  \tau_i \tau_j \projv \cdot \left( \mat{m}_e^{-1} \vect{r}_j \right) - |\Omega_e|^{-1} |\Gamma_{e,j}| \bn_i \cdot \bn_j - \tau_i \delta_{ij} \Big) 
	, \\
	&
	{\widehat{f}}^e_i :=  |\Gamma_{e,i}| \Big( |\Omega_e|^{-1}  \bn_i \cdot \vect{z}_e - \tau_i \projv \cdot \left( \mat{m}_e^{-1} \mat{b}_e\right) - t_i \, \chi_{\Nset}(i)  \Big),
	\end{align}
\end{subequations}
for $i,j \in \Bset$ and with $\delta_{ij}$ denoting the Kronecker delta.

As discussed in Remark~\ref{rk:projection}, the use of the projection operator in the numerical flux of Equation~\eqref{eq:EBENumFlux} is required to obtain a second-order method. If the projection is not introduced, the resulting method is only first-order~\cite{oikawa2015hybridized}. This minor difference is exploited in Section~\ref{sc:adaptivity} to devise an error indicator that can be used to drive an automatic mesh adaptive process. Therefore, it is of interest here to study the difference in the global system of Equation~\ref{HDG-Poisson-globalSystem} induced by the introduction of the projection operator. The next result shows that only the matrix $\mat{m}_e$ changes if the projection is not considered.

\begin{lemma} \label{rk:projDifference}
Let us consider a face $\Gamma_{e,i}$ of a triangular or tetrahedral cell $\Omega_e$ and $u_e^h \in  \VoneElem$. Then, the following equality holds
\begin{equation} \label{eq:projectionGlobalProblem}
\int_{\Gamma_{e,i}} \Pu u_e^h d\Gamma = \int_{\Gamma_{e,i}} u_e^h d\Gamma.
\end{equation}

\begin{proof}
The first integral of Equation~\eqref{eq:projectionGlobalProblem} can be written as
\begin{equation}
\int_{\Gamma_{e,i}} \Pu u_e^h d\Gamma = |\Gamma_{e,i}| \Pu u_e^h = |\Gamma_{e,i}| \projv \cdot \bue
\end{equation}
because $\Pu u_e^h$ is constant within each face.

The second integral of Equation~\eqref{eq:projectionGlobalProblem} can be easily computed using the expression of the linear shape functions used to define the approximation $u_e^h$, namely
\begin{equation}
\int_{\Gamma_{e,i}} u_e^h d\Gamma = \sum_{l=1}^{\nen} \left( \int_{\Gamma_{e,i}} N_l d\Gamma \right) u_l = \sum_{l=1}^{\nen} \left( \frac{1}{\nfn^{e,i}} |\Gamma_{e,i}| \chi_{\mathcal{F}_{e,i}}(l) \right) u_l.
\end{equation}
where $N_l$ and $u_l$ denote the linear shape function and nodal value associated with the $l$-th node.

Equality~\eqref{eq:projectionGlobalProblem} follows from the definition of $\projv$ introduced in Equation~\eqref{eq:vectPoissonProj}.
\end{proof}
\end{lemma}

\section{Second-order FCFV for the Stokes equation}
\label{sc:Stokes}

\subsection{Problem statement}
\label{sc:StokesStatement}

The strong form of the Stokes problem can be written in the partitioned domain and after splitting the second-order momentum conservation equation into a system of two first-order equations, as
\begin{equation} \label{eq:StokesBrokenMixed}
\left\{\begin{aligned}
\bL + \sqrt{\nu} \grad\bu & = \bm{0} &&\text{in $\Omega_e$, and for $e=1,\dotsc ,\numel$,}\\
\grad\cdot\bigl(\sqrt{\nu} \bL + p\Insd \bigr) &= \bm{s}         &&\text{in $\Omega_e$, and for $e=1,\dotsc ,\numel$,}\\
\grad\cdot\bu &= 0  &&\text{in $\Omega_e$, and for $e=1,\dotsc ,\numel$,}\\
\bu &= \bu_D     &&\text{on $\Gamma_D$,}\\
\bn \cdot \bigl(\sqrt{\nu} \bL + p \Insd\bigr)  &= - \bm{t}  &&\text{on $\Gamma_N$,}\\
\jump{\bu \otimes \bn} &=\bm{0}  &&\text{on $\Gamma$,}\\
\jump{\bn \cdot \bigl(\sqrt{\nu} \bL + p \Insd\bigr)} &= \bm{0}  &&\text{on $\Gamma$.}\\
\end{aligned} \right.
\end{equation}
where $\nu>0$ is the viscosity and the last two equations enforce the continuity of the velocity and the normal flux across the interface $\Gamma$ respectively. 

\begin{remark}
To simplify the presentation, this work considers the traditional velocity-pressure HDG formulation of the Stokes equation~\cite{HDGNEFEMstokes}, where the vector $\bm{t}$ does not correspond to the boundary traction and it is usually called a pseudo-traction~\cite{donea2003finite}. It is worth emphasising that the so-called Cauchy formulation could also be employed here, using the formulation proposed in~\cite{giacomini2018superconvergent}. This formulation, contrary to other existing HDG methods, guarantees optimal convergence even for low order approximations. This idea was also exploited in the context of linear elasticity to obtain optimal convergence for low order approximations in an HDG context~\cite{HDGelasticityVoigt,FCFVelas,sevilla2019hdg}.
\end{remark}

\subsection{Strong form of the local and global problems}
\label{sc:StokesStrong}

As usually done in HDG methods~\cite{Nguyen-CNP:10,oikawa2016analysis,giacomini2018superconvergent} and FCFV methods, the strong form of the problem given by Equation~\eqref{eq:StokesBrokenMixed} is split into the local and global problems. The local problem
\begin{equation} \label{eq:StokesLocalStrong}
\left\{\begin{aligned}
\bL_e + \sqrt{\nu} \grad\bu_e & = \bm{0} &&\text{in $\Omega_e$,}\\
\grad\cdot\bigl(\sqrt{\nu} \bL_e + p_e\Insd \bigr) &= \bm{s}         &&\text{in $\Omega_e$,}\\
\grad\cdot\bu_e &= 0  &&\text{in $\Omega_e$}\\
\bu_e &= \bu_D     &&\text{on $\partial\Omega_e \cap \Gamma_D$,}\\
\bu_e &= \bhu     &&\text{on $\partial\Omega_e \setminus \Gamma_D$.}\\                       
\end{aligned} \right.
\end{equation}
is a pure Dirichlet problem and, therefore, requires the introduction of a solvability constraint for the pressure~\cite{Nguyen-CNP:10}, namely
\begin{equation}\label{eq:pressureConstraint}
\frac{1}{|\partial\Omega_e|} \int_{\partial\Omega_e}  p_e d\Gamma  = \rho_e
\end{equation}
where $\rho_e$ is the mean value of the pressure on the boundary of the cell $\Omega_e$.

In addition, a compatibility condition is induced by the free-divergence character of the velocity, namely
\begin{equation}\label{eq:divergenceFreeConstraint}
\int_{\partial\Omega_e\setminus\Gamma_D} \bhu \cdot \bn_e d\Gamma + \int_{\partial\Omega_e\cap\Gamma_D} \bu_D \cdot \bn_e  d\Gamma = 0.
\end{equation}

As the continuity of the solution is automatically imposed by the introduction of the velocity on the cell faces as an independent variable that is uniquely defined on each face, the global problem imposes the continuity of the normal flux across the interface and the Neumann boundary conditions, that is
\begin{equation} \label{eq:StokesGlobalStrong}
\left\{\begin{aligned}
\jump{\bn \cdot \bigl(\sqrt{\nu} \bL + p \Insd\bigr)} &= \bm{0}  &&\text{on $\Gamma$,}\\
\bn \cdot \bigl(\sqrt{\nu} \bL + p \Insd\bigr)  &= - \bm{t}  &&\text{on $\Gamma_N$.}\\  		
\end{aligned} \right.
\end{equation}

\subsection{Second-order FCFV weak formulation}
\label{sc:StokesWeak}

Following the same rationale presented for the Poisson problem, the discrete weak formulation of the local Stokes problem is: find $(\bu_e^h,\bL_e^h,p_e^h) \in  [\VoneElem]^{\nsd} \times [\VzeroElem]^{\nsd\times\nsd} \times [\VzeroElem]$ such that
\begin{subequations}\label{HDG-Stokes-DlocalSymm}
	\begin{align}
	&
	-\int_{\Omega_e} \bL^h_e d\Omega = \int_{\partial\Omega_e\cap\Gamma_D}  \sqrt{\nu} \bn_e \otimes \bu_D d\Gamma + \int_{\partial\Omega_e\setminus\Gamma_D}  \sqrt{\nu} \bn_e \otimes \bhu^h d\Gamma, \label{eq:HDG-Stokes-Dlocal-LSymm} 
	\\
	&
	\int_{\partial \Omega_e} \tau_e \bw \cdot \Pu \bu^h_e  d\Gamma
	=  \int_{\Omega_e}\bw \cdot \bm{s} d\Omega   + \int_{\partial\Omega_e\cap\Gamma_D}  \tau_e \bw \cdot \bu_D d\Gamma + \int_{\partial\Omega_e\setminus\Gamma_D}  \tau_e \bw \cdot \bhu^h d\Gamma,
 	 \label{eq:HDG-Stokes-Dlocal-MomSymm}
	\\
	&
	\int_{\partial\Omega_e\setminus\Gamma_D} \bhu \cdot \bn_e d\Gamma + \int_{\partial\Omega_e\cap\Gamma_D} \bu_D \cdot \bn_e  d\Gamma = 0 , 
	\label{eq:HDG-Stokes-Dlocal-USymm} 
	\\
	&
	\frac{1}{|\partial\Omega_e|} \int_{\partial\Omega_e}  p^h_e d\Gamma  = \rho_e , \label{eq:HDG-Stokes-Dlocal-PSymm}
	\end{align}
\end{subequations}
for all $\bw \in [\VoneElem]^{\nsd}$. It is worth emphasising that in Equations~\eqref{eq:HDG-Stokes-Dlocal-LSymm} and \eqref{eq:HDG-Stokes-Dlocal-USymm}, a constant test function has been arbitrarily chosen in the spaces $[\VzeroElem]^{\nsd \times \nsd}$ and $\VzeroElem$ respectively. 
It is also worth noting that Equation~\eqref{eq:HDG-Stokes-Dlocal-USymm}  is exactly the compatibility condition introduced in Equation~\eqref{eq:divergenceFreeConstraint}. As done in the standard FCFV method, this equation is then removed from the local problem and imposed only in the global problem. 

The weak form of the local problem has been introduced after using the following definition of the numerical flux
\begin{equation} \label{eq:traceStokes}
\bn_e \! \cdot \!\bigl(\widehat{\sqrt{\nu} \bL^h_e \!+\! p^h_e \Insd} \bigr) \!:= \!
\begin{cases}
\bn_e \! \cdot \! \bigl( \sqrt{\nu} \bL^h_e \!+\! p^h_e \Insd \bigr) \!+\! \tau_e (\Pu \bu^h_e \!-\! \bu_D) & \text{on $\partial\Omega_e\cap\Gamma_D$,} \\
\bn_e \!\cdot \!\bigl( \sqrt{\nu} \bL^h_e \!+\! p^h_e \Insd \bigr) \!+\! \tau_e (\Pu \bu^h_e \!-\! \bhu^h) & \text{elsewhere.}  
\end{cases}
\end{equation}

The discrete global problem that accounts for the transmission conditions, the Neumann boundary condition and the incompressibility constraint reads: find $(\bhu^h,\rho_e) \in [\VzeroFaces]^{\nsd} \times \mathbb{R}$ such that
\begin{subequations}\label{HDG-Stokes-Dglobal}
	\begin{align}
	&
	\sum_{e=1}^{\numel} 
	\int_{\partial\Omega_e\setminus\Gamma_D} \left(  \sqrt{\nu} \bn_e\cdot \bL^h_e + p^h_e \bn_e + \tau_e (\Pu \bu^h_e - \bhu^h) \right) d\Gamma 
	= -\sum_{e=1}^{\numel}  \int_{\partial\Omega_e\cap\Gamma_N} \bm{t}  d\Gamma,
	\label{eq:HDG-Stokes-Dglobal1}
	\\
	&
	\int_{\partial\Omega_e\setminus\Gamma_D} \bhu \cdot \bn_e d\Gamma + \int_{\partial\Omega_e\cap\Gamma_D} \bu_D \cdot \bn_e  d\Gamma = 0 \quad \text{ for } e=1,\dotsc,\numel.
	\label{eq:weakCompatibilityStokes}
	\end{align}
\end{subequations}
It is worth noting that a constant test function in the space $[\VzeroFaces]^{\nsd}$ has been used in Equation~\eqref{eq:HDG-Stokes-Dglobal1}.

\subsection{Second-order FCFV discretisation}
\label{sc:StokesDiscrete}

Using the notation introduced in Section~\ref{sc:PoissonDiscrete}, the discrete local Stokes problem provides explicit expressions of the velocity, its gradient and the pressure in terms of the velocity on the faces and the mean pressure. The expressions are
\begin{subequations}\label{HDG-Stokes-localK0}
	\begin{align}
	&
	 \Le =  -| \Omega_e |^{-1}\sqrt{\nu}\mat{Z}_e    -| \Omega_e |^{-1} \sqrt{\nu} \sum_{j \in \Bset} | \Gamma_{e,j} | \bn_j \otimes \buHj  , \label{eq:HDG-Stokes-Dlocal-LK0} 
	\\
	& \bue =  \mat{M}_e^{-1} \mat{B}_e + \mat{M}_e^{-1} \sum_{j\in\Bset} \tau_j \mat{R}_j \buHj, \label{eq:HDG-Stokes-Dlocal-MomK0} 
	\\
	&
	\pe = \rho_e, \label{eq:HDG-Stokes-Dlocal-PK0}
	\end{align}
\end{subequations}
for $e=1,\dotsc,\numel$, where
\begin{equation} \label{eq:stokesPrecomp}
	\mat{B}_e := \mat{F}_e  + \sum_{j\in\Dset} \tau_j \mat{D}_j
	, \qquad
	\mat{Z}_e := \sum_{j\in\Dset} |\Gamma_{e,j}| \bn_j \otimes \bu_{D,j} .
\end{equation}
and
\begin{align} 
	(\mat{M}_e)_{IJ} & := \Insd \sum_{k=1}^{\numfa^e} \frac{1}{\nfn^{e,i}} |\Gamma_{e,k}| \tau_k \chi_{\mathcal{F}_{e,k}}(I) (\projc)_J, &  \mat{F}_I & := \frac{1}{\nen} \bm{s}_e |\Omega_e|,      \label{eq:matsStokes1}  \\
	(\mat{D}_j)_I & := \frac{1}{\nfn^{e,i}} \bu_{D,j} |\Gamma_{e,j}|, &  (\mat{R}_j)_{IJ} & := \Insd \frac{1}{\nfn^{e,i}} |\Gamma_{e,j}| \delta_{Ij}.    \label{eq:matsStokes2}
\end{align}

Analogously, the discretisation of the global problem of Equation~\eqref{HDG-Stokes-Dglobal} after inserting the expressions of Equation~\eqref{HDG-Stokes-localK0}, leads to a system of equations
\begin{equation} \label{eq:globalSystemStokes}
\begin{bmatrix}
\mat{\widehat{K}}_{\hu \hu} & \mat{\widehat{K}}_{\hu \rho} \\
\mat{\widehat{K}}_{\hu \rho}^T & \mat{0}_{\numel}
\end{bmatrix}
\begin{Bmatrix}
\vect{\hat{u}}  \\
\bm{\rho}
\end{Bmatrix}
=
\begin{Bmatrix}
\vect{\hat{f}}_{\hu} \\
\vect{\hat{f}}_{\rho}
\end{Bmatrix} ,
\end{equation}
where the global matrix and right hand side are obtained by assembling the elemental contributions given by
\begin{subequations}\label{HDG-Stokes-globalSystem}
	\begin{align}
	&
	(\mat{\widehat{K}}_{\hu \hu})^e_{i,j} :=  |\Gamma_{e,i}| \left[ \tau_i \tau_j \projm \left( \mat{M}_e^{-1} \mat{R}_j \right) - \nu | \Omega_e |^{-1} |\Gamma_{e,j}| (\bn_i \cdot \bn_j)\mat{I}_{\nsd}  - \tau_i \delta_{ij} \mat{I}_{\nsd} \right] 
	, \\
	&
	(\mat{\widehat{K}}_{\hu \rho})^e_i := |\Gamma_{e,i}| \bn_i , \\
	&
	(\vect{\hat{f}}_{\hu})^e_i :=  |\Gamma_{e,i}| \left( \nu | \Omega_e |^{-1} \bn_i \cdot \mat{Z}_e  - \tau_i \projm \left( \mat{M}_e^{-1} \mat{B}_e \right)  - \bm{t}_i \, \chi_{\Nset}(i) \right) , \\
	&
	(\hat{\text{f}}_{\rho})^e := -\sum_{j \in \Dset} | \Gamma_{e,j} | \bu_{D,j} \cdot \bn_j,
	\end{align}
\end{subequations}
for $i,j \in \Bset$. The matrix $\projm$, introduced to account for the projection of the solution on a space of constant functions on the face $\Gamma_{e,i}$, is defined as
\begin{equation} \label{eq:matStokesProj-2D}
\projm = 
\begin{bmatrix}
\projv^T &  \bm{0}_{1 \times 3}  \\
\bm{0}_{1 \times 3} &  \projv^T
\end{bmatrix}
\end{equation}
and
\begin{equation} \label{eq:matStokesProj-3D}
\projm = 
\begin{bmatrix}
\projv^T &  \bm{0}_{1 \times 4} &  \bm{0}_{1 \times 4}  \\
\bm{0}_{1 \times 4} &  \projv^T &  \bm{0}_{1 \times 4}  \\
\bm{0}_{1 \times 4}  &  \bm{0}_{1 \times 4}&  \projv^T
\end{bmatrix}
\end{equation}
in two and three dimensions respectively. The vector $\projv$ was introduced in Equation~\eqref{eq:vectPoissonProj}.

\section{Mesh adaptivity}
\label{sc:adaptivity}

As mentioned in Section~\ref{sc:PoissonWeak}, the proposed FCFV can be seen as a particular case of the hybridised DG method with reduced stabilisation~\cite{oikawa2015hybridized} and therefore provides second-order convergence for the solution. The key aspect is the projection of the solution over a space of constant functions. Without this projection the method is only first-order accurate. This small difference in the formulation is exploited here to devise an error indicator. Noting $u$ the solution of the proposed FCFV methodology and $\tilde{u}$ the solution of the method where the projection is not performed, the following error indicator is proposed for the Poisson problem
\begin{equation} \label{eq:errorMeasure}
E_e = \left[ \frac{1}{|\Omega_e|} \int_{\Omega_e} \left( \tilde{u} - u \right)^2 d\Omega \right]^{1/2},
\end{equation}
with a similar definition for the Stokes problem.

To compute the desired cell size, the error indicator of Equation~\eqref{eq:errorMeasure} is combined with the a priori local error estimate for elliptic problems~\cite{diez1999unified} given by
\begin{equation} \label{eq:aPrioriError}
\varepsilon_e = \| u - u_h \|_{\eltwo(\Omega_e)} \leq C h_e^{1 + \nsd/2},
\end{equation}
for a constant degree of approximation. By using Richardson extrapolation, the following desired cell size is computed, for a desired error $\varepsilon$,
\begin{equation} \label{eq:newElemSize}
h_e^\star = h_e \left(\frac{\varepsilon}{E_e}\right)^{2 + \nsd/2},
\end{equation}
where $h_e$ is the characteristic cell size used to perform the computation of both $u$ and $\tilde{u}$.

It is worth noting that the difference between $u$ and $\tilde{u}$ is only due to the use of the projection of the solution over a space of constant functions. Therefore, the majority of the calculations required to assemble the global system of equations can be re-used, substantially reducing the computational effort required to compute the error indicator of Equation~\eqref{eq:errorMeasure}. In fact, as detailed in Remark~\ref{rk:projDifference}, the only difference between both formulations is in the matrix $\mat{m}_e$ in Equation~\eqref{HDG-Poisson-globalSystem} and $\mat{M}_e$ in Equation~\eqref{HDG-Stokes-globalSystem} for the Poisson and Stokes problems respectively.

\section{Numerical studies}
\label{sc:studies}

This Section presents a series of numerical experiments designed to test the optimal convergence properties of the proposed technique and to compare its performance with the recently proposed first-order FCFV~\cite{FCFV2018}. Numerical experiments are also presented to illustrate the accuracy of the method in terms of the stabilisation parameter and the distortion and the stretching of the meshes. A detailed description of the model problems considered for both Poisson and Stokes problems in two and three dimensions can be found in~\cite{FCFV2018} and they are omitted here for brevity. 

\subsection{Optimal convergence of the second-order FCFV scheme for Poisson equation}
\label{sc:PoissonVerification}

\begin{figure}[!tb]
	\centering
	\subfigure[$u$]{\includegraphics[width=0.48\textwidth]{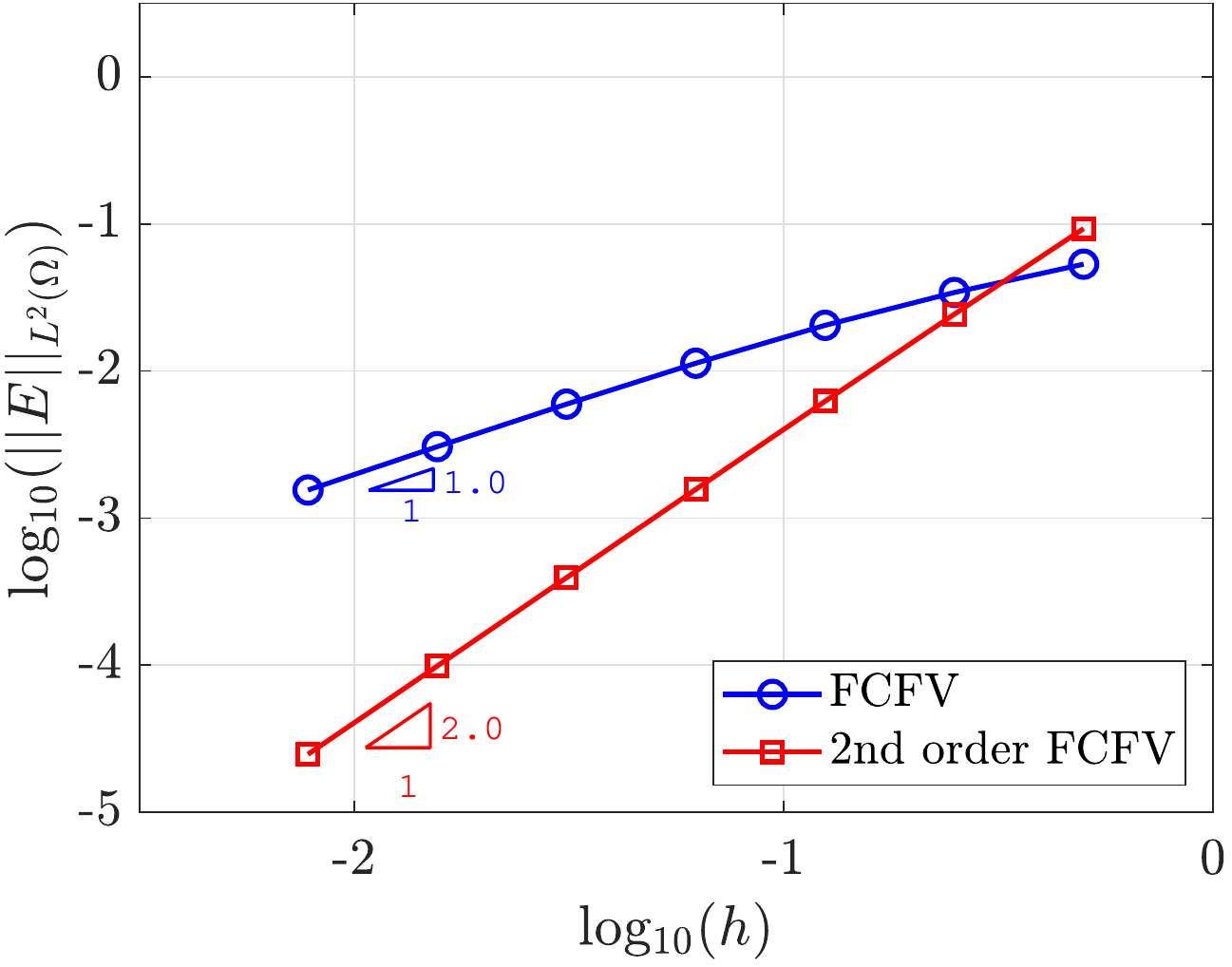}}
	\subfigure[$\bq$]{\includegraphics[width=0.48\textwidth]{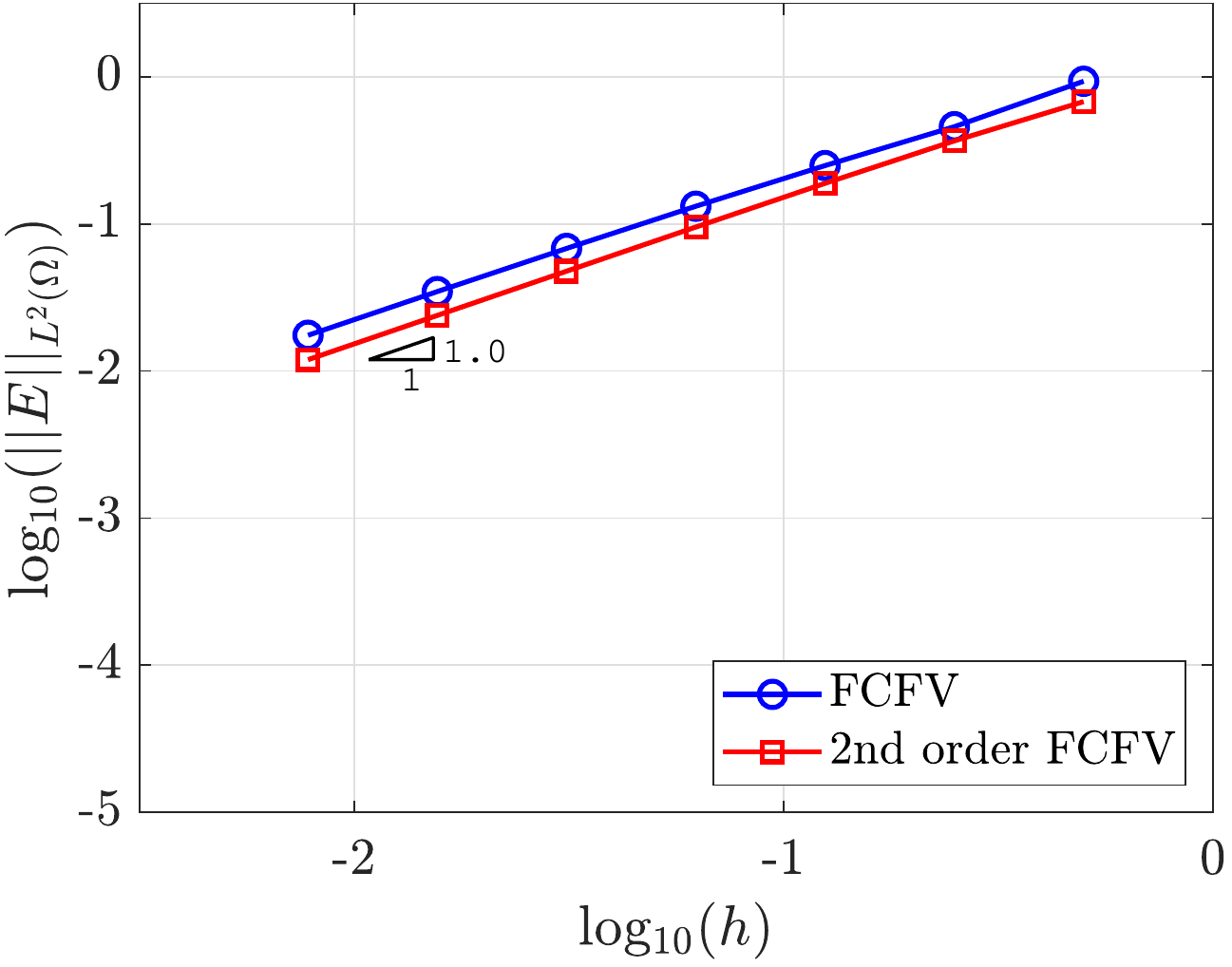}}
	\caption{Mesh convergence of the error of the solution and its gradient in the \eltwo($\Omega$) norm for two dimensional Poisson problem.}
	\label{fig:Poisson_Conv_2D}
\end{figure}

A mesh convergence study is performed for the Poisson problem using a series of successively refined triangular meshes. Figure~\ref{fig:Poisson_Conv_2D} shows the relative $\eltwo(\Omega)$ norm of the error of the solution and its gradient as a function of the characteristic mesh size $h$. 
The results clearly show the optimal second-order convergence of the error of the solution and the first-order convergence of the error of the solution gradient. The results of the original FCFV are also included, clearly showing the gain in accuracy achieved for the solution $u$ for a given spatial discretisation. The gradient of the solution is only marginally more accurate as the approximation space for this variable is not changed with respect to the original FCFV.

\begin{figure}[!tb]
	\centering
	\subfigure[$u$]{\includegraphics[width=0.48\textwidth]{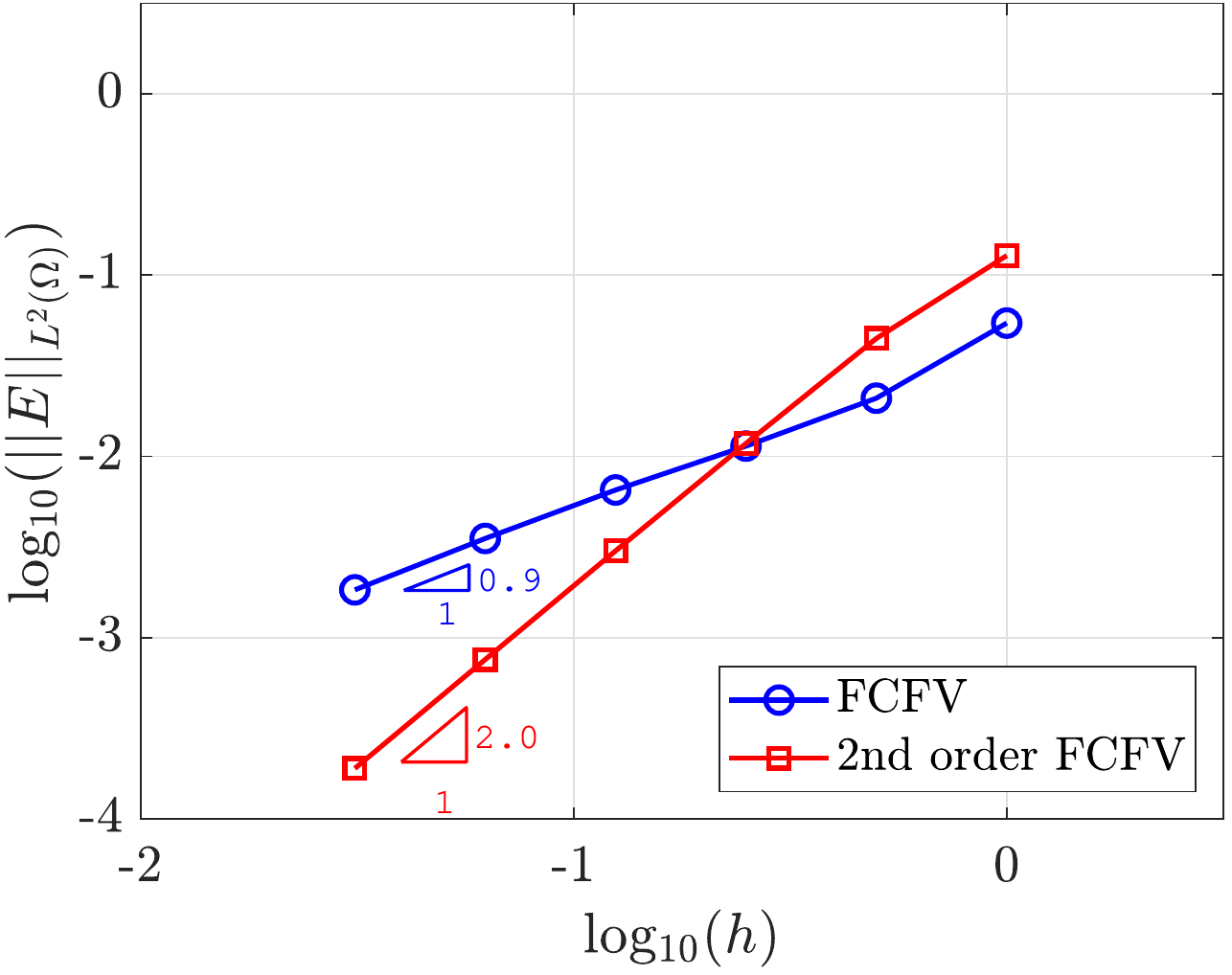}}
	\subfigure[$\bq$]{\includegraphics[width=0.48\textwidth]{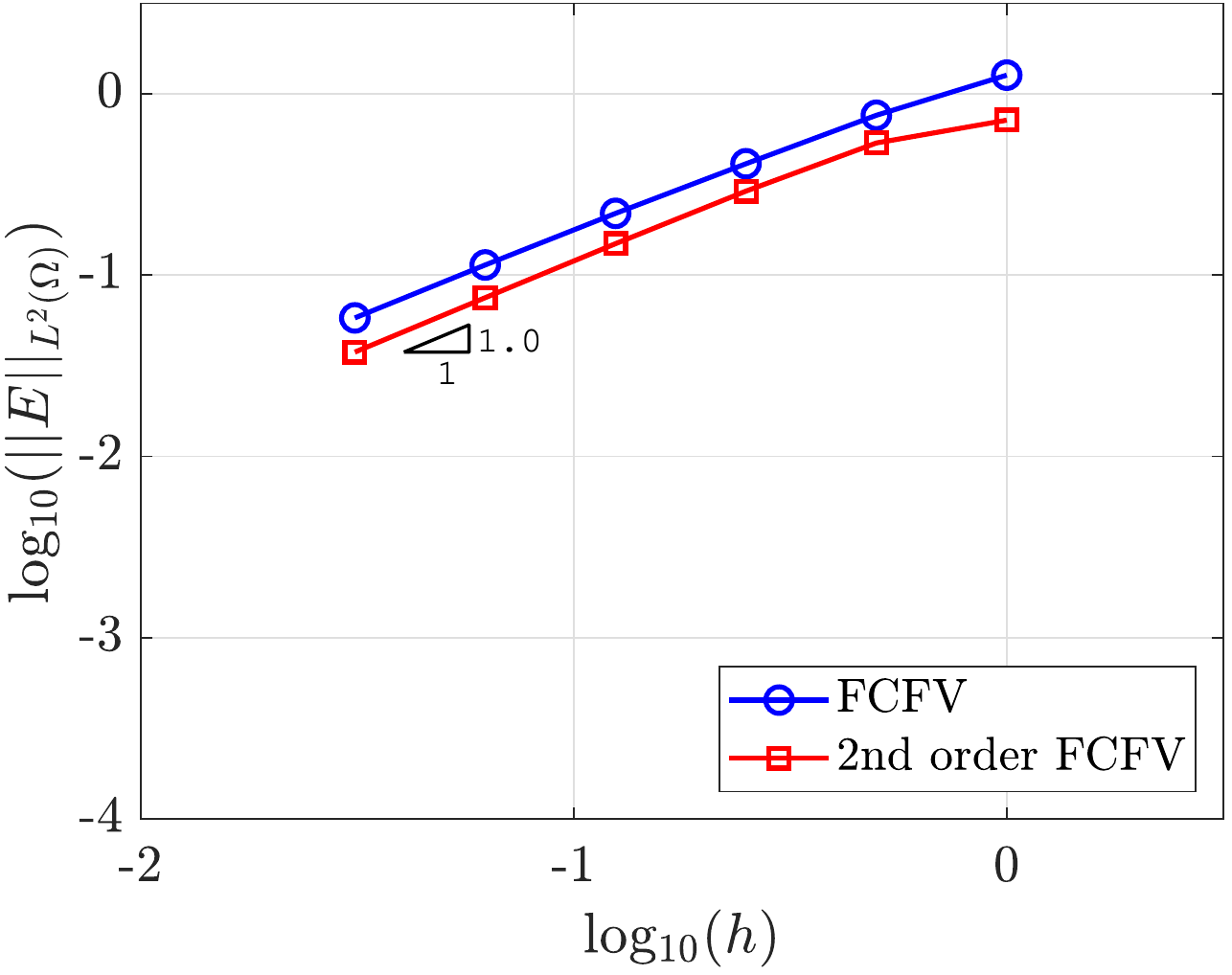}}
	\caption{Mesh convergence of the error of the solution and its gradient in the \eltwo($\Omega$) norm for three dimensional Poisson problem.}
	\label{fig:Poisson_Conv_3D}
\end{figure}

The results of the analogous study in three dimensions are shown in Figure~\ref{fig:Poisson_Conv_3D}, demonstrating the optimal convergence of the error for both the solution and its gradient as well as the increased accuracy with respect to the first-order FCFV.

\subsection{Optimal convergence of the second-order FCFV scheme for Stokes equation}
\label{sc:StokesVerification}

\begin{figure}[!tb]
	\centering	
	\subfigure[$\bu$]{\includegraphics[width=0.32\textwidth]{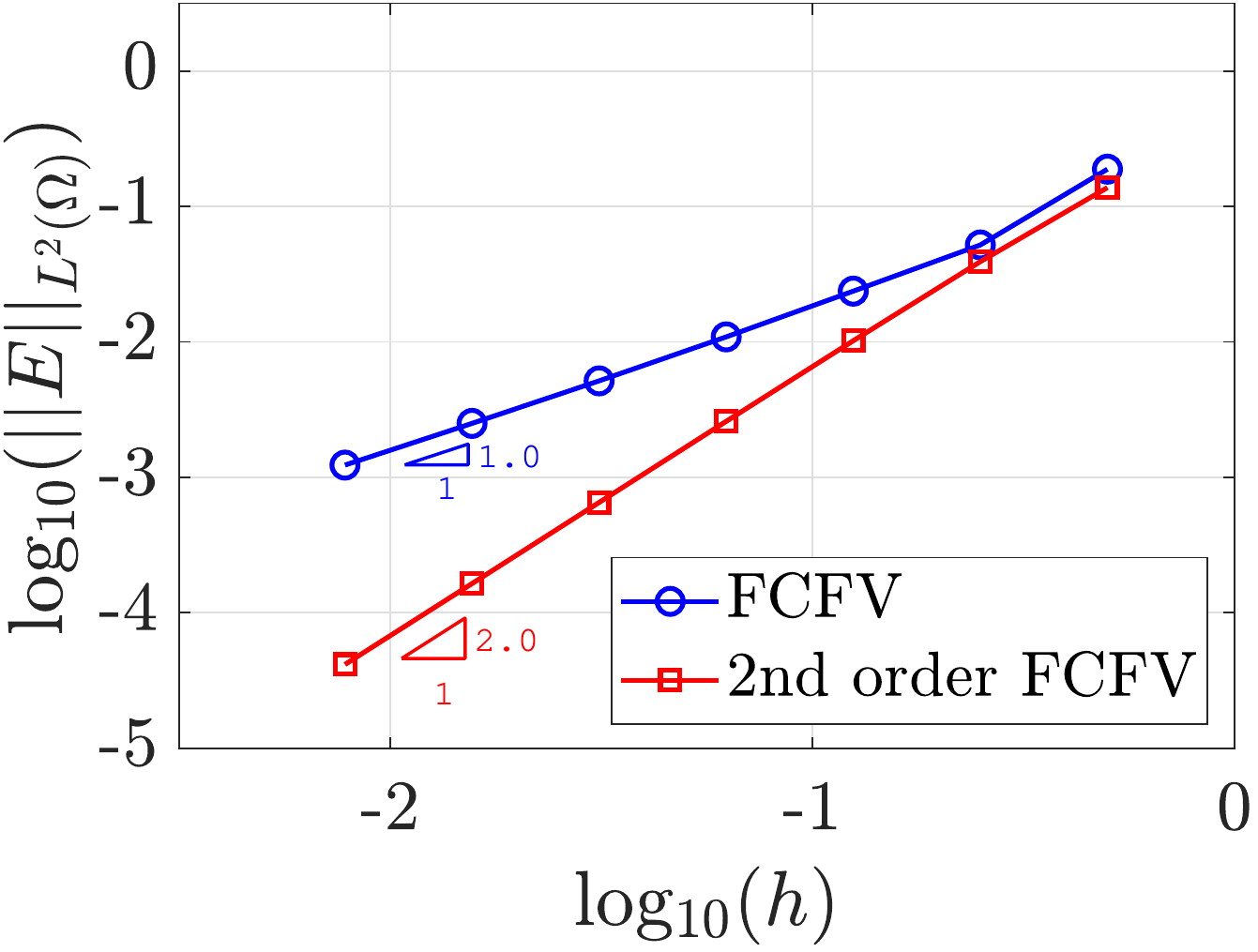}}
	\subfigure[$\bL$]{\includegraphics[width=0.32\textwidth]{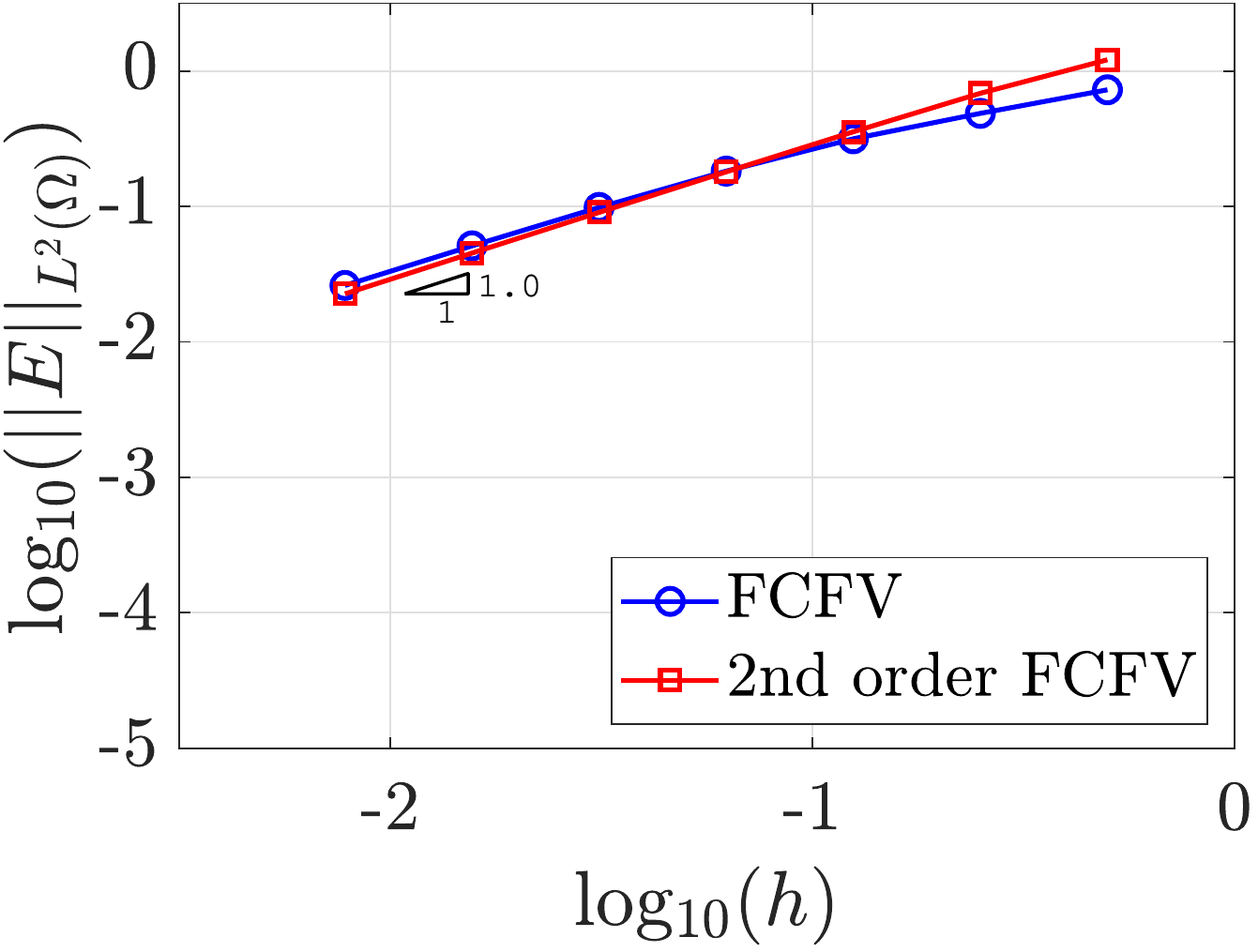}}
	\subfigure[$p$]{\includegraphics[width=0.32\textwidth]{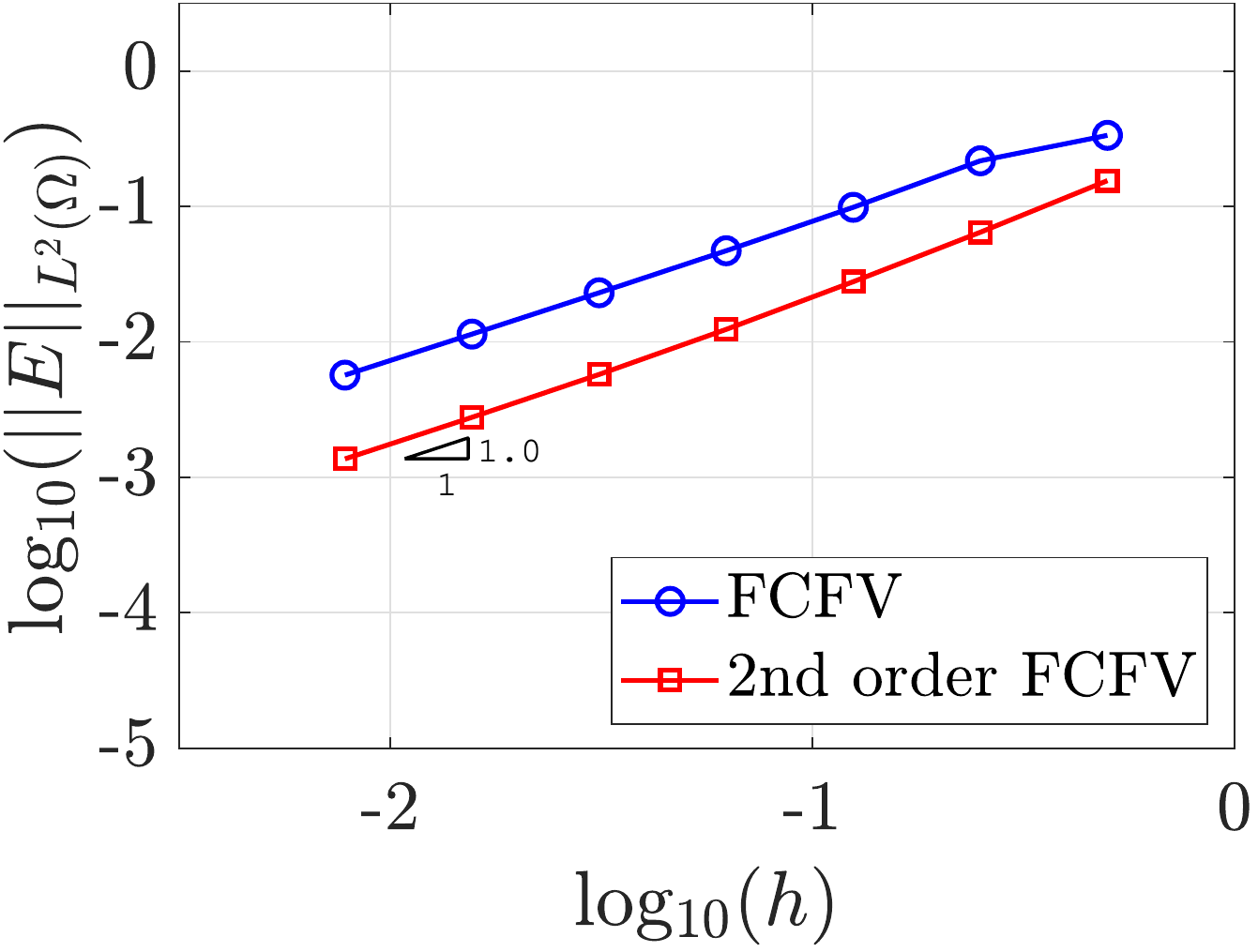}}
	\caption{Mesh convergence of the error of the velocity, its gradient and the pressure in the \eltwo($\Omega$) norm for two dimensional Stokes problem.}
	\label{fig:Stokes_Conv_2D}
\end{figure}

\begin{figure}[!tb]
	\centering	
	\subfigure[$\bu$]{\includegraphics[width=0.32\textwidth]{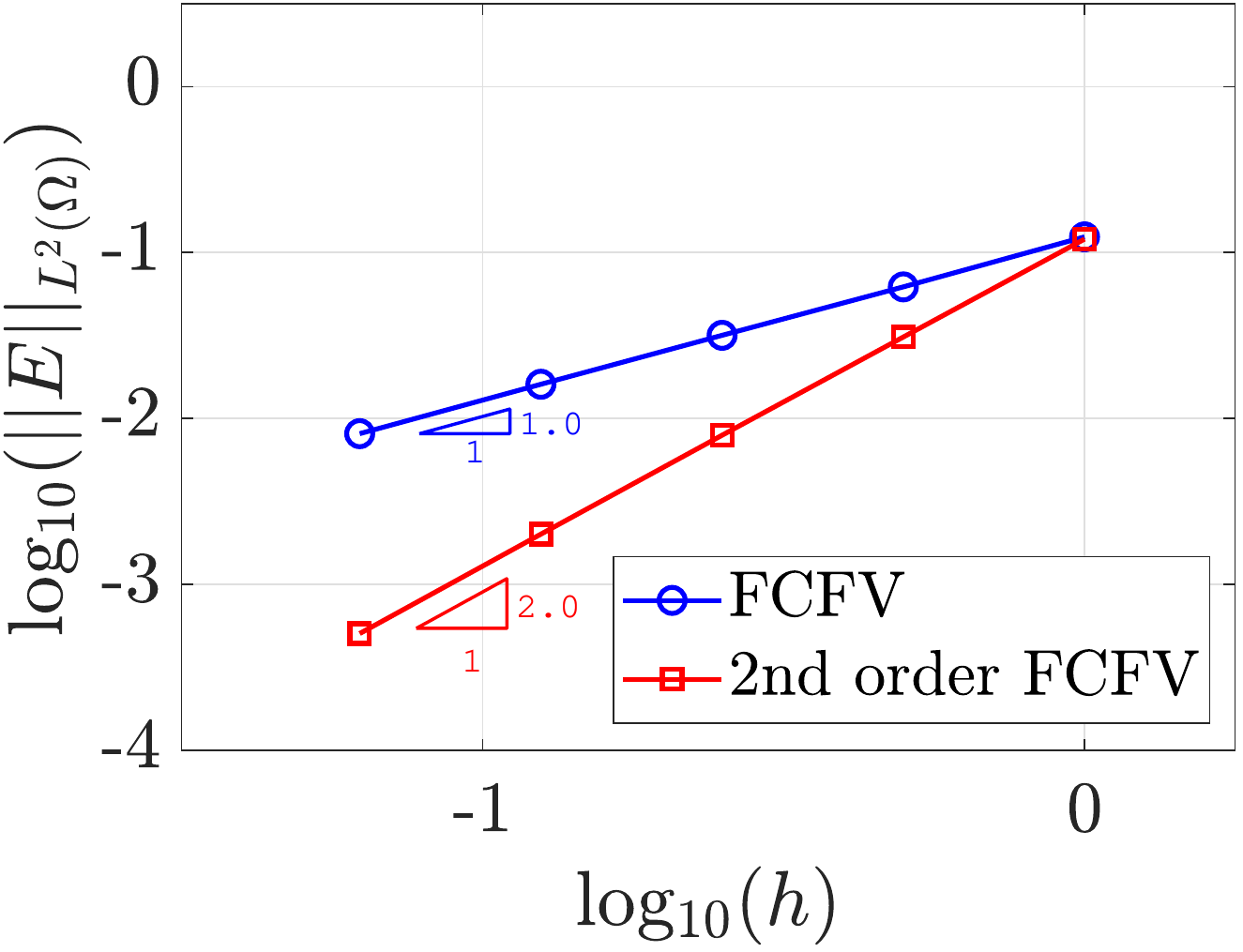}}
	\subfigure[$\bL$]{\includegraphics[width=0.32\textwidth]{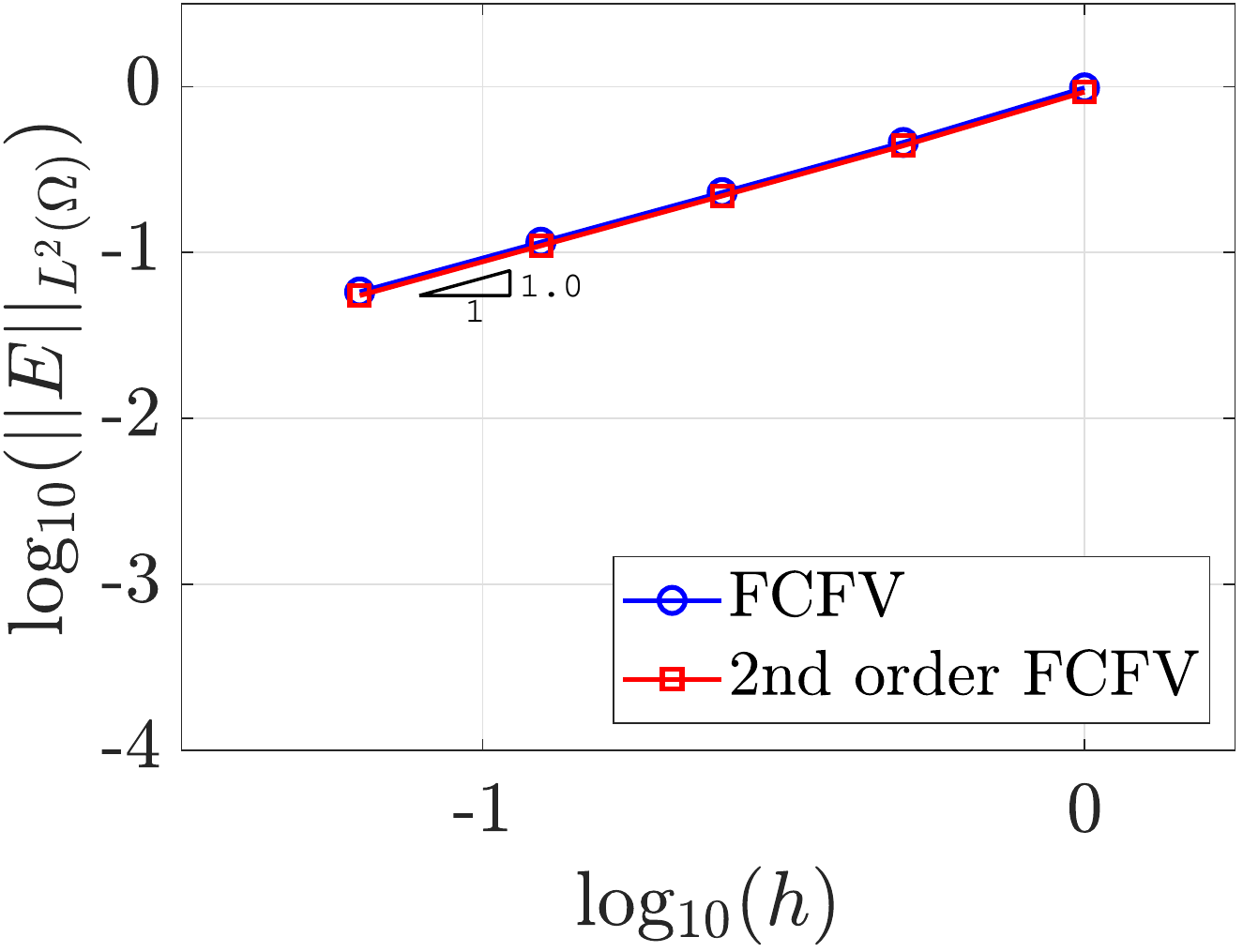}}
	\subfigure[$p$]{\includegraphics[width=0.32\textwidth]{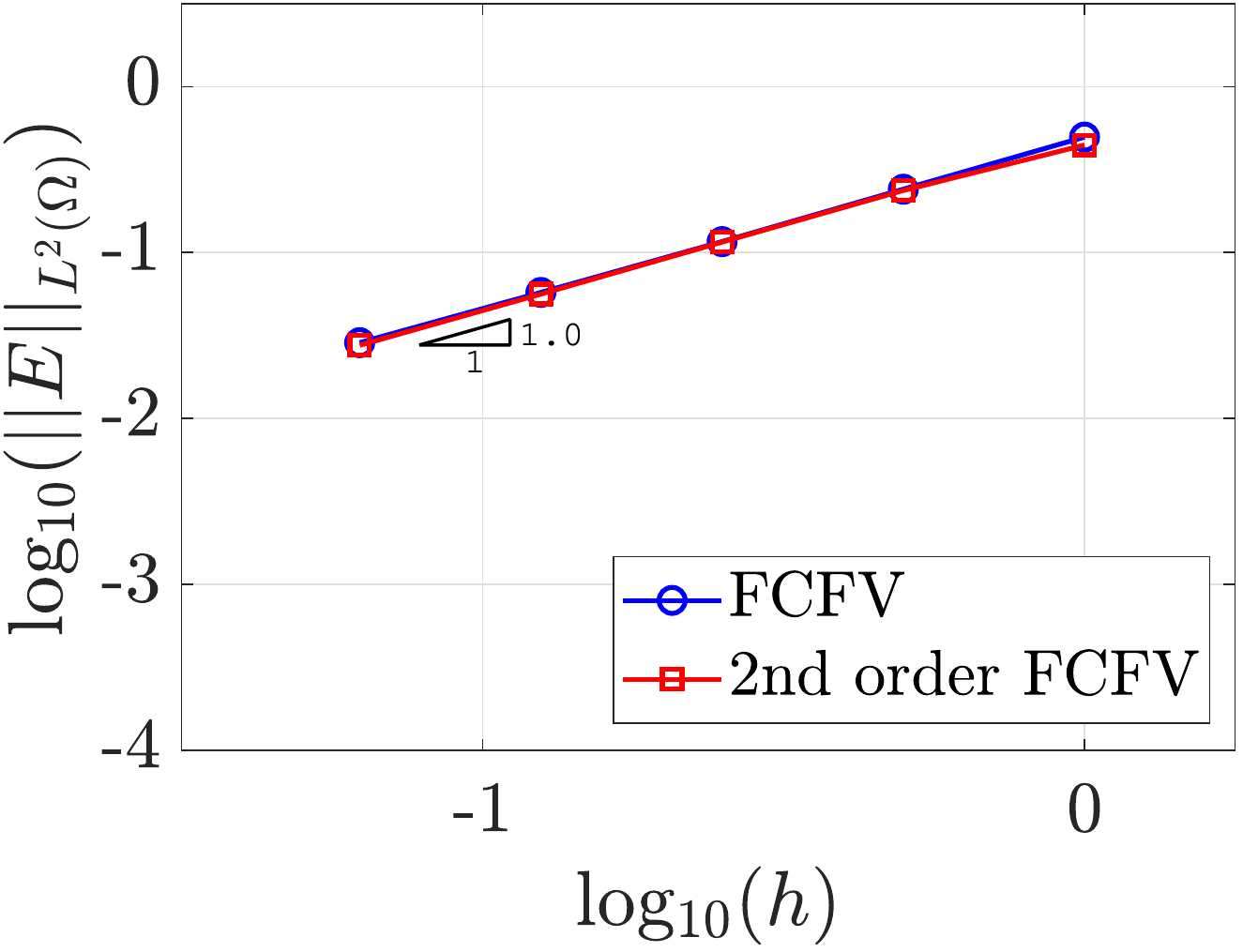}}
	\caption{Mesh convergence of the error of the velocity, its gradient and the pressure in the \eltwo($\Omega$) norm for three dimensional Stokes problem.}
	\label{fig:Stokes_Conv_3D}
\end{figure}

A mesh convergence is next performed for the Stokes problem in two and three dimensions. Figures~\ref{fig:Stokes_Conv_2D} and \ref{fig:Stokes_Conv_3D} show relative $\eltwo(\Omega)$ norm of the error of the velocity, its gradient and the pressure as a function of the characteristic mesh size $h$. 

The error of the velocity converges with second-order accuracy whereas first-order convergence is observed for its gradient, with an important gain on accuracy in the pressure. For the three dimensional test case, the optimal convergence is again observed for all the variables.

\subsection{Computational cost}
\label{sc:ComputationalCost}

The convergence studies performed in the previous section show an important gain in accuracy of the proposed second-order FCFV method when compared to the original FCFV in the same mesh. In this section both methods are compared in terms of the computational time. 

\begin{figure}[!tb]
	\centering	
	\subfigure[$\bu$]{\includegraphics[width=0.32\textwidth]{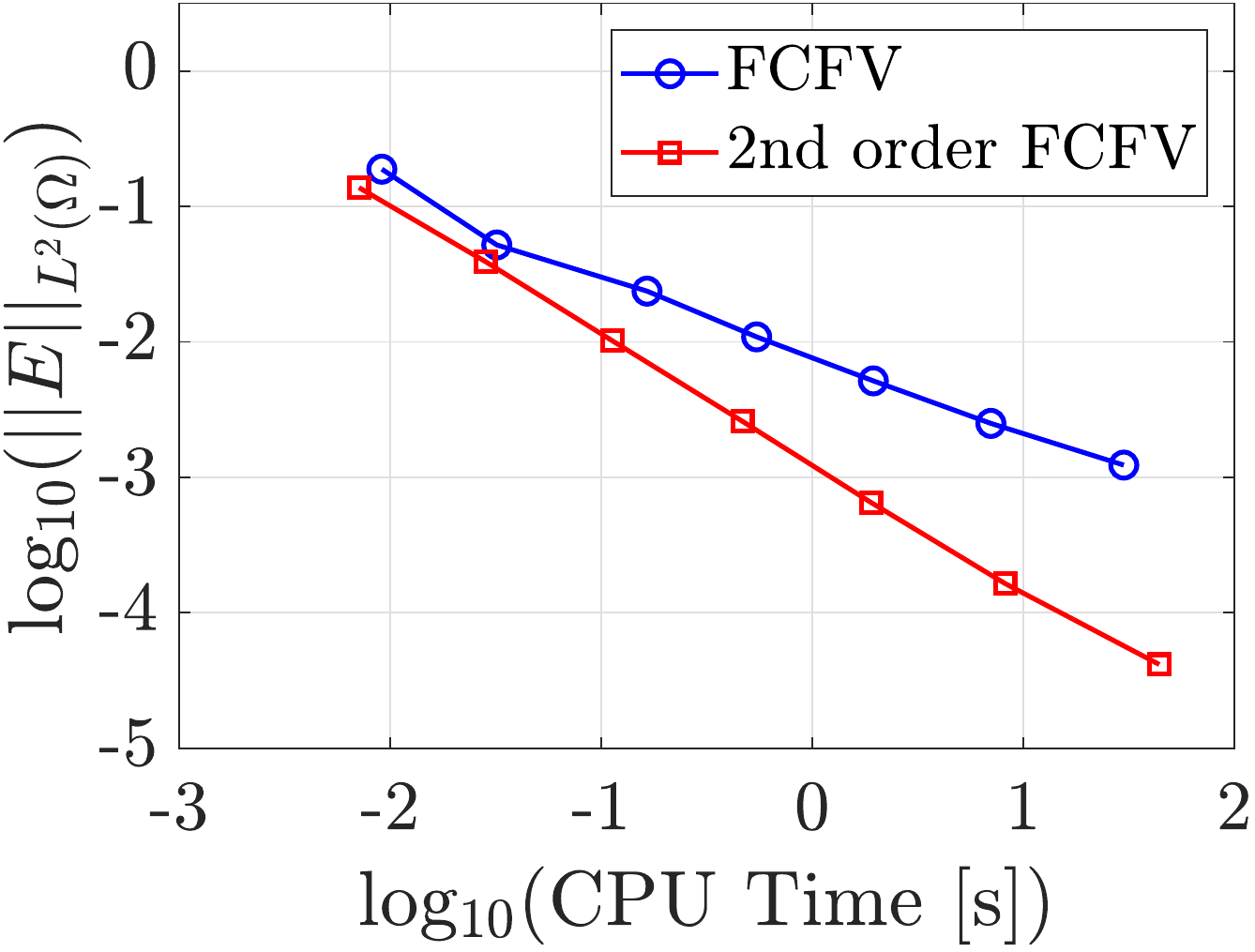}}
	\subfigure[$\bL$]{\includegraphics[width=0.32\textwidth]{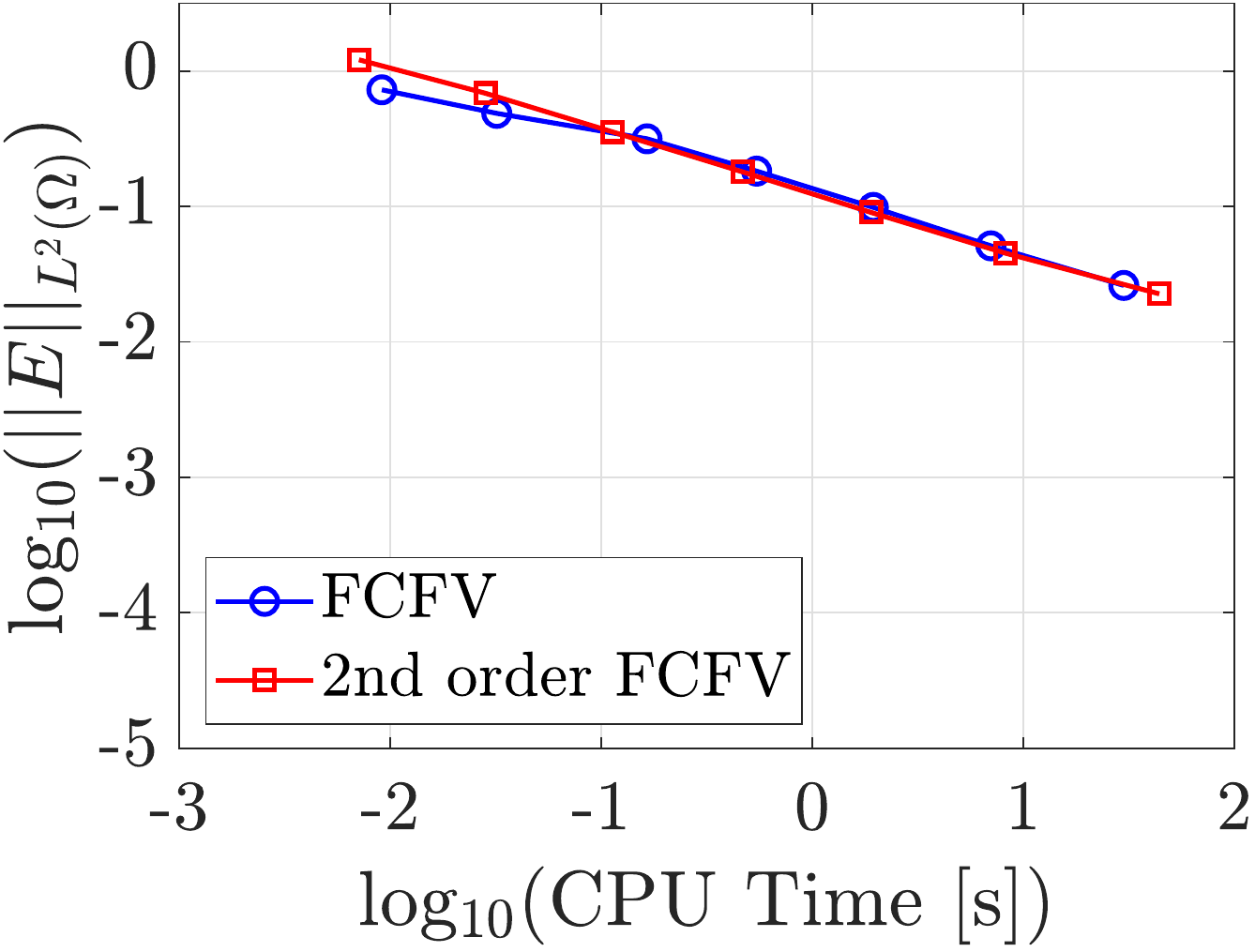}}
	\subfigure[$p$]{\includegraphics[width=0.32\textwidth]{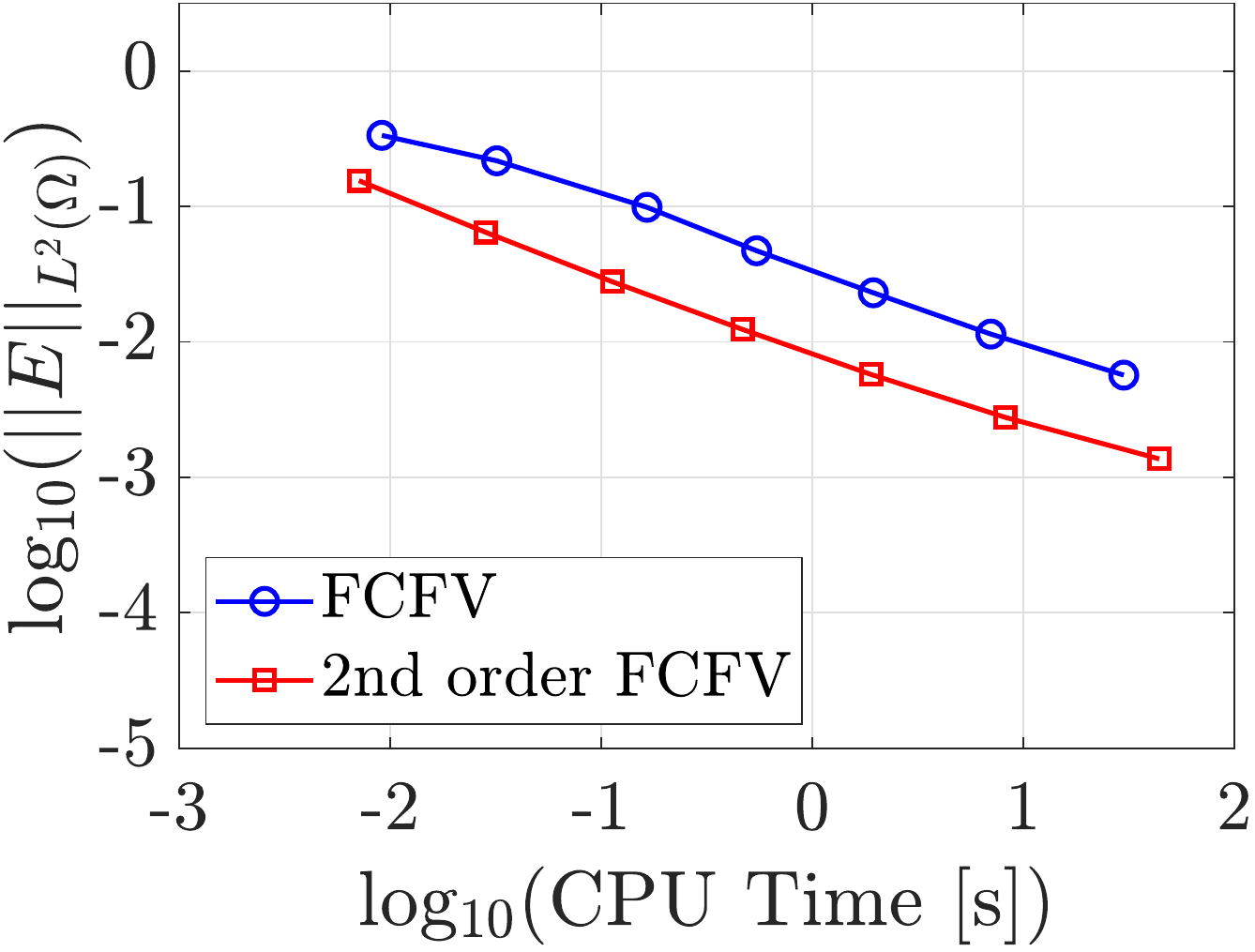}}
	\caption{Error of the velocity, its gradient and the pressure in \eltwo($\Omega$) norm as a function of the CPU time for two dimensional Stokes problem.}
	\label{fig:Stokes_Cpu_2D}
\end{figure}

\begin{figure}[!tb]
	\centering	
	\subfigure[$\bu$]{\includegraphics[width=0.32\textwidth]{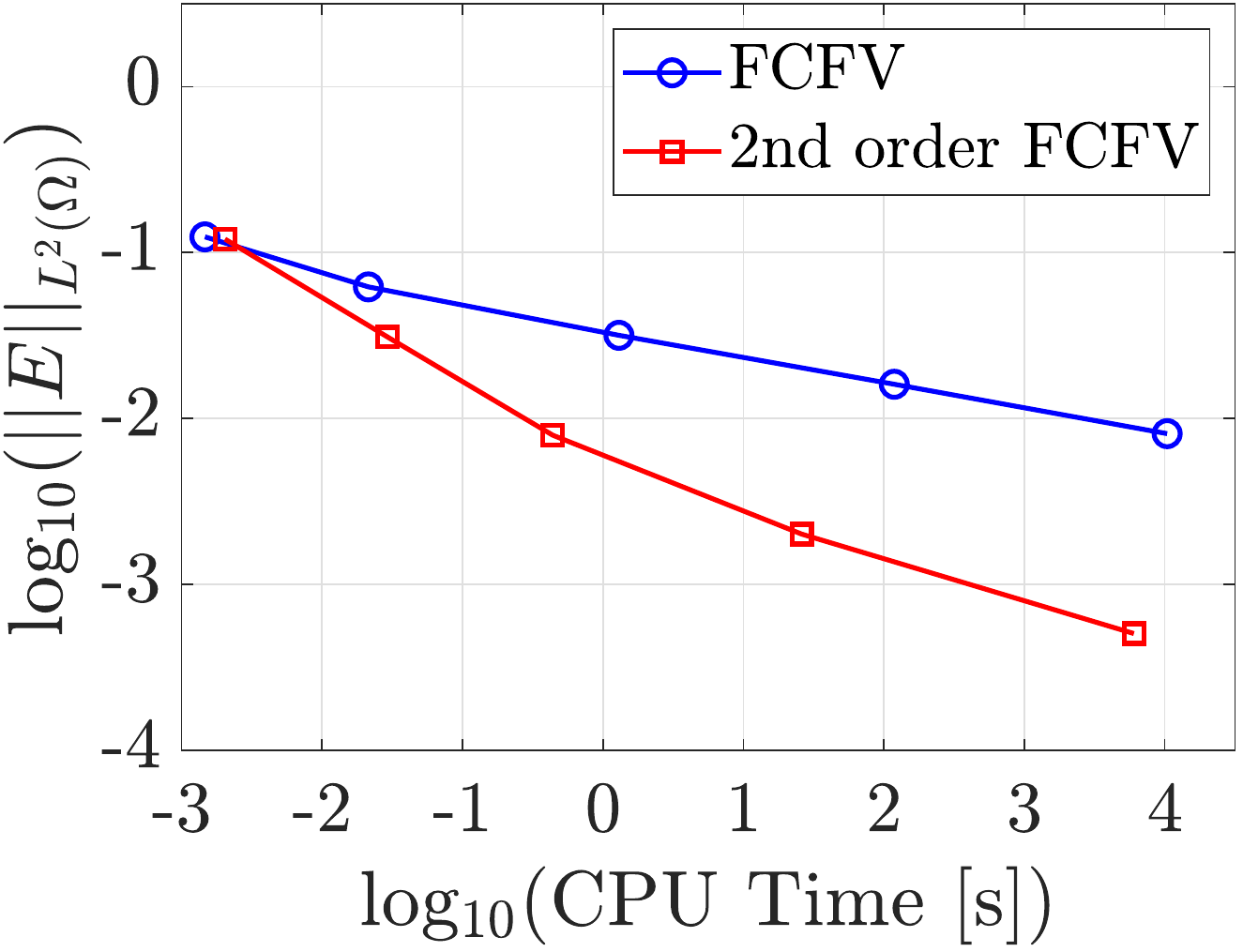}}
	\subfigure[$\bL$]{\includegraphics[width=0.32\textwidth]{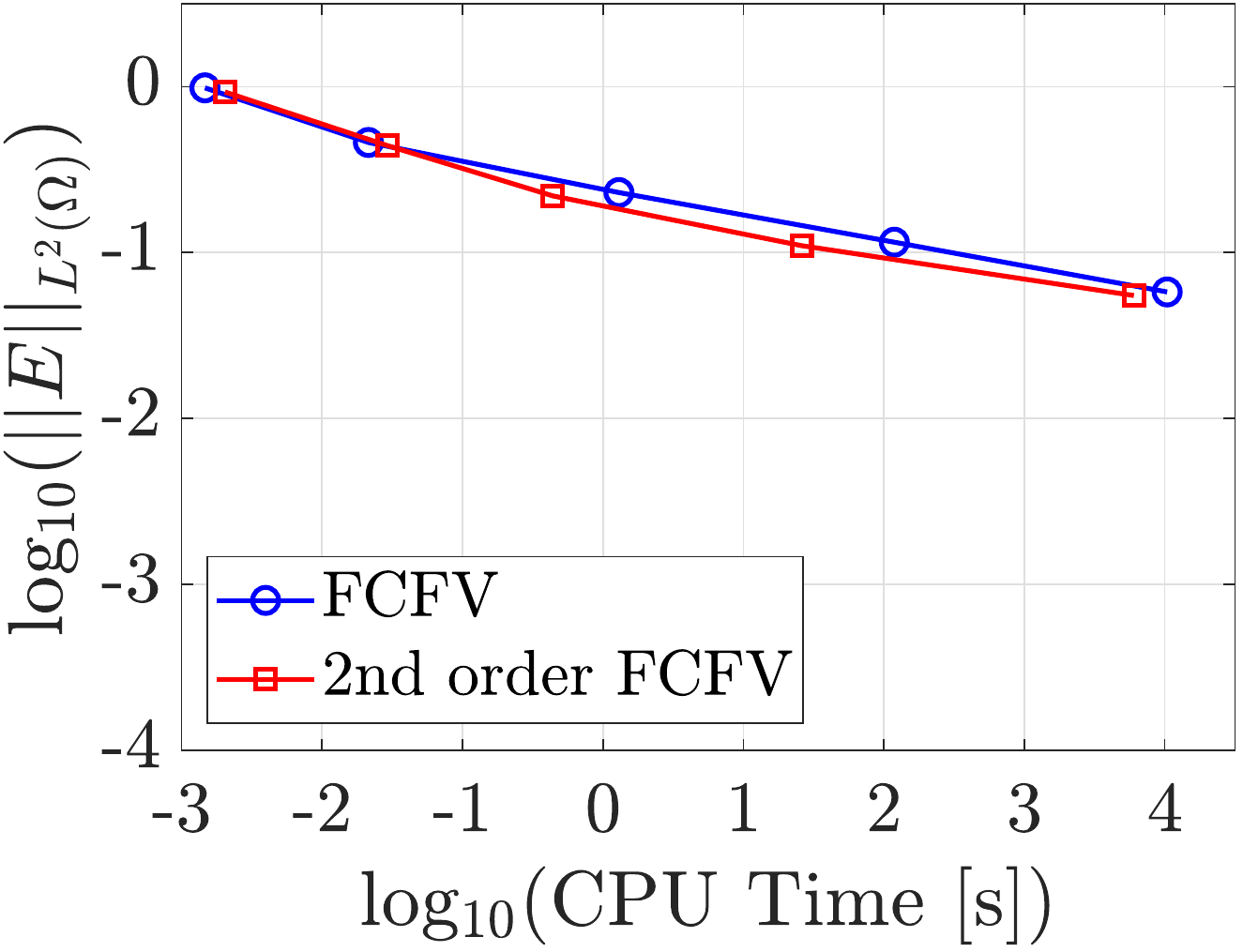}}
	\subfigure[$p$]{\includegraphics[width=0.32\textwidth]{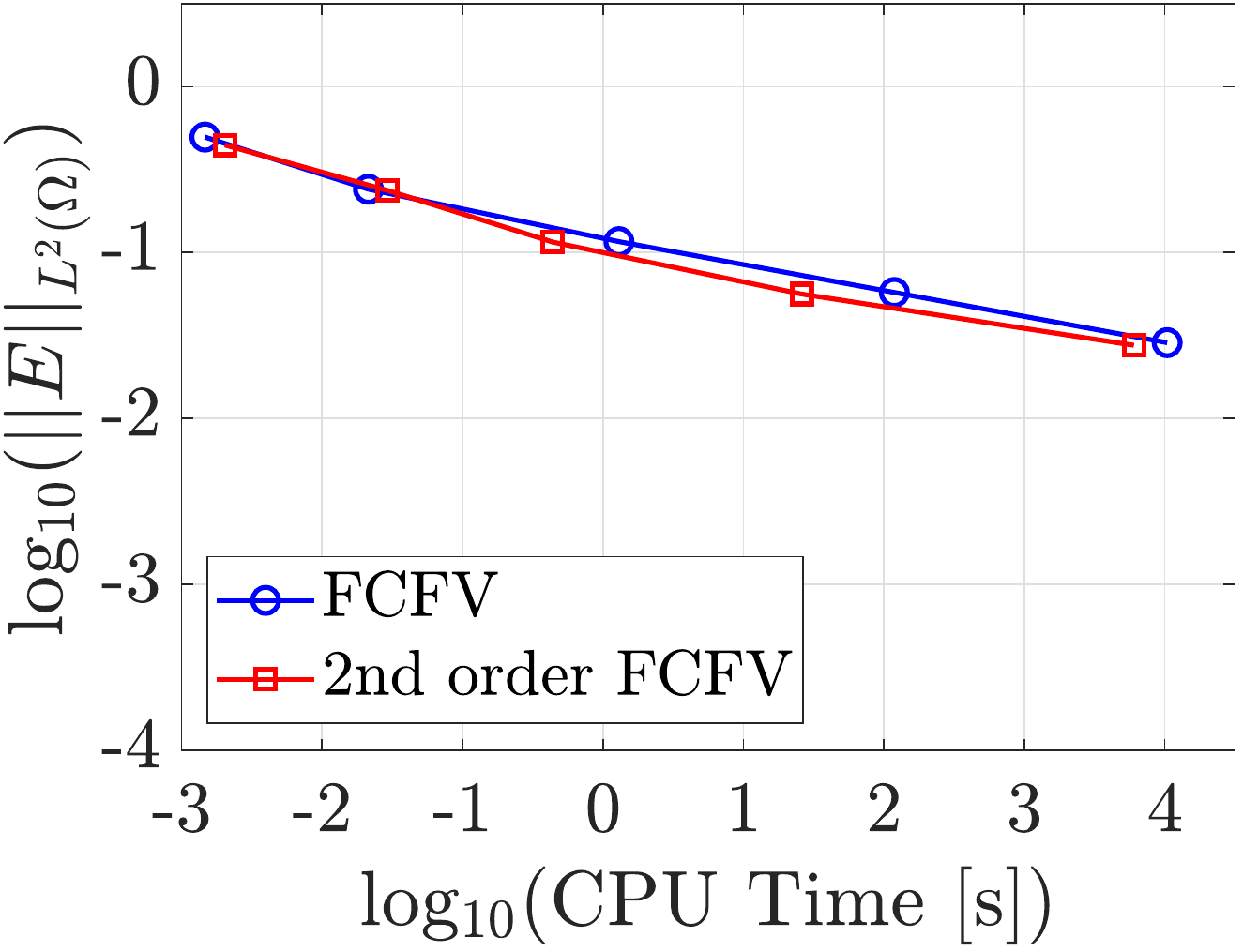}}
	\caption{Error of the velocity, its gradient and the pressure in \eltwo($\Omega$) norm as a function of the CPU time for three dimensional Stokes problem.}
	\label{fig:Stokes_Cpu_3D}
\end{figure}

Figure~\ref{fig:Stokes_Cpu_2D} shows relative $\eltwo(\Omega)$ norm of the error of the velocity, its gradient and the pressure as a function of the CPU time for the two dimensional Stokes problem. 
The results show that the proposed second-order FCFV is able to produce more accurate results for the velocity and the pressure using the same CPU time when compared to the first-order FCFV, whereas similar results are obtained for the veloocity gradient. In three dimensions similar conclusions are obtained, as shown in Figure~\ref{fig:Stokes_Cpu_3D}. 
The results show that the proposed second-order FCFV provides the same accuracy as the original first-order FCFV with orders of magnitude less CPU time when the velocity is of interest. For instance, an error in the velocity field of the order of 1\% is obtained in less than 1 second with the second-order FCFV whereas the first-order FCFV requires 2.7 hours.

The results for the Poisson problem, not displayed here for brevity, show the same advantages for the proposed second-order FCFV.

\subsection{Influence of the stabilisation parameter}
\label{sc:InfluenceTau}

\begin{figure}[!tb]
	\centering
	\subfigure[2D]{\includegraphics[width=0.48\textwidth]{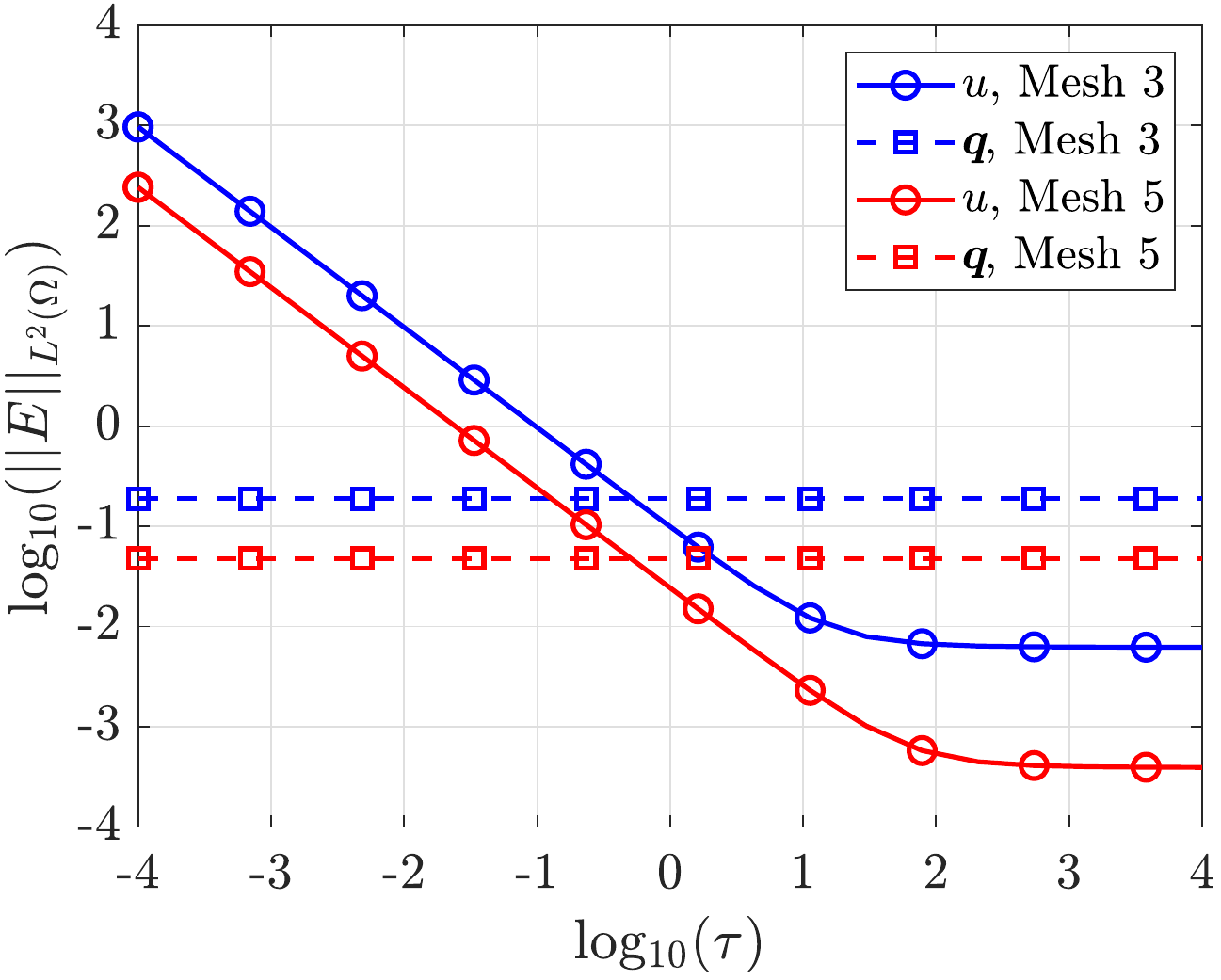}}
	\subfigure[3D]{\includegraphics[width=0.48\textwidth]{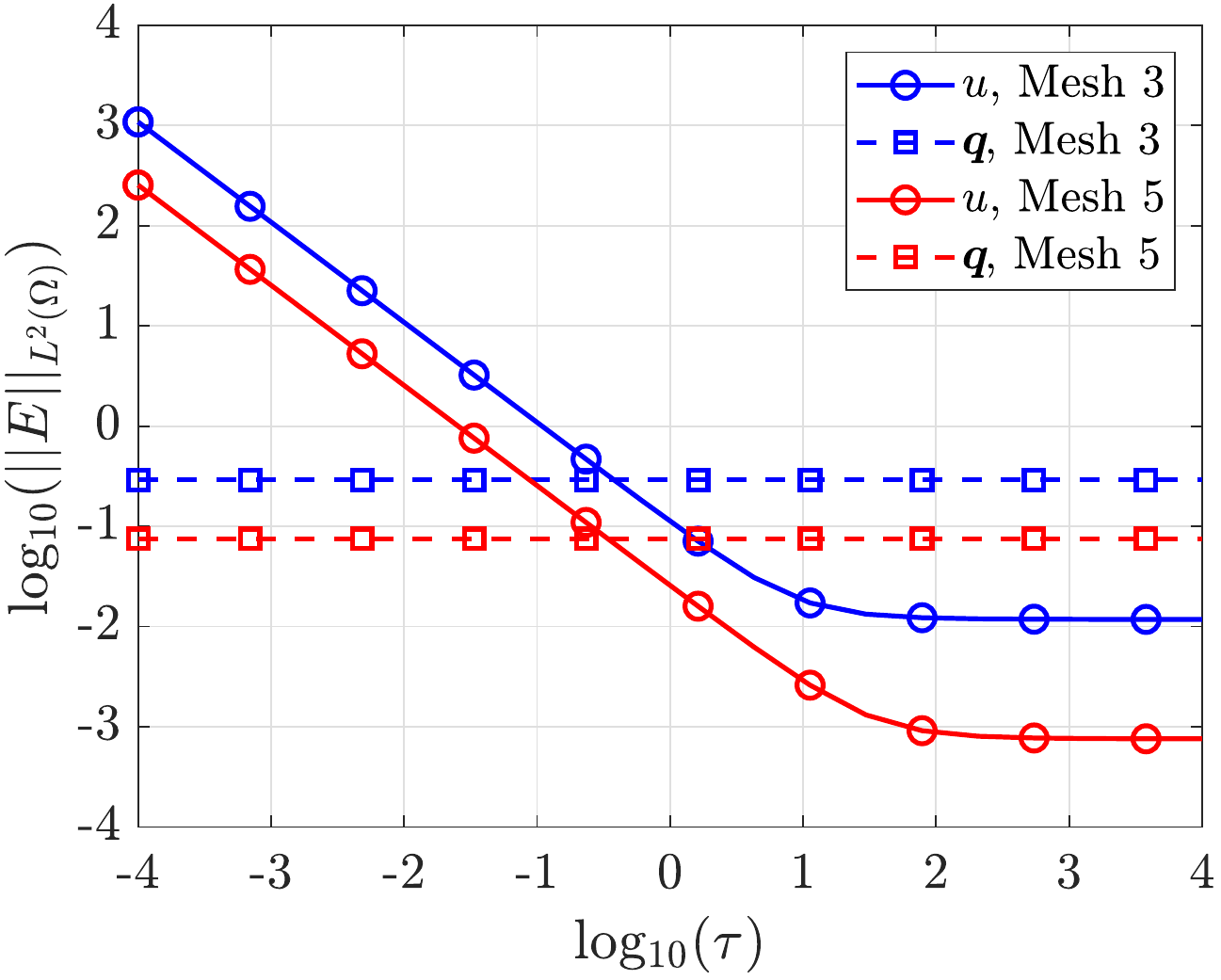}}
	\caption{Error of the solution and its gradient in \eltwo($\Omega$) norm as a function of stabilisation parameter $\tau$ for Poisson problem.}
	\label{fig:Poisson_Tau}
\end{figure}

The next study considers the influence of the stabilisation parameter $\tau$ in the accuracy of the proposed second-order FCFV method. Figure~\ref{fig:Poisson_Tau} shows the relative error, measured in the $\eltwo(\Omega)$ norm, as a function of the stabilisation parameter for the two and three dimensional Poisson problem and for two levels of mesh refinement. 
The results show that a maximum accuracy in the solution is achieved for a value of the stabilisation parameter of $10^2$ or larger, whereas the error of the solution gradient seems insensitive to the choice of this parameter. The behaviour is almost identical in both two and three dimensions. 

\begin{figure}[!tb]
	\centering
	\subfigure[2D]{\includegraphics[width=0.48\textwidth]{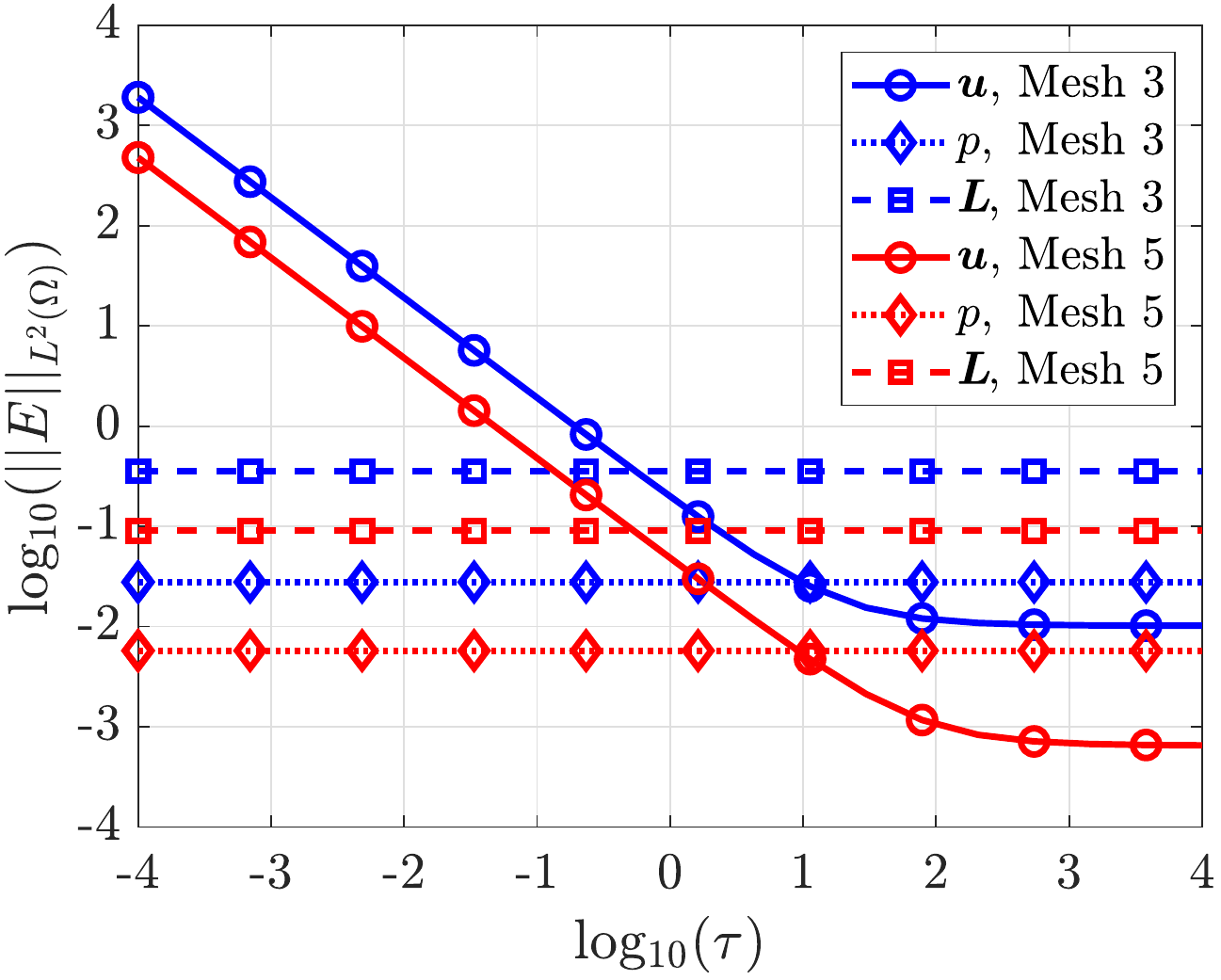}}
	\subfigure[3D]{\includegraphics[width=0.48\textwidth]{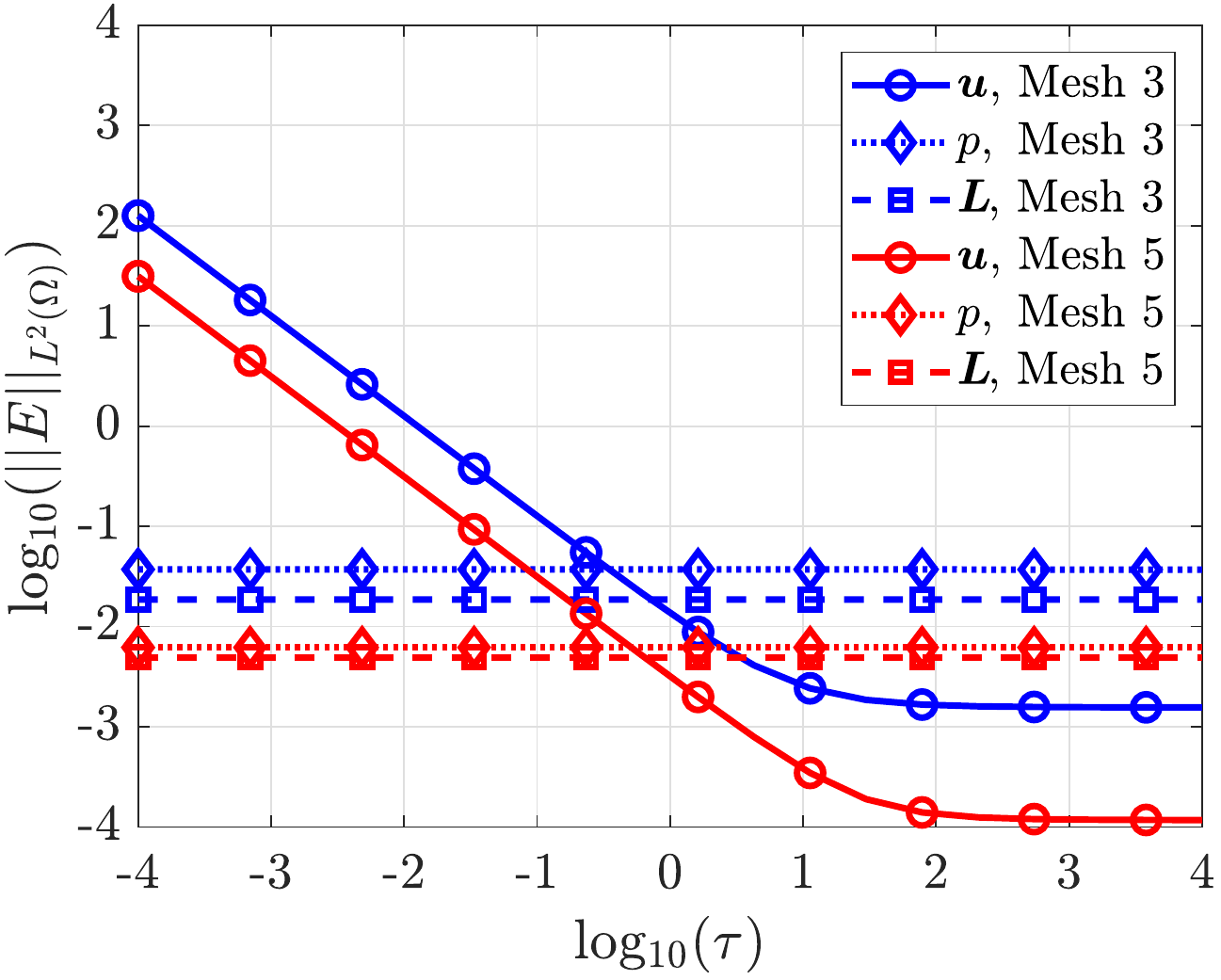}}
	\caption{Error of the velocity, its gradient and the pressure in \eltwo($\Omega$) norm as a function of stabilisation parameter $\tau$ for Stokes problem.}
	\label{fig:Stokes_Tau}
\end{figure}

For the Stokes problem, similar conclusions are obtained, as illustrated in Figure~\ref{fig:Stokes_Tau}.
In this case, the velocity gradient and the pressure are insensitive to the stabilisation parameter $\tau$, whereas the velocity requires a value of $10^2$ or larger to provide the maximum accuracy. 

It is worth emphasising that the value required to achieve the maximum accuracy of the solution is the same for two and three dimensional problems and for Poisson and Stokes problems. In addition, compared to the results presented in~\cite{FCFV2018} for the first-order FCFV, the proposed second-order FCFV is less sensitive to a particular choice of the stabilisation parameter.

\subsection{Influence of the cell distortion and stretching}
\label{sc:InfluenceDistortion}

\begin{figure}[!tb]
	\centering
	\subfigure[Distorted cells]{\includegraphics[width=0.32\textwidth]{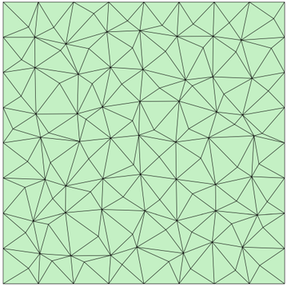}}
	\subfigure[Stretched cells]{\includegraphics[width=0.32\textwidth]{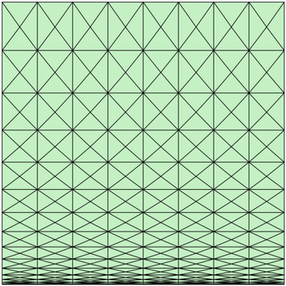}}
	\caption{Mesh with (a) distorted and (b) stretched cells to test the sensitivity of the FCFV to mesh quality.}
	\label{fig:mesh2D_TRI_H3_PerturbStretch}
\end{figure}

The last study considers the solution of Poisson and Stokes problems in meshes involving distorted and stretched cells. To illustrate the type of cells tested, Figure~\ref{fig:mesh2D_TRI_H3_PerturbStretch} shows the mesh corresponding to the third level of refinement where the cells have been randomnly disorted, as explained in~\cite{FCFV2018} and with stretched cells with a stretching factor of 100.

\begin{figure}[!tb]
	\centering
	\subfigure[2D]{\includegraphics[width=0.48\textwidth]{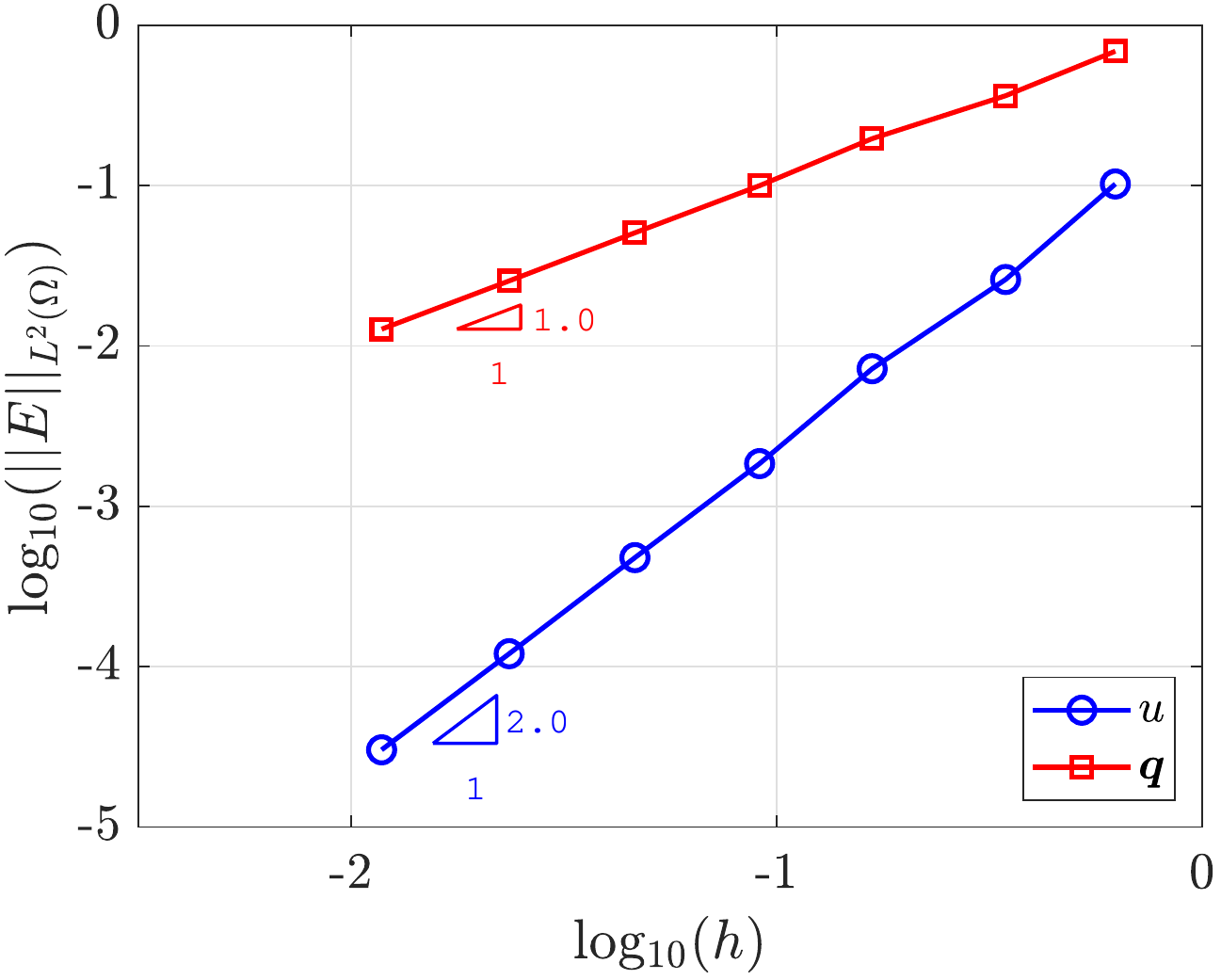}}
	\subfigure[3D]{\includegraphics[width=0.48\textwidth]{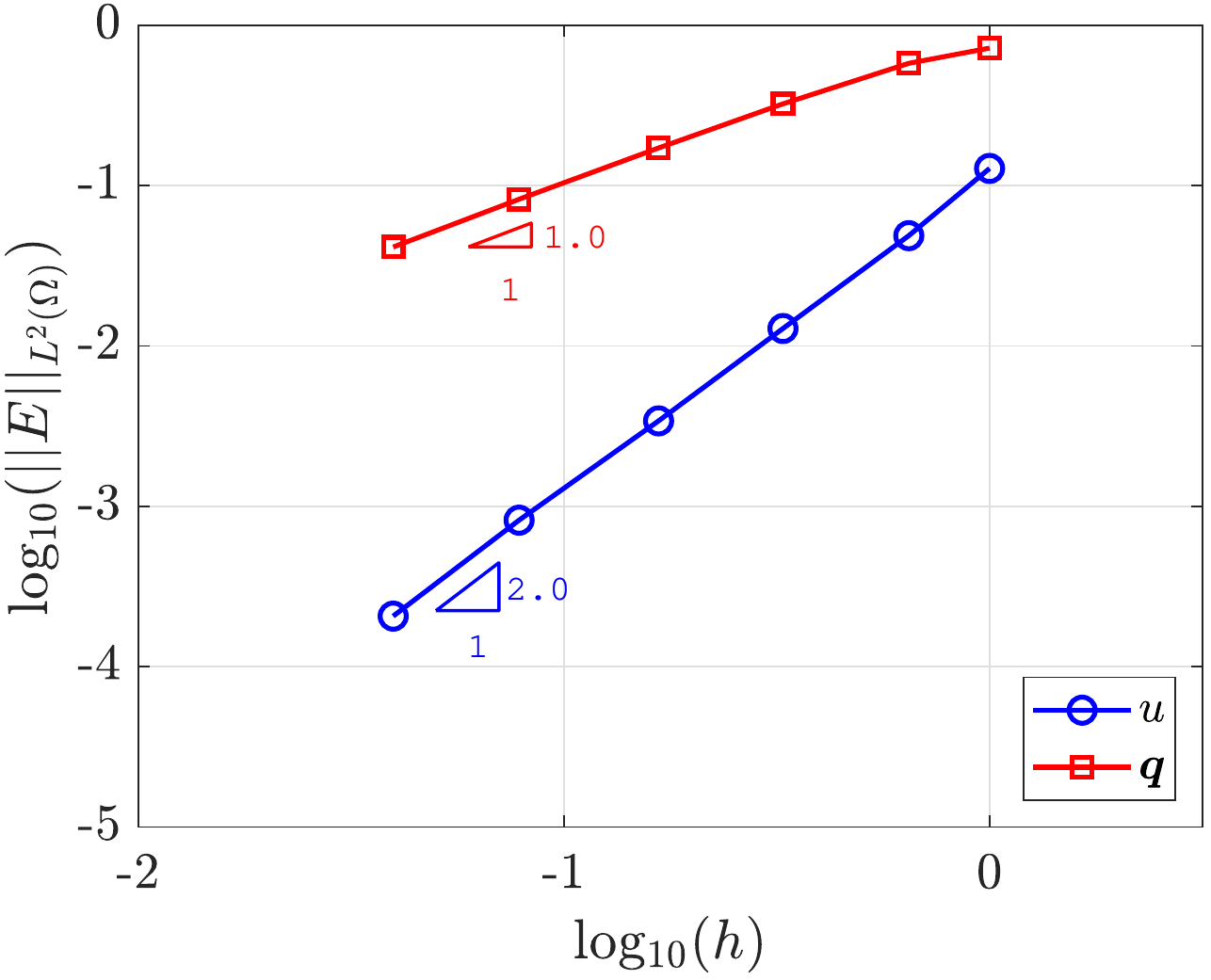}}
	\caption{Mesh convergence of the error of the solution and its gradient in \eltwo($\Omega$) norm for 2D and 3D Poisson problem with irregular mesh.}
	\label{fig:Poisson_Conv_Pert}
\end{figure}

Figure~\ref{fig:Poisson_Conv_Pert} shows a mesh convergence study for the Poisson problem in two and three dimensional meshes that have been distorted by randomly moving the interior nodes (i.e. the nodes not on the boundary of the domain) as described in~\cite{FCFV2018}. 
The results are almost identical to the ones obtained for regular meshes and displayed in Figures~\ref{fig:Poisson_Conv_2D} and \ref{fig:Poisson_Conv_3D}, showing that the proposed method is insensitive to mesh distortion. 

\begin{figure}[!tb]
	\centering
	\subfigure[2D]{\includegraphics[width=0.48\textwidth]{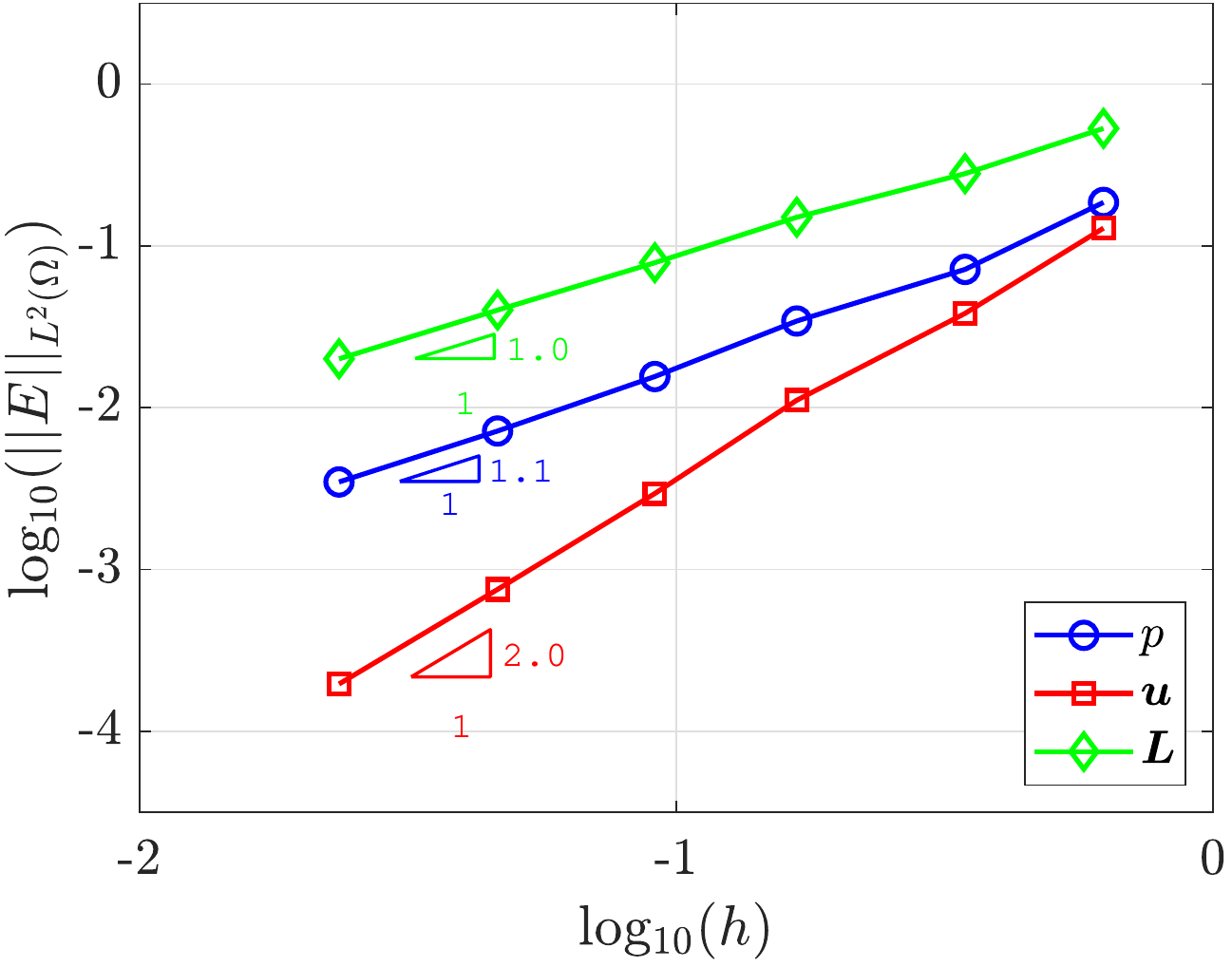}}
	\subfigure[3D]{\includegraphics[width=0.48\textwidth]{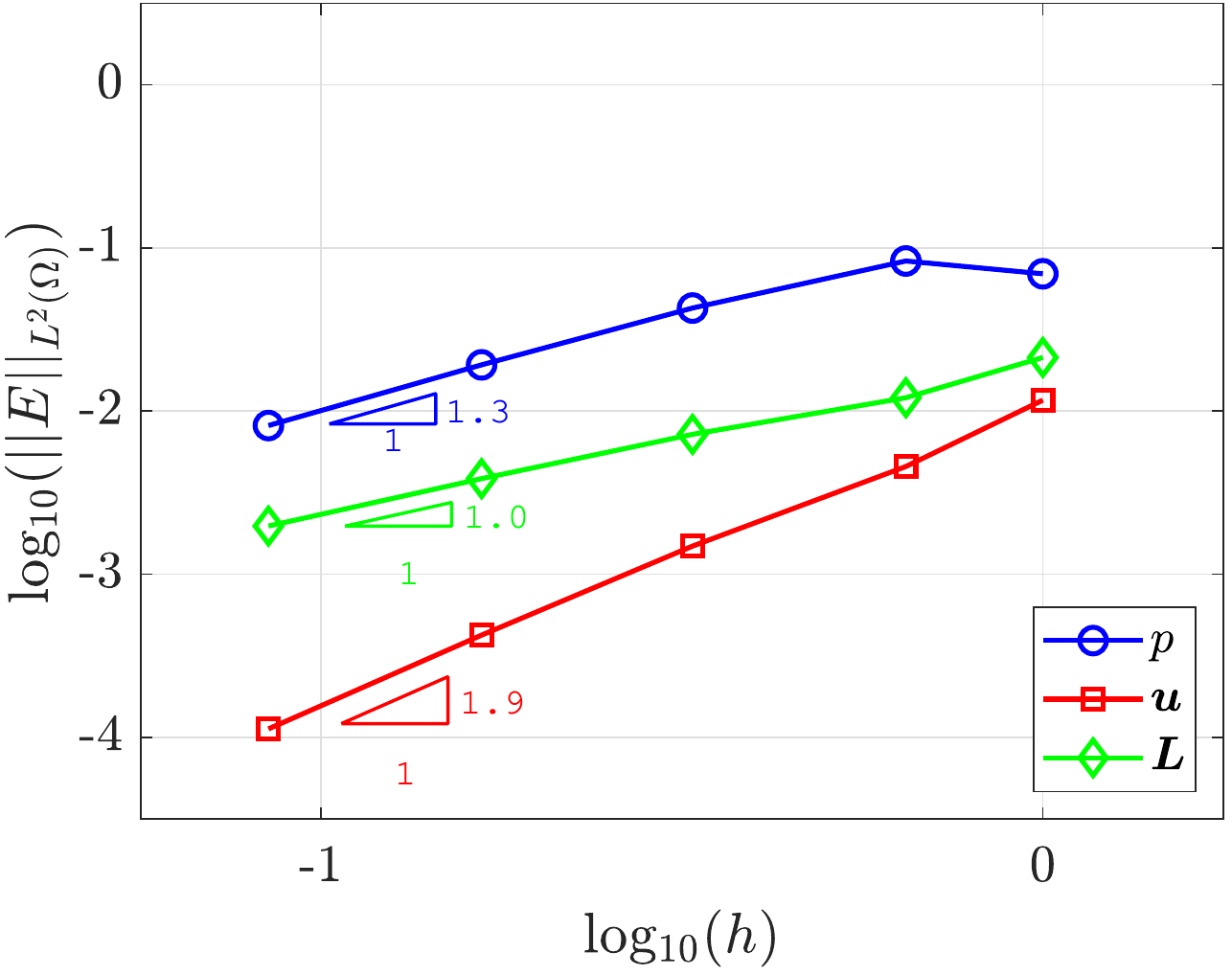}}
	\caption{Mesh convergence of the error of the velocity, its gradient and the pressure in \eltwo($\Omega$) norm for 2D and 3D Stokes problem with irregular mesh.}
	\label{fig:Stokes_Conv_Pert}
\end{figure}

The same conclusions are also obtained for the Stokes problem in two and three dimensions as shown in Figure~\ref{fig:Stokes_Conv_Pert}.

\begin{figure}[!tb]
	\centering
	\subfigure[2D]{\includegraphics[width=0.48\textwidth]{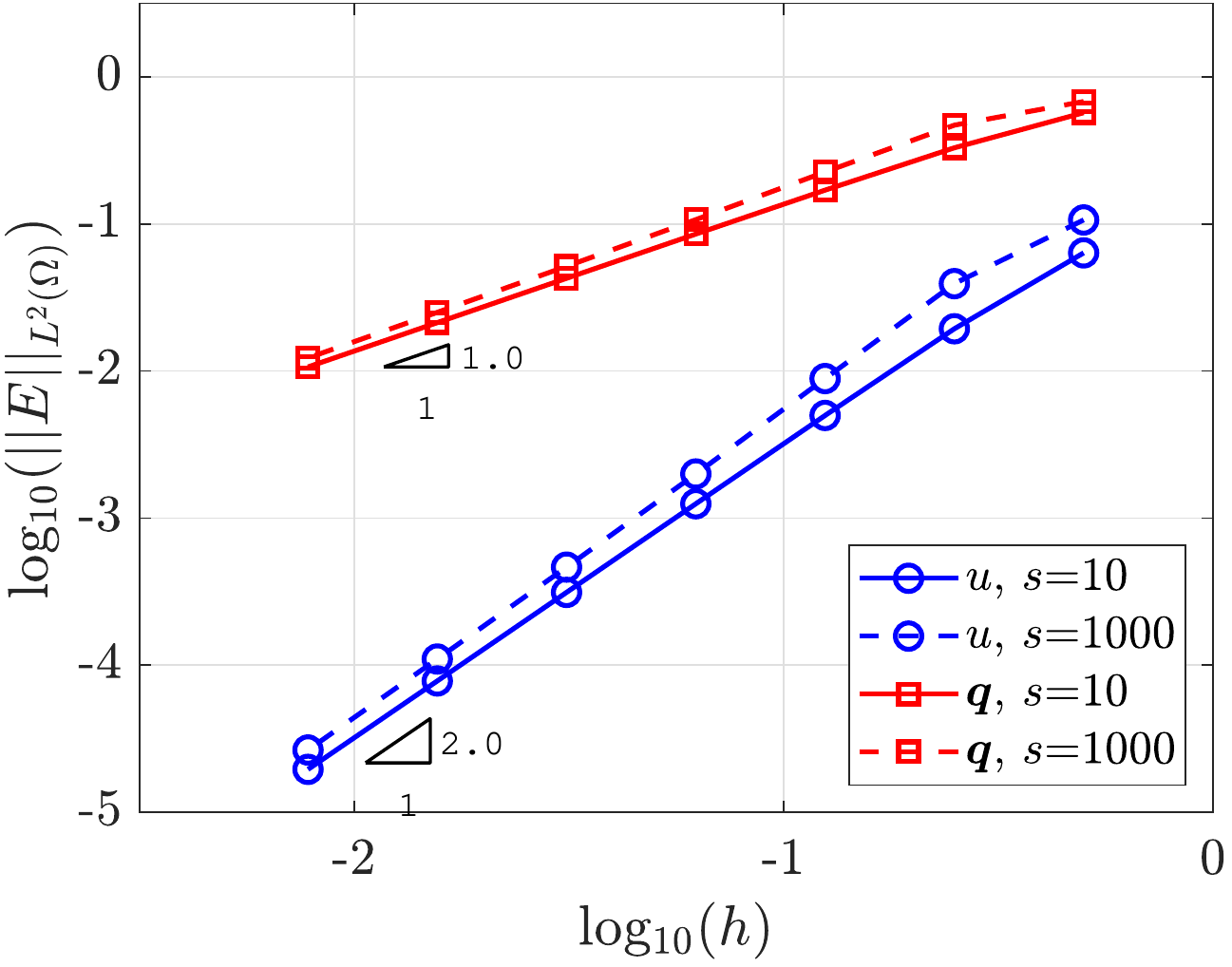}}
	\subfigure[3D]{\includegraphics[width=0.48\textwidth]{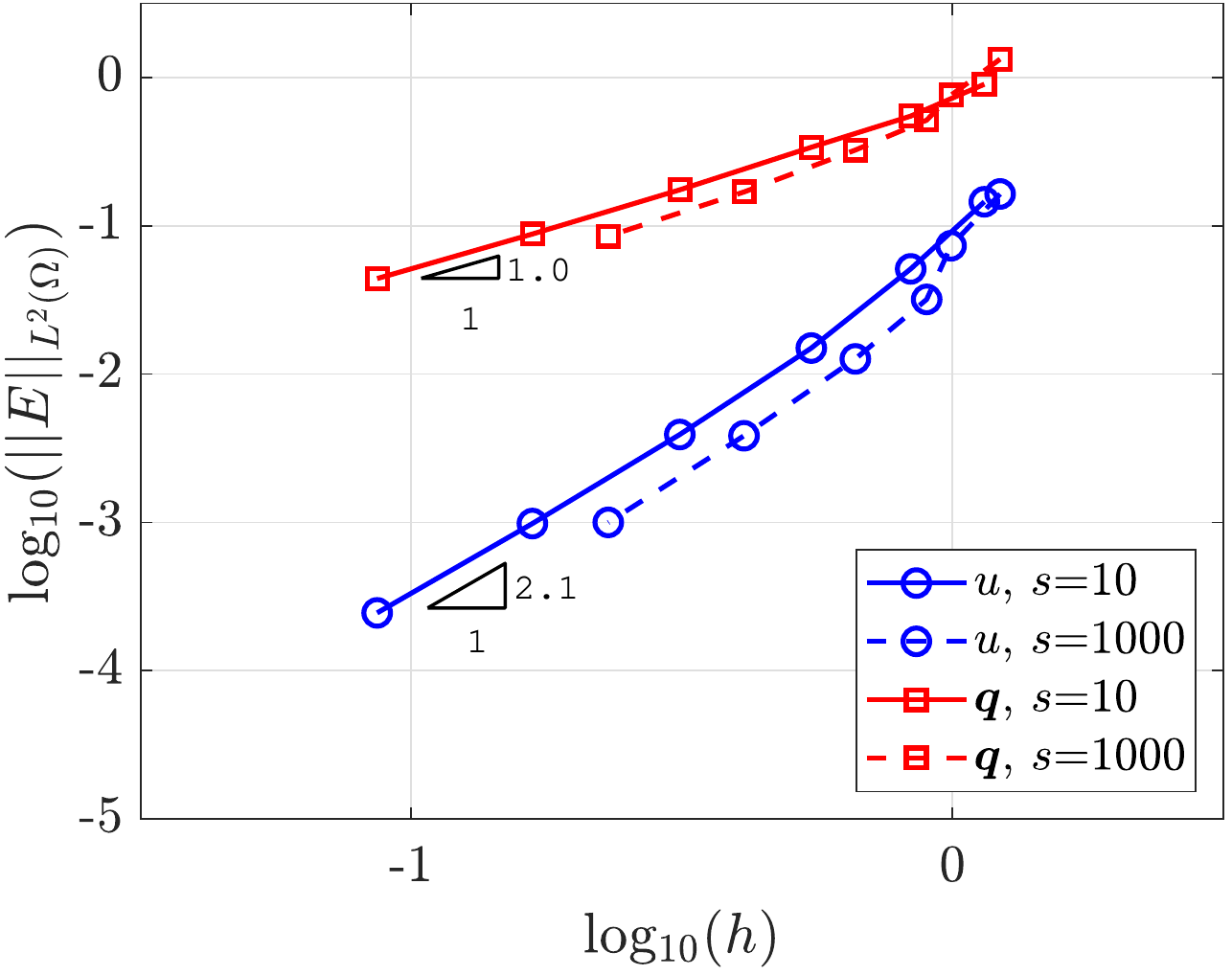}}
	\caption{Mesh convergence of the error of the solution and its gradient in \eltwo($\Omega$) norm for 2D and 3D Poisson problem with stretched meshes with stretching factor $s=10$ and $s=1,000$.}
	\label{fig:Poisson_Conv_Stret}
\end{figure}

Next, the influence of the cell stretching is considered. Two cases involving a maximum stretching factor, $s$, of 10 and 1,000 are considered, where this factor is measured in each cell as the ratio between the longest and the shortest edges. Figure~\ref{fig:Poisson_Conv_Stret} shows the relative $\eltwo(\Omega)$ norm of the error of the solution and its gradient as a function of the characteristic mesh size for the Poisson problem solved in stretched meshes in two and three dimensions.
Almost identical results are observed for both stretching factors. In the two dimensional problem a marginal lower error is observed for the mesh with stretching factor $s=10$ whereas for the three dimensional problem the mesh with stretching factor $s=1,000$ produces a slightly lower error. 

\begin{figure}[!tb]
	\centering
	\subfigure[2D]{\includegraphics[width=0.48\textwidth]{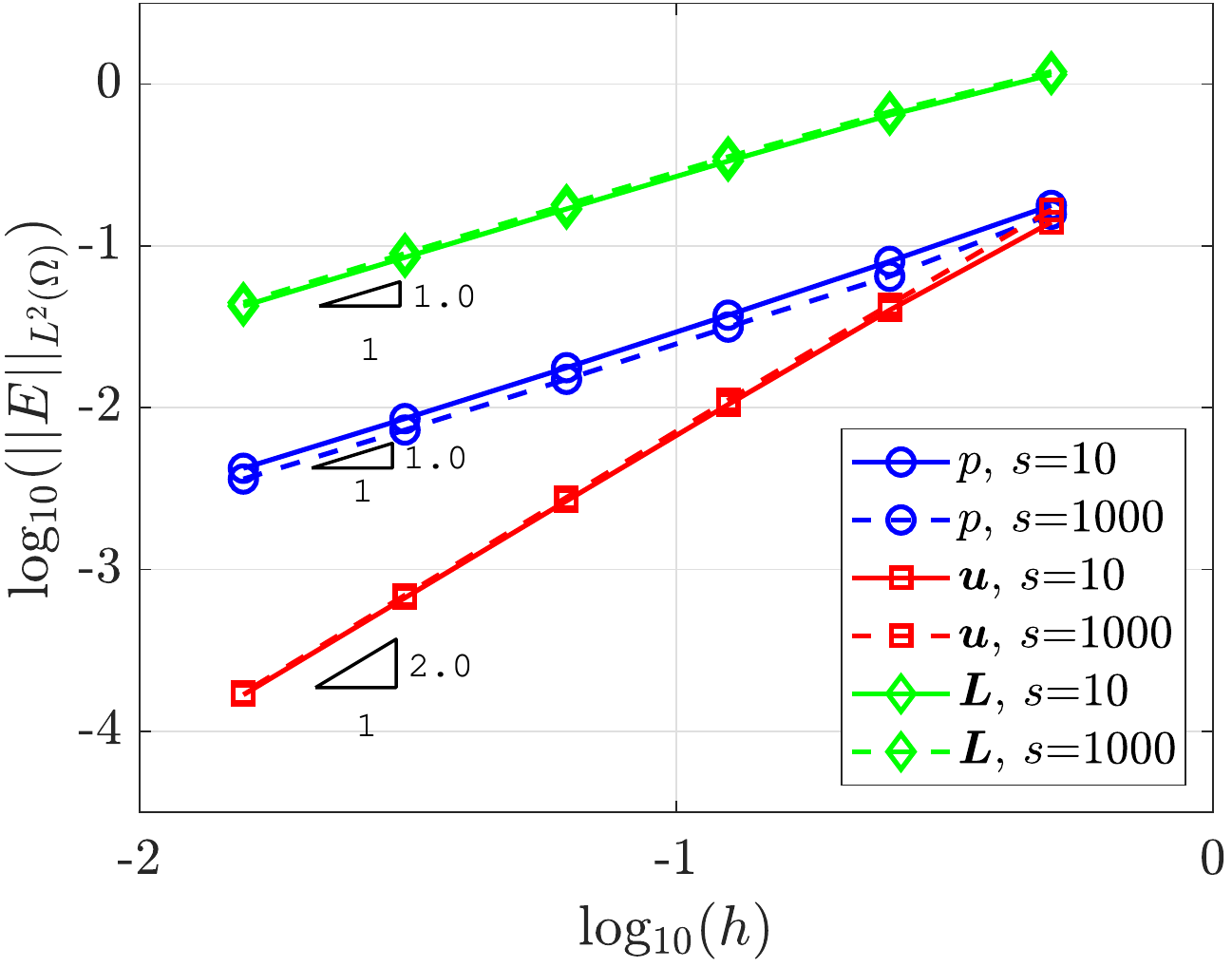}}
	\subfigure[3D]{\includegraphics[width=0.48\textwidth]{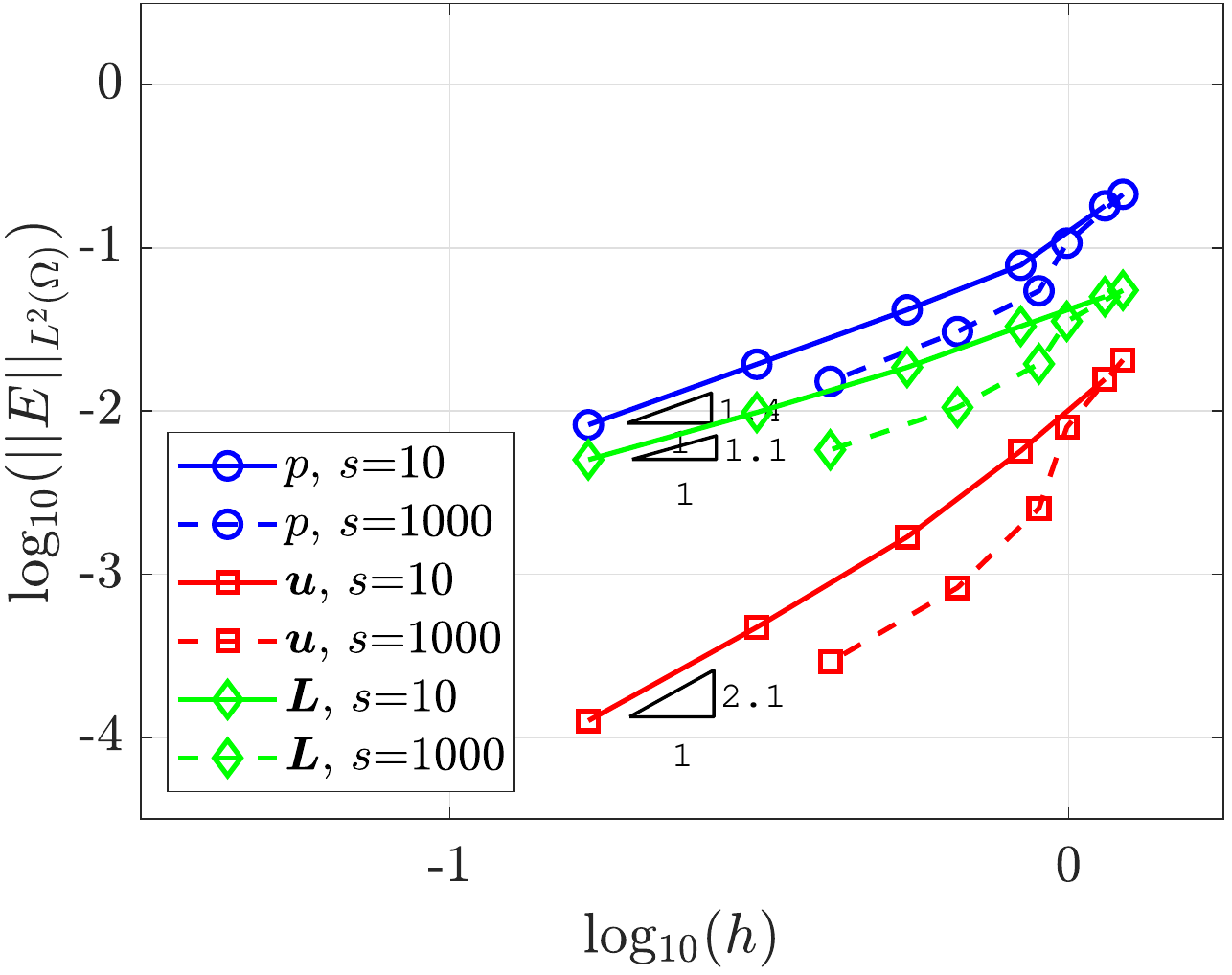}}
	\caption{Mesh convergence of the error of the velocity, its gradient and the pressure in \eltwo($\Omega$) norm for 2D and 3D Stokes problem with stretched meshes with stretching factor $s=10$ and $s=1,000$.}
	\label{fig:Stokes_Conv_Stret}
\end{figure}

The results for the Stokes problem in two and three dimensions are displayed in Figure~\ref{fig:Stokes_Conv_Stret}, confirming the conclusions observed for the Poisson problem.

It is worth noting that contrary to other FV methods, the proposed second-order FCFV method not only shows an accuracy that is insensitive to mesh distortion and stretching but also preserves the optimal rate of convergence in all the variables, i.e. the solution and its gradient for the Poisson problem and the velocity, its gradient and the pressure for the Stokes problem.

\section{Numerical examples}
\label{sc:examples}

\subsection{Irrotational flow past a full aircraft}
\label{sc:PoissonBoeing} 

\begin{figure}[!tb]
	\centering
	\subfigure[Velocity]{\includegraphics[width=0.48\textwidth]{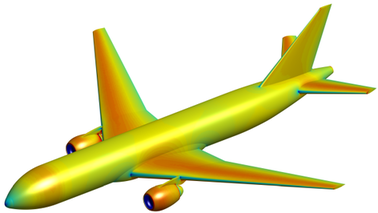}}
	\subfigure[Pressure]{\includegraphics[width=0.48\textwidth]{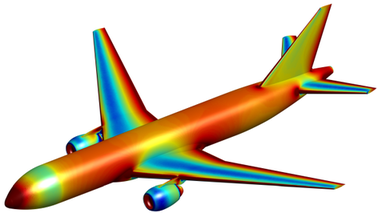}}
	\caption{Magnitude of the velocity and pressure distribution for the irrotational flow around a full aircraft configuration.}
	\label{fig:boeingK0K1}
\end{figure}

To show the ability of the proposed method to efficiently solve large scale problems involving complex geometries, the irrotational flow around a full aircraft is considered. A tetrahedral mesh with 5,125,998 cells is considered, leading to a global system of 11,283,113 equations to find the solution on the cell faces. The magnitude of the velocity, computed from the gradient of the solution, and the pressure, computed from the Bernoulli equation are displayed in Figure~\ref{fig:boeingK0K1}.

The solution using the proposed second-order FCFV took 5.1 minutes for the computation of all the elemental matrices and the assembly and 6.4 minutes for solving the global system of equations using a direct method. The developed code is written in Matlab and the computation was performed in an Intel$^{\tiny{\textregistered}}$
Xeon$^{\tiny{\textregistered}}$ CPU $@$ 3.70GHz and 32GB main memory available.

It is worth noting that the time recorded for assembling the system is slightly higher than the times reported in~\cite{FCFV2018} for the first-order FCFV due to the extra operations required by the second-order method for the computation and assembly of the global matrix. The time required by the proposed second-order FCFV for solving the system is almost identical to the time employed by the first-order method due to the global matrix having the exact same size and structure, with the same number of non-zero elements. Finally, it is worth emphasising that, as shown in the numerical experiments of the previous Section, the extra cost induced by the second-order FCFV for a given spatial discretisation leads to a substantial gain in accuracy and, in general, for a given error, the second-order FCFV is more efficient.

\subsection{Stokes flow past a sphere}
\label{sc:stokesSphere} 

The next example considers the Stokes flow around a sphere. This classical three dimensional example is used to compare the accuracy and performance of the proposed second-order FCFV against the original first-order FCFV method for a problem involving large three dimensional meshes. The domain of interest is $\Omega = \left( [-7,15] \times [-5,5] \times [-5,5] \right)\setminus \mathcal{B}_{1,\mathbf{0}}$, where $\mathcal{B}_{1,\mathbf{0}}$ is the unit ball with centre at the origin. Homogeneous Neumann boundary conditions are imposed on the outlet part of the boundary, corresponding to $x=15$, whereas Dirichlet boundary conditions, corresponding to the exact solution, are imposed on the rest of the boundary.

Six unstructured tetrahedral meshes are considered with 3,107, 10,680, 43,682, 204,099, 686,853 and 2,516,099 cells respectively. The size of the global system of equations to be solved to obtain the velocity on the cell faces is 20,711, 72,249, 299,276, 1,409,916, 4,765,776 and 17,513,075 respectively.

\begin{figure}[!tb]
	\centering
	\subfigure[$\eltwo(\Omega)$ error convergence]{\includegraphics[width=0.48\textwidth]{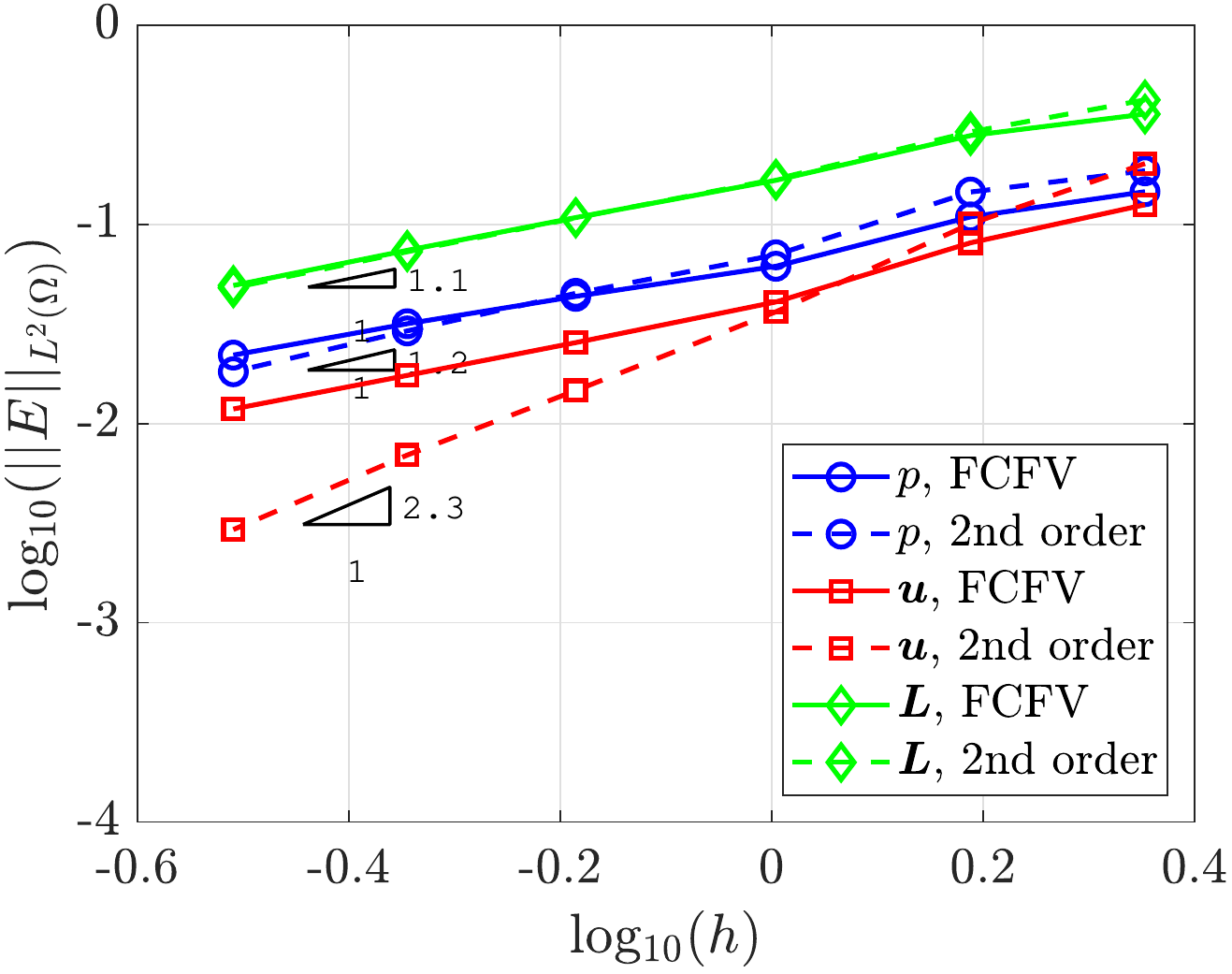}}
	\subfigure[Drag convergence]{\includegraphics[width=0.48\textwidth]{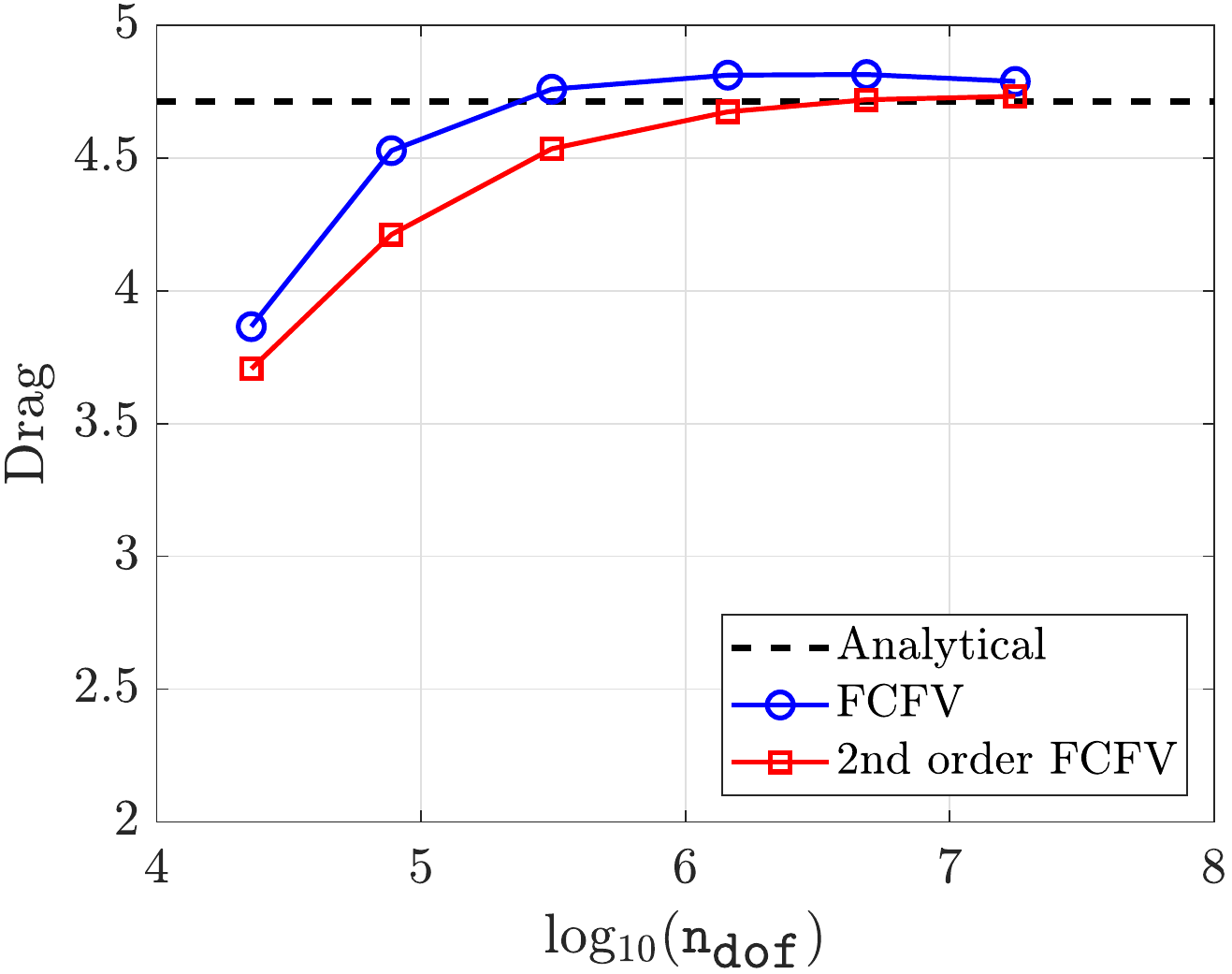}}
	\caption{Comparison between first and second-order FCFV for the Stokes flow past a sphere. (a) Convergence of the error of velocity, pressure and gradient of the velocity in the $\eltwo(\Omega)$ norm and (b) convergence of the drag.}
	\label{fig:stokes3D_Sphere}
\end{figure}

Figure~\ref{fig:stokes3D_Sphere}(a) shows the convergence of the error for the velocity, pressure and gradient of the velocity as the mesh is refined. The results are compared to the original first-order FCFV and clearly show the advantage of the proposed method by providing second-order convergence on the velocity. In this example the accuracy for the other variables is almost identical, with a marginal gain observed in the computation of the pressure in the finest meshes. 
Figure~\ref{fig:stokes3D_Sphere}(b) shows the convergence of the drag as the mesh is refined. The advantages of the proposed second-order FCFV are observed as the convergence to the exact value is faster than with the original first-order method. 

\subsection{Mesh adaptivity for the Poisson problem} \label{sc:adaptivityPoisson}

This Section presents a numerical example to illustrate the strategy described in Section~\ref{sc:adaptivity} to perform an automatic mesh adaptive process by combining the results of the first-order and second-order FCFV methods. A two dimensional Poisson problem with known analytical solution is considered in $\Omega = [0,1]^2$. The source term and Dirichlet boundary conditions are selected so that the analytical solution is given by
\begin{equation}
u(x_1,x_2) = \exp\Big(-100\left[ (x_1-0.7)^2 + (x_2-0.7)^2 \right] \Big)
\end{equation}
and the desired accuracy is $\varepsilon = 10^{-2}$.

\begin{figure}[!tb]
	\centering
	\subfigure[Mesh 1]{\includegraphics[width=0.32\textwidth]{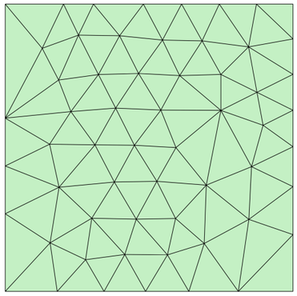}}
	\subfigure[Mesh 3]{\includegraphics[width=0.32\textwidth]{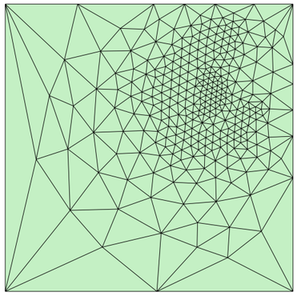}}
	\subfigure[Mesh 6]{\includegraphics[width=0.32\textwidth]{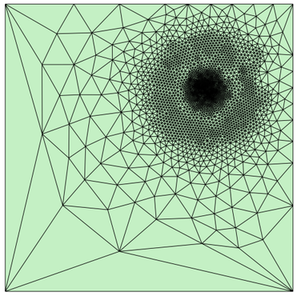}}
	\caption{Three meshes used in the automatic mesh adaptive process for the Poisson problem.}
	\label{fig:Poisson_Adap_Mesh}
\end{figure}

\begin{figure}[!tb]
	\centering
	\subfigure[Mesh 1]{\includegraphics[width=0.32\textwidth]{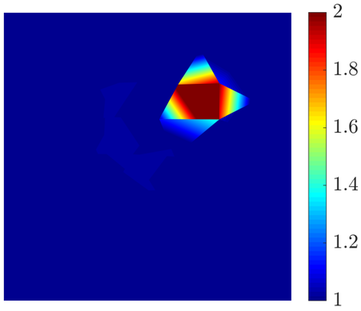}}
	\subfigure[Mesh 3]{\includegraphics[width=0.32\textwidth]{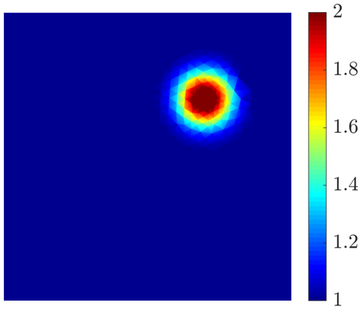}}
	\subfigure[Mesh 6]{\includegraphics[width=0.32\textwidth]{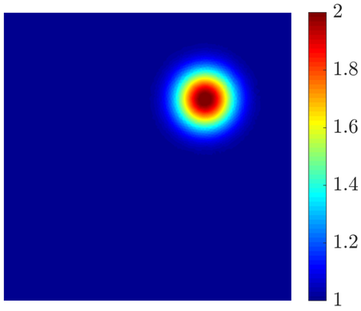}}
	\caption{Three FCFV approximations corresponding to the meshes of Figure~\ref{fig:Poisson_Adap_Mesh}.}
	\label{fig:Poisson_Adap_Sol}
\end{figure}

The process starts with the coarse mesh represented in Figure~\ref{fig:Poisson_Adap_Mesh}(a). The approximation with the proposed FCFV on the coarsest mesh is depicted in Figure~\ref{fig:Poisson_Adap_Sol}(a).
By comparing the approximation of the second-order FCFV with the approximation of the first-order FCFV method (i.e. without the projection operator), a desired cell size is computed, using Equation~\eqref{eq:newElemSize}, and a new mesh is generated. The mesh and the second-order FCFV approximation after two iterations of the mesh adaptive process are displayed in Figures~\ref{fig:Poisson_Adap_Mesh}(b) and \ref{fig:Poisson_Adap_Sol}(b) respectively. It can be clearly observed how the mesh is coarsened in the regions where the approximation is almost constant, whereas the mesh density is increased in the regions where the approximation changes rapidly. The adaptive process finishes in five iterations, when the desired error is achieved. The final mesh and FCFV approximation are represented in Figures~\ref{fig:Poisson_Adap_Mesh}(c) and \ref{fig:Poisson_Adap_Sol}(c) respectively.

\begin{figure}[!tb]
	\centering
	\subfigure[Exact error, \nit=0]{\includegraphics[width=0.32\textwidth]{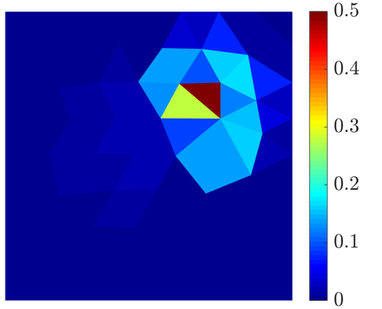}}
	\subfigure[Exact error, \nit=2]{\includegraphics[width=0.32\textwidth]{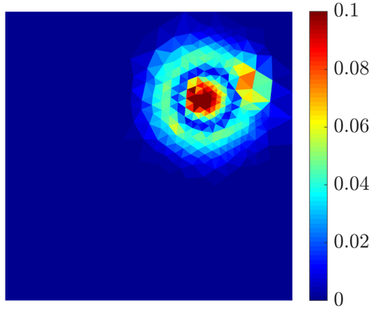}}
	\subfigure[Exact error, \nit=5]{\includegraphics[width=0.32\textwidth]{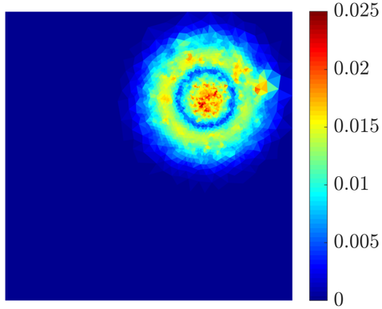}}
	\subfigure[Error indicator, \nit=0]{\includegraphics[width=0.32\textwidth]{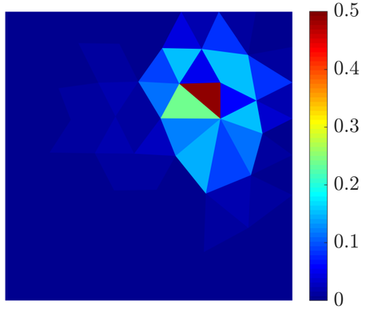}}	
	\subfigure[Error indicator, \nit=2]{\includegraphics[width=0.32\textwidth]{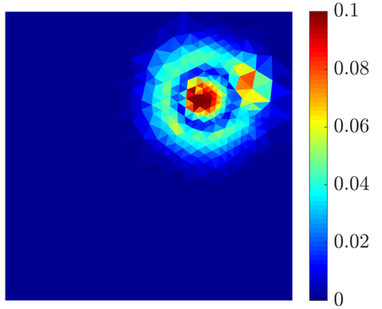}}
	\subfigure[Error indicator, \nit=5]{\includegraphics[width=0.32\textwidth]{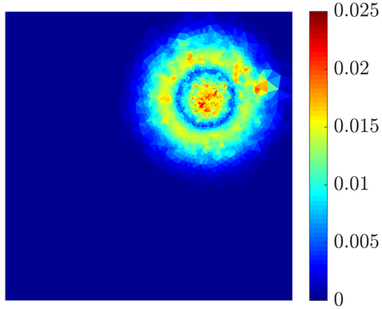}}
	\caption{Exact error and indicator map in the $\eltwo(\Omega_e)$ norm for the three stages of the adaptive process corresponding to the meshes and approximations shown in Figures~\ref{fig:Poisson_Adap_Mesh} and \ref{fig:Poisson_Adap_Sol}.}
	\label{fig:Poisson_Adap_Err}
\end{figure}

As described in Section~\ref{sc:adaptivity}, the adaptive process is driven by an error indicator that results from computing the relative difference in each cell, measured in the $\eltwo(\Omega_e)$ norm. To illustrate the efficiency of the proposed error indicator, Figure~\ref{fig:Poisson_Adap_Err} shows both the error indicator and the exact error for the three iterations, $\nit$, of the adaptive process corresponding to the meshes and approximations shown in Figures~\ref{fig:Poisson_Adap_Mesh} and \ref{fig:Poisson_Adap_Sol}.

\begin{figure}[!tb]
	\centering
	\subfigure[Error indicator]{\includegraphics[width=0.48\textwidth]{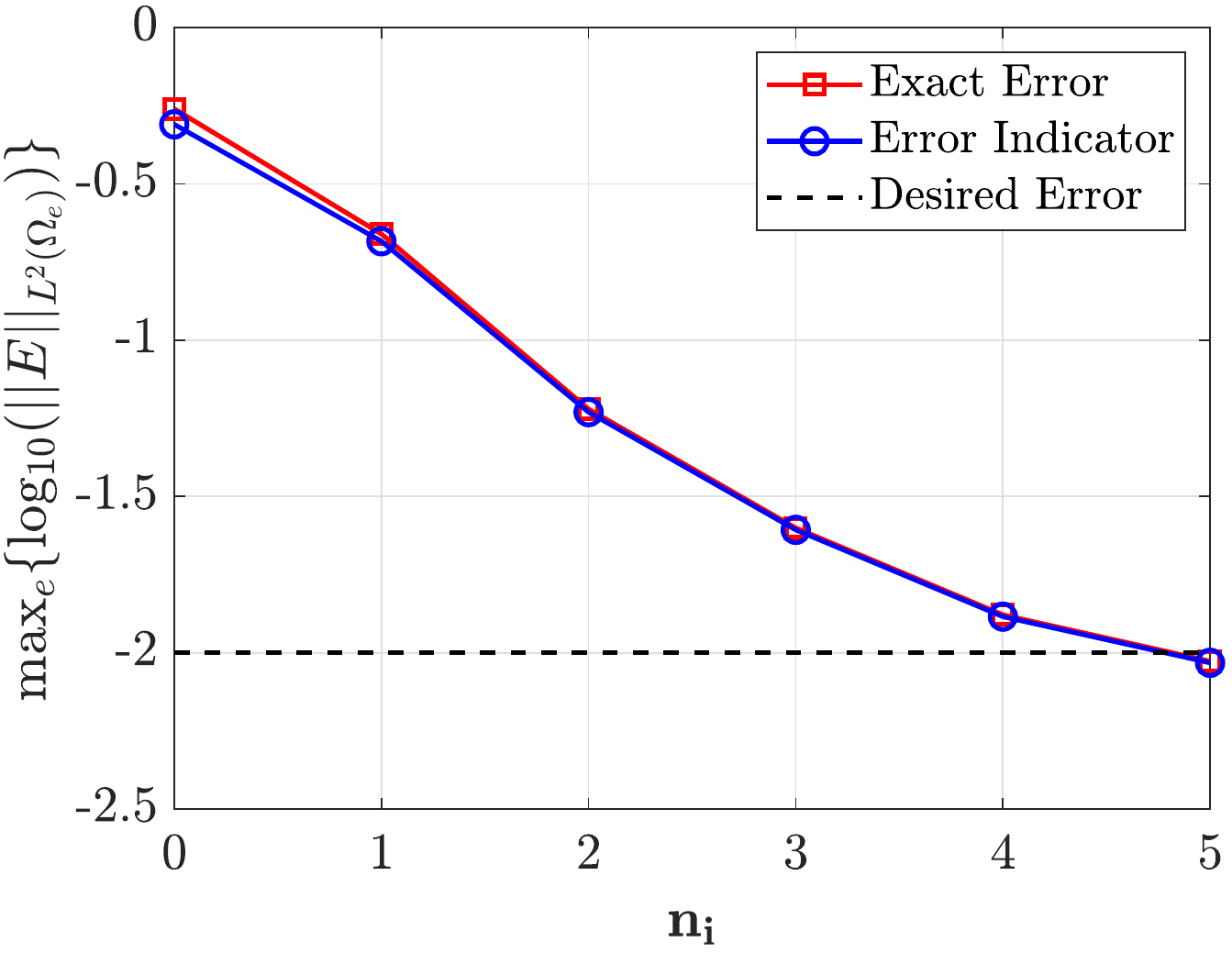}}
	\subfigure[Indicator efficiency]{\includegraphics[width=0.48\textwidth]{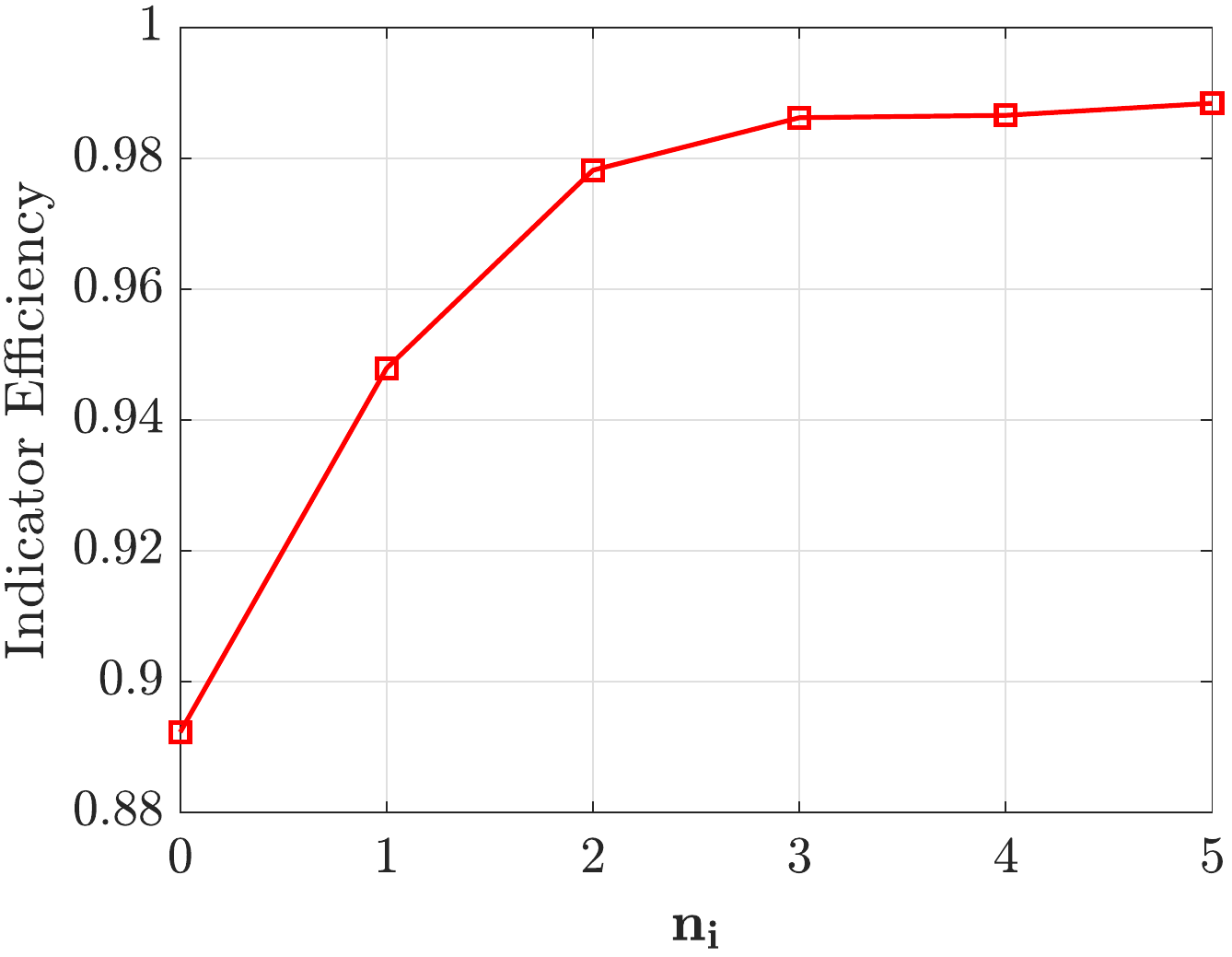}}
	\caption{(a) Maximum value of the indicator and the exact error over all the cells and (b) indicator efficiency.}
	\label{fig:Poisson_Adap_CovEff}
\end{figure}

To further illustrate the performance of the error indicator and the automatic adaptive process, Figure~\ref{fig:Poisson_Adap_CovEff}(a) shows the maximum value of the indicator and the exact error over all the cells as a function of the number of iterations of the mesh adaptive process.
The efficiency of the error indicator, defined as the ratio between the exact error and the indicator is also displayed in Figure~\ref{fig:Poisson_Adap_CovEff}(b), clearly illustrating the suitability of the proposed technique to drive an automatic mesh adaptive process.

\subsection{Mesh adaptivity for the Stokes problem} \label{sc:adaptivityStokes}

\begin{figure}[!tb]
	\centering
	\includegraphics[width=0.7\textwidth]{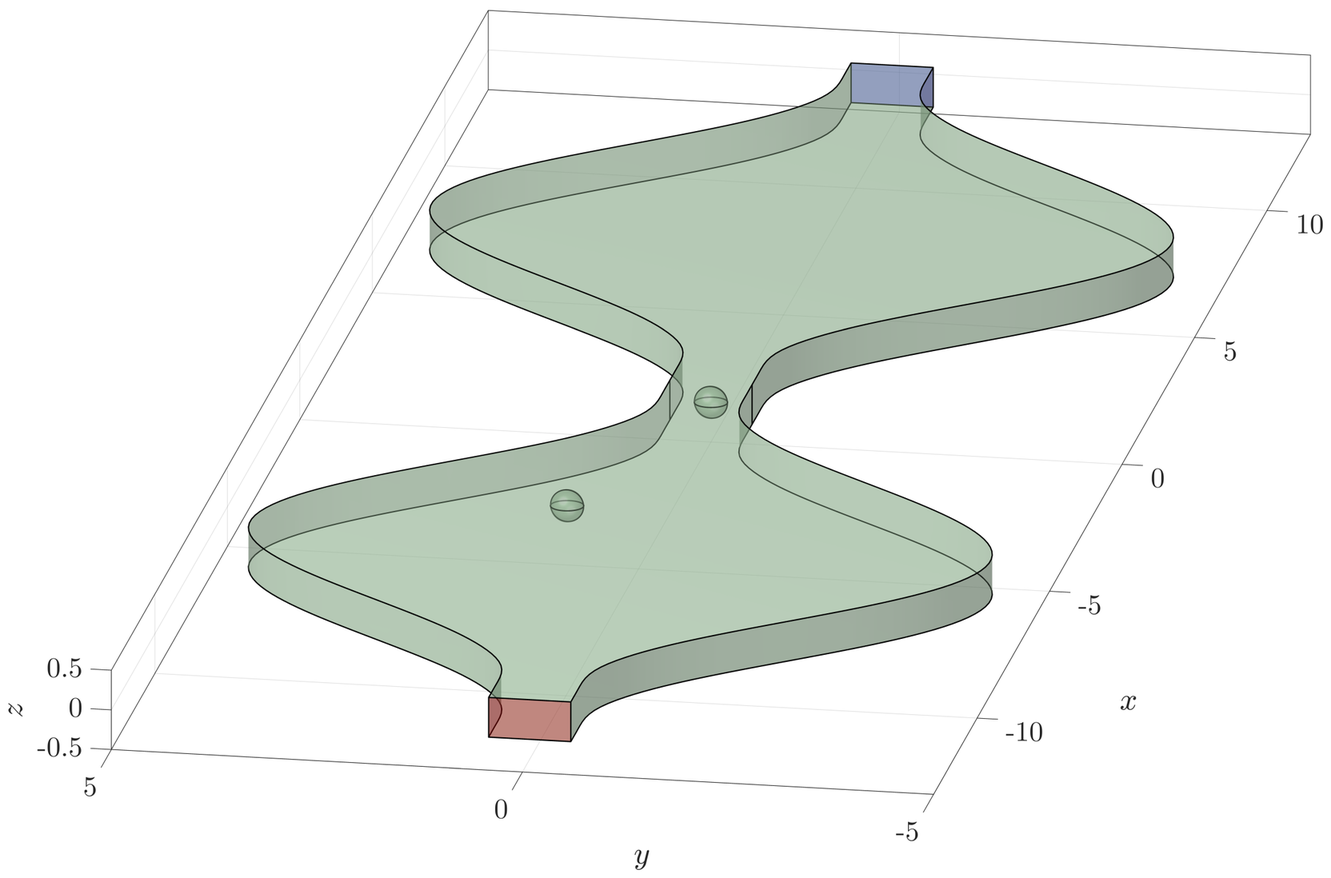}
	\caption{Geometry description for the computation of the Stokes flow in a corrugated channel with two spheres.}
	\label{fig:coarrugatedChannelGeometry}
\end{figure}

The last example involves the solution of the Stokes equations in three dimensions for the complex geometry, taken from~\cite{yang2017hydrodynamic}, depicted in Figure~\ref{fig:coarrugatedChannelGeometry}. 
The corrugated channel has a height of 0.5$\mu$m and the curved profile is obtained by repeating the expression
\begin{equation}
y = 
\begin{cases}
\frac{1}{2} \left( f_\omega + f_n \right) + \frac{1}{2} \left( f_\omega - f_n \right) \cos \left( \frac{8 \pi (2x-L)}{7L} \right) & \text{if} \quad \abs{x} < \frac{15}{16}L, \\
f_n & \text{if} \quad \frac{15}{16}L \leq \abs{x} \leq L,
\end{cases}
\end{equation}
where $f_n = 0.5\mu$m, $f_\omega = 4.5\mu$m and $L=12.5\mu$m.
Two spheres of radius 0.2$\mu$m are placed inside the channel. The first sphere, with centre $(0,0,0.25)\mu$m, is placed in the middle of the channel, where the cross section is minimum, and it is expected to produce a major flow disturbance. The second sphere, with centre at $(-3.75,1,0.25)\mu$m, is situated in a region with larger cross section and it is expected to disturb much less the flow. This scenario is utilised to show the ability of the proposed method to drive the adaptivity for a problem involving an incompressible flow in a  complex geometry.

A Dirichlet boundary condition, corresponding to a velocity inlet given by $\bu_D(x,y,z) = 64(y^2-1/4)(z^2-1/16)$, is introduced at one end of the channel, at $x=-L$, depicted in red in Figure~\ref{fig:coarrugatedChannelGeometry}.  A homogeneous Neumann boundary condition is imposed at the outlet, at $x=L$, depicted in blue in Figure~\ref{fig:coarrugatedChannelGeometry}. Homogeneous Dirichlet boundary conditions are imposed on the rest of the boundary, corresponding to material walls. 

\begin{figure}[!tb]
	\centering
	\subfigure[Mesh 1]{\includegraphics[width=0.32\textwidth]{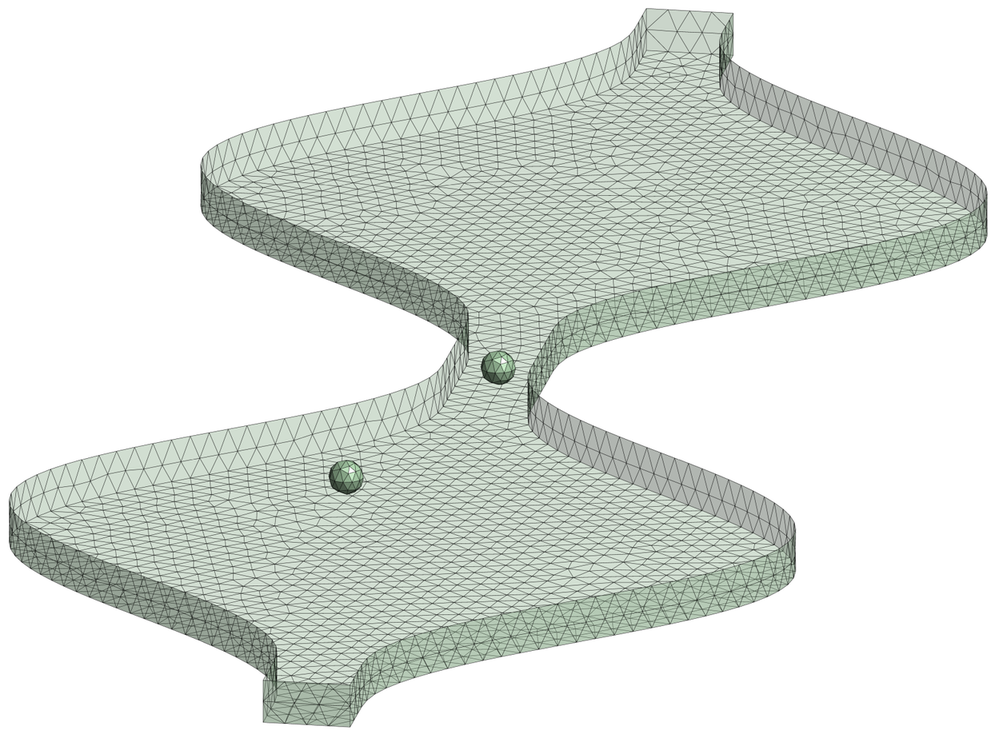}}
	\subfigure[Mesh 4]{\includegraphics[width=0.32\textwidth]{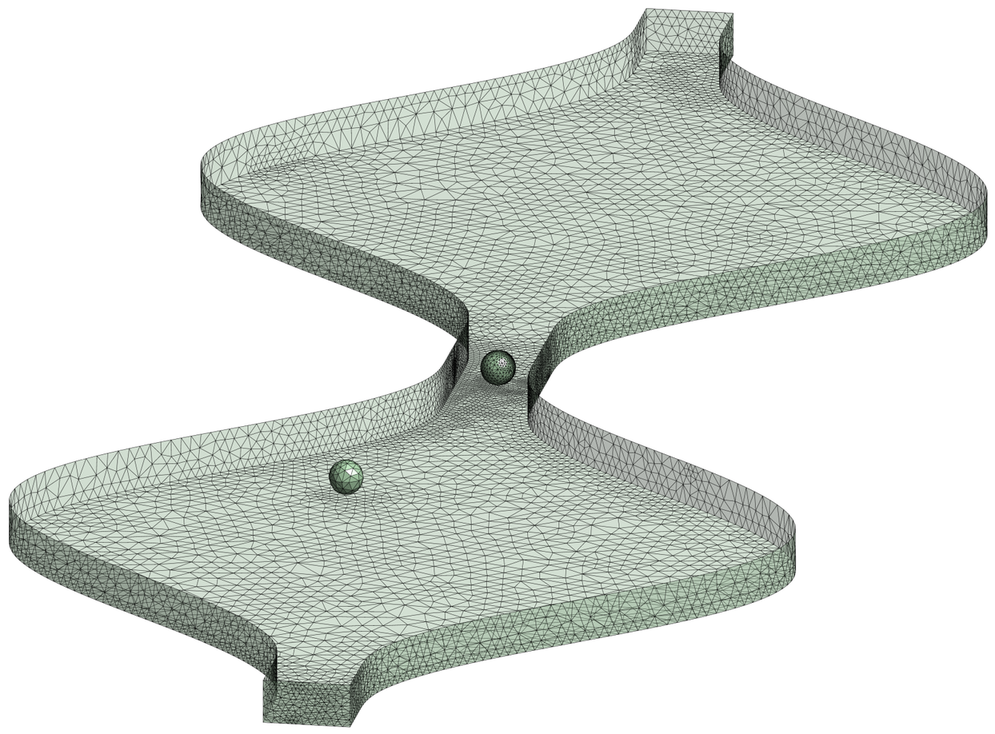}}
	\subfigure[Mesh 7]{\includegraphics[width=0.32\textwidth]{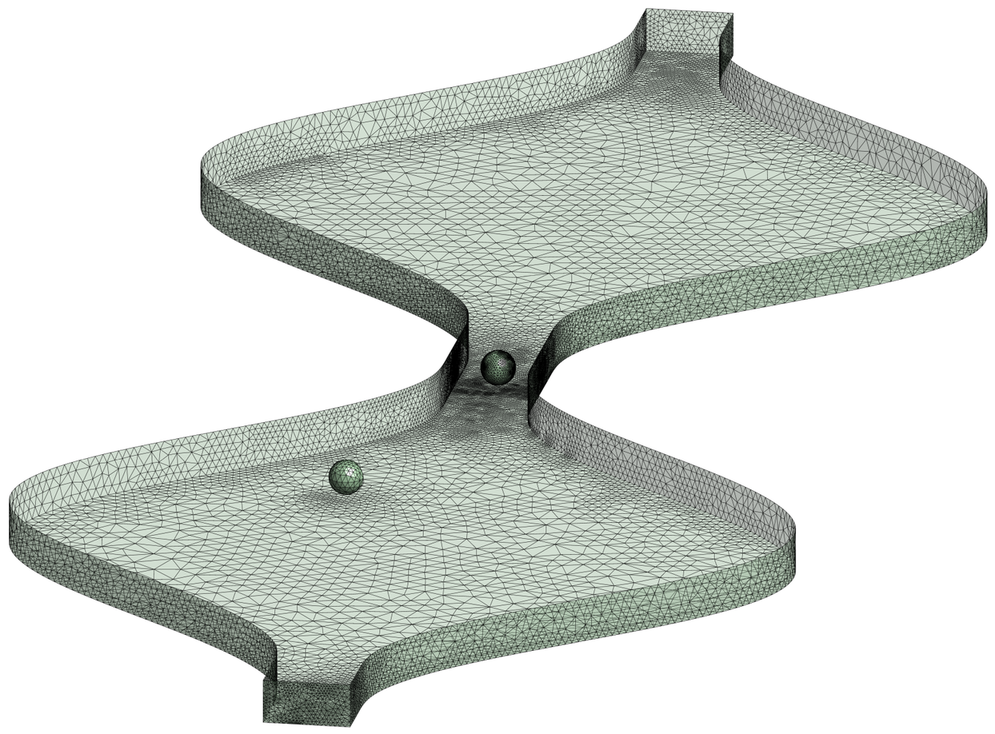}}
	\caption{Three of the meshes used for the adaptive computation of the Stokes flow in a corrugated channel with two spheres.}
	\label{fig:coarrugatedChannelMesh}
\end{figure}

The initial mesh, shown in Figure~\ref{fig:coarrugatedChannelMesh} (a), has 37,415 tetrahedral cells. The initial mesh is generated with a required element size of $0.05\mu$m on the surfaces defining two spheres, to ensure an appropriate geometric representation. An element size of $0.5\mu$m is imposed in the domain, with a smooth transition between these two values. 

\begin{figure}[!tb]
	\centering
	\subfigure[Mesh 1]{\includegraphics[width=0.32\textwidth]{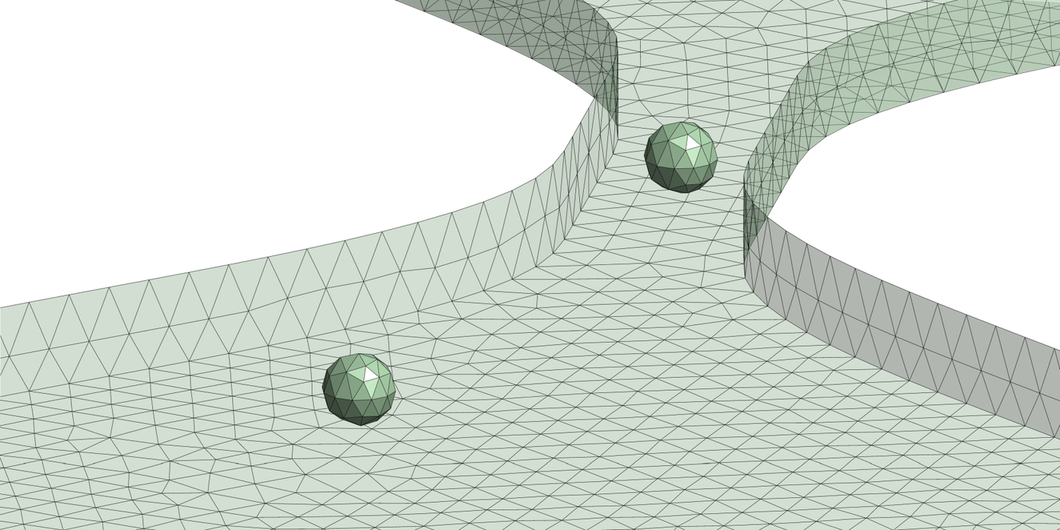}}
	\subfigure[Mesh 4]{\includegraphics[width=0.32\textwidth]{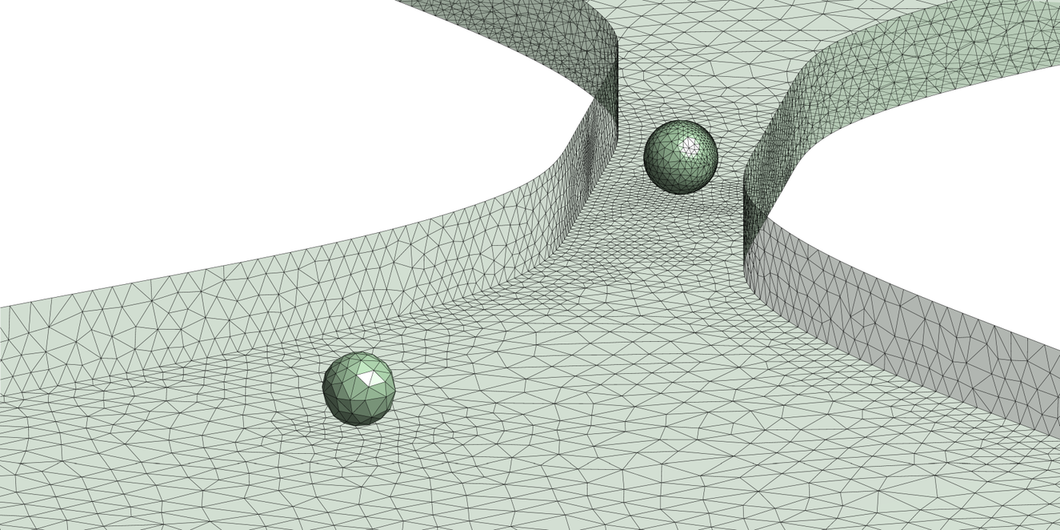}}
	\subfigure[Mesh 7]{\includegraphics[width=0.32\textwidth]{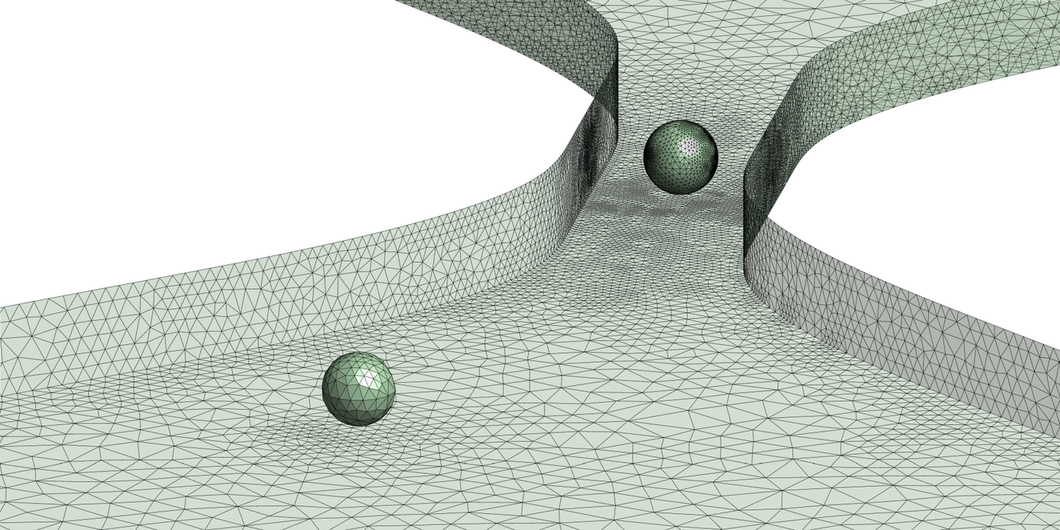}}
	\caption{Detailed view of three of the meshes shown in Figure~\ref{fig:coarrugatedChannelMesh}.}
	\label{fig:coarrugatedChannelMeshZoom}
\end{figure}

Next, the automatic adaptive process described in Section~\ref{sc:adaptivity} is applied, with a desired accuracy of $\varepsilon = 5 \times 10^{-2}$. Convergence of the estimated error to the desired tolerance is achieved in this example after seven iterations of the adaptive process. The resulting meshes in the fourth and seven iteration are depicted in Figures~\ref{fig:coarrugatedChannelMesh} (b) and (c) respectively. The fourth mesh has 61,871 cells, with a minimum and maximum element size of $0.010\mu$m and $0.594\mu$m respectively. The last mesh has 116,913 cells, with a minimum and maximum element size of $0.003\mu$m and $0.520\mu$m respectively. A detailed view of the three meshes of Figure~\ref{fig:coarrugatedChannelMesh} near the two spheres is shown in Figure~\ref{fig:coarrugatedChannelMeshZoom}.

\begin{figure}[!tb]
	\centering
	\subfigure[Velocity]{\includegraphics[width=0.48\textwidth]{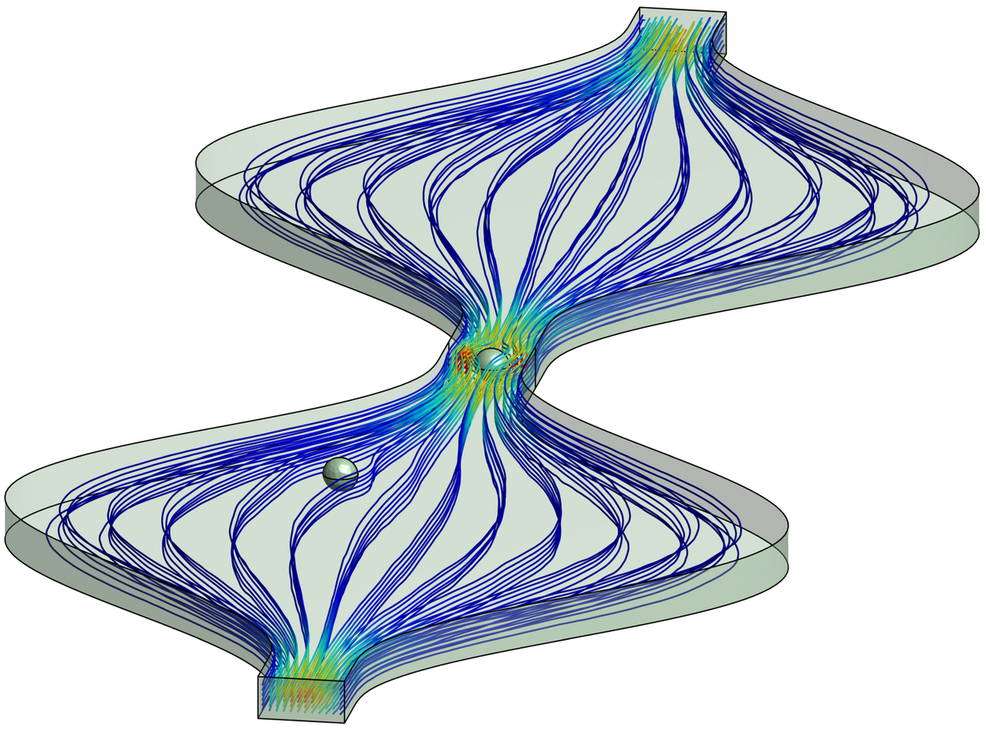}}
	\subfigure[Pressure]{\includegraphics[width=0.48\textwidth]{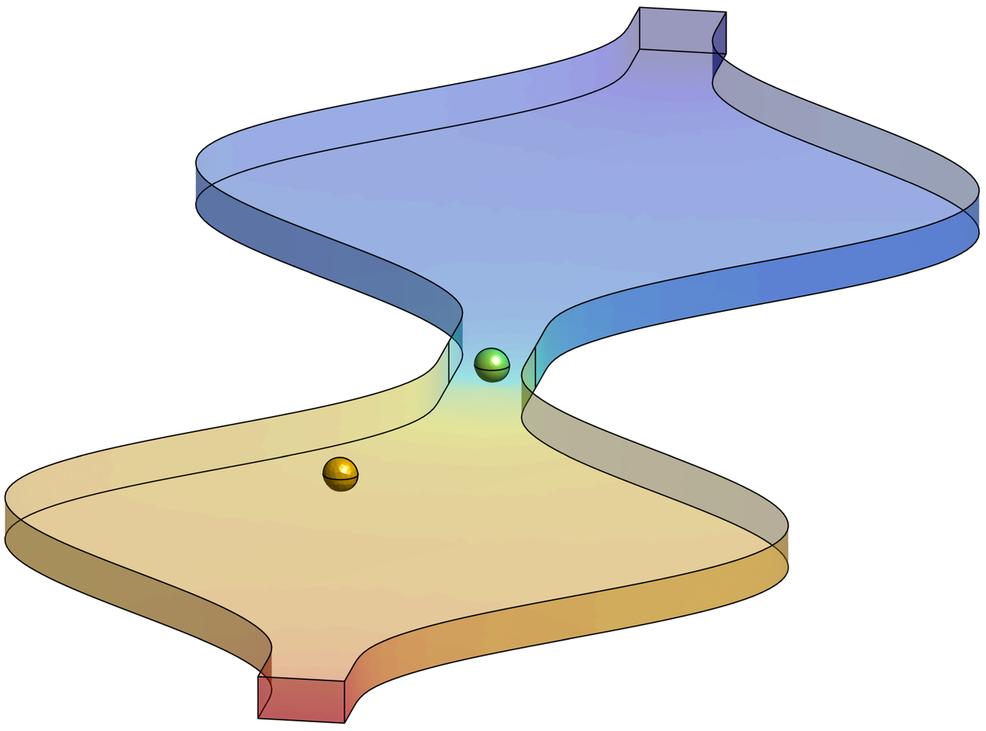}}	
	\caption{(a) Velocity streamlines and (b) pressure field, obtained in the mesh of Figure~\ref{fig:coarrugatedChannelMesh} (c), at the last iteration of the adaptive computation.}
	\label{fig:coarrugatedChannelSolution}
\end{figure}

\begin{figure}[!tb]
	\centering
	\subfigure[Velocity]{\includegraphics[width=0.48\textwidth]{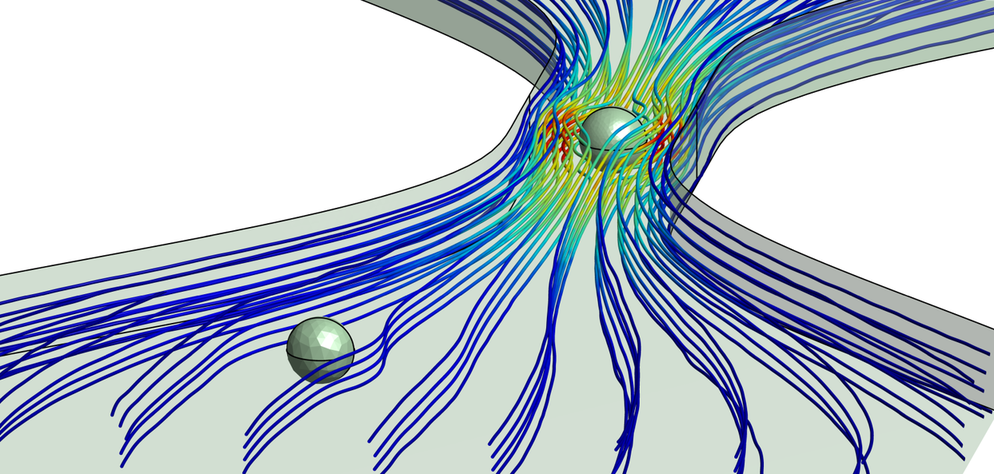}}
	\subfigure[Pressure]{\includegraphics[width=0.48\textwidth]{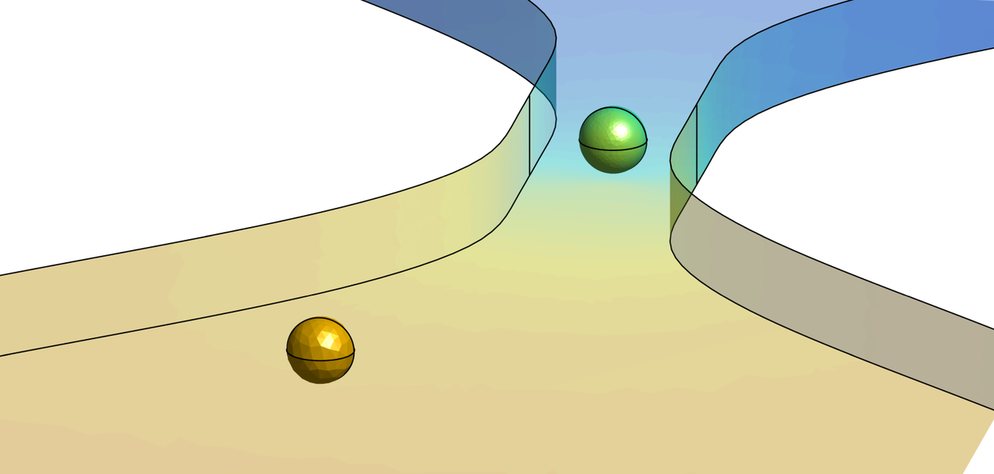}}	
	\caption{Detailed view of the (a) velocity streamlines and (b) pressure field, obtained in the mesh of Figure~\ref{fig:coarrugatedChannelMesh} (c), at the last iteration of the adaptive computation.}
	\label{fig:coarrugatedChannelSolutionZoom}
\end{figure}

The velocity streamlines and the pressure field obtained using the last mesh are represented in Figure~\ref{fig:coarrugatedChannelSolution}.
A detailed view of the velocity streamlines and the pressure field near the sphere surfaces is shown in Figure~\ref{fig:coarrugatedChannelSolutionZoom}. This Figure offers a qualitative comparison of the major disturbance caused by one sphere compared to the other. 
%

\section{Concluding remarks}
\label{sc:Conclusion}

This paper proposes a second-order FCFV method for the solution of scalar and vector elliptic problems in two and three dimensions. 

The proposed method preserves the attractive properties of the original first-order FCFV method, namely the first-order convergence of the gradient of the solution, without the need of a reconstruction, and the insensitivity to mesh distortion and stretching. It also satisfies the LBB condition in the context of incompressible flows. Contrary to the original FCFV, the proposed method guarantees second-order convergence of the solution. Numerical experiments show an increased performance when compared to the first-order method in terms of the CPU time required to achieve a desired accuracy as well as a lower sensitivity to the choice of the stabilisation parameter. A combination of first-order and second-order schemes is used to devise an error indicator that can be used to drive a mesh adaptivity process. An extensive set of numerical experiments has been used to demonstrate the optimal approximation properties of the method and more complex problems demonstrate its potential for large scale three dimensional simulations, including a Stokes problem where an automatic mesh adaptive process is employed.

\section*{Acknowledgements}

This work was partially supported by the European Union's Horizon 2020 research and innovation programme under the Marie Sk\l odowska-Curie Actions (Grant numbers: 675919 and 764636). The first author gratefully acknowledges the financial support of the Erasmus+ programme.
The second, third and last author were also supported by the Spanish Ministry of Economy and Competitiveness (Grant agreement No. DPI2017-85139-C2-2-R). 
The second and last authors are grateful for the financial support provided by the Generalitat de Catalunya (Grant agreement No. 2017-SGR-1278).

\bibliographystyle{abbrv}
\bibliography{Ref-HDG}

\end{document}